\documentclass[bj,preprint]{imsart}

\RequirePackage[OT1]{fontenc}
\RequirePackage{amsthm,amsmath,natbib}
\RequirePackage[colorlinks,citecolor=blue,urlcolor=blue]{hyperref}

\usepackage{color}  
\definecolor{AliceBlue}{rgb}{0.94,0.97,1.00}
\definecolor{AntiqueWhite1}{rgb}{1.00,0.94,0.86}
\definecolor{AntiqueWhite2}{rgb}{0.93,0.87,0.80}
\definecolor{AntiqueWhite3}{rgb}{0.80,0.75,0.69}
\definecolor{AntiqueWhite4}{rgb}{0.55,0.51,0.47}
\definecolor{AntiqueWhite}{rgb}{0.98,0.92,0.84}
\definecolor{BlanchedAlmond}{rgb}{1.00,0.92,0.80}
\definecolor{BlueViolet}{rgb}{0.54,0.17,0.89}
\definecolor{CadetBlue1}{rgb}{0.60,0.96,1.00}
\definecolor{CadetBlue2}{rgb}{0.56,0.90,0.93}
\definecolor{CadetBlue3}{rgb}{0.48,0.77,0.80}
\definecolor{CadetBlue4}{rgb}{0.33,0.53,0.55}
\definecolor{CadetBlue}{rgb}{0.37,0.62,0.63}
\definecolor{CornflowerBlue}{rgb}{0.39,0.58,0.93}
\definecolor{DarkBlue}{rgb}{0.00,0.00,0.55}
\definecolor{DarkCyan}{rgb}{0.00,0.55,0.55}
\definecolor{DarkGoldenrod1}{rgb}{1.00,0.73,0.06}
\definecolor{DarkGoldenrod2}{rgb}{0.93,0.68,0.05}
\definecolor{DarkGoldenrod3}{rgb}{0.80,0.58,0.05}
\definecolor{DarkGoldenrod4}{rgb}{0.55,0.40,0.03}
\definecolor{DarkGoldenrod}{rgb}{0.72,0.53,0.04}
\definecolor{DarkGray}{rgb}{0.66,0.66,0.66}
\definecolor{DarkGreen}{rgb}{0.00,0.39,0.00}
\definecolor{DarkGrey}{rgb}{0.66,0.66,0.66}
\definecolor{DarkKhaki}{rgb}{0.74,0.72,0.42}
\definecolor{DarkMagenta}{rgb}{0.55,0.00,0.55}
\definecolor{DarkOliveGreen1}{rgb}{0.79,1.00,0.44}
\definecolor{DarkOliveGreen2}{rgb}{0.74,0.93,0.41}
\definecolor{DarkOliveGreen3}{rgb}{0.64,0.80,0.35}
\definecolor{DarkOliveGreen4}{rgb}{0.43,0.55,0.24}
\definecolor{DarkOliveGreen}{rgb}{0.33,0.42,0.18}
\definecolor{DarkOrange1}{rgb}{1.00,0.50,0.00}
\definecolor{DarkOrange2}{rgb}{0.93,0.46,0.00}
\definecolor{DarkOrange3}{rgb}{0.80,0.40,0.00}
\definecolor{DarkOrange4}{rgb}{0.55,0.27,0.00}
\definecolor{DarkOrange}{rgb}{1.00,0.55,0.00}
\definecolor{DarkOrchid1}{rgb}{0.75,0.24,1.00}
\definecolor{DarkOrchid2}{rgb}{0.70,0.23,0.93}
\definecolor{DarkOrchid3}{rgb}{0.60,0.20,0.80}
\definecolor{DarkOrchid4}{rgb}{0.41,0.13,0.55}
\definecolor{DarkOrchid}{rgb}{0.60,0.20,0.80}
\definecolor{DarkRed}{rgb}{0.55,0.00,0.00}
\definecolor{DarkSalmon}{rgb}{0.91,0.59,0.48}
\definecolor{DarkSeaGreen1}{rgb}{0.76,1.00,0.76}
\definecolor{DarkSeaGreen2}{rgb}{0.71,0.93,0.71}
\definecolor{DarkSeaGreen3}{rgb}{0.61,0.80,0.61}
\definecolor{DarkSeaGreen4}{rgb}{0.41,0.55,0.41}
\definecolor{DarkSeaGreen}{rgb}{0.56,0.74,0.56}
\definecolor{DarkSlateBlue}{rgb}{0.28,0.24,0.55}
\definecolor{DarkSlateGray1}{rgb}{0.59,1.00,1.00}
\definecolor{DarkSlateGray2}{rgb}{0.55,0.93,0.93}
\definecolor{DarkSlateGray3}{rgb}{0.47,0.80,0.80}
\definecolor{DarkSlateGray4}{rgb}{0.32,0.55,0.55}
\definecolor{DarkSlateGray}{rgb}{0.18,0.31,0.31}
\definecolor{DarkSlateGrey}{rgb}{0.18,0.31,0.31}
\definecolor{DarkTurquoise}{rgb}{0.00,0.81,0.82}
\definecolor{DarkViolet}{rgb}{0.58,0.00,0.83}
\definecolor{DeepPink1}{rgb}{1.00,0.08,0.58}
\definecolor{DeepPink2}{rgb}{0.93,0.07,0.54}
\definecolor{DeepPink3}{rgb}{0.80,0.06,0.46}
\definecolor{DeepPink4}{rgb}{0.55,0.04,0.31}
\definecolor{DeepPink}{rgb}{1.00,0.08,0.58}
\definecolor{DeepSkyBlue1}{rgb}{0.00,0.75,1.00}
\definecolor{DeepSkyBlue2}{rgb}{0.00,0.70,0.93}
\definecolor{DeepSkyBlue3}{rgb}{0.00,0.60,0.80}
\definecolor{DeepSkyBlue4}{rgb}{0.00,0.41,0.55}
\definecolor{DeepSkyBlue}{rgb}{0.00,0.75,1.00}
\definecolor{DimGray}{rgb}{0.41,0.41,0.41}
\definecolor{DimGrey}{rgb}{0.41,0.41,0.41}
\definecolor{DodgerBlue1}{rgb}{0.12,0.56,1.00}
\definecolor{DodgerBlue2}{rgb}{0.11,0.53,0.93}
\definecolor{DodgerBlue3}{rgb}{0.09,0.45,0.80}
\definecolor{DodgerBlue4}{rgb}{0.06,0.31,0.55}
\definecolor{DodgerBlue}{rgb}{0.12,0.56,1.00}
\definecolor{FloralWhite}{rgb}{1.00,0.98,0.94}
\definecolor{ForestGreen}{rgb}{0.13,0.55,0.13}
\definecolor{GhostWhite}{rgb}{0.97,0.97,1.00}
\definecolor{GreenYellow}{rgb}{0.68,1.00,0.18}
\definecolor{HotPink1}{rgb}{1.00,0.43,0.71}
\definecolor{HotPink2}{rgb}{0.93,0.42,0.65}
\definecolor{HotPink3}{rgb}{0.80,0.38,0.56}
\definecolor{HotPink4}{rgb}{0.55,0.23,0.38}
\definecolor{HotPink}{rgb}{1.00,0.41,0.71}
\definecolor{IndianRed1}{rgb}{1.00,0.42,0.42}
\definecolor{IndianRed2}{rgb}{0.93,0.39,0.39}
\definecolor{IndianRed3}{rgb}{0.80,0.33,0.33}
\definecolor{IndianRed4}{rgb}{0.55,0.23,0.23}
\definecolor{IndianRed}{rgb}{0.80,0.36,0.36}
\definecolor{LavenderBlush1}{rgb}{1.00,0.94,0.96}
\definecolor{LavenderBlush2}{rgb}{0.93,0.88,0.90}
\definecolor{LavenderBlush3}{rgb}{0.80,0.76,0.77}
\definecolor{LavenderBlush4}{rgb}{0.55,0.51,0.53}
\definecolor{LavenderBlush}{rgb}{1.00,0.94,0.96}
\definecolor{LawnGreen}{rgb}{0.49,0.99,0.00}
\definecolor{LemonChiffon1}{rgb}{1.00,0.98,0.80}
\definecolor{LemonChiffon2}{rgb}{0.93,0.91,0.75}
\definecolor{LemonChiffon3}{rgb}{0.80,0.79,0.65}
\definecolor{LemonChiffon4}{rgb}{0.55,0.54,0.44}
\definecolor{LemonChiffon}{rgb}{1.00,0.98,0.80}
\definecolor{LightBlue1}{rgb}{0.75,0.94,1.00}
\definecolor{LightBlue2}{rgb}{0.70,0.87,0.93}
\definecolor{LightBlue3}{rgb}{0.60,0.75,0.80}
\definecolor{LightBlue4}{rgb}{0.41,0.51,0.55}
\definecolor{LightBlue}{rgb}{0.68,0.85,0.90}
\definecolor{LightCoral}{rgb}{0.94,0.50,0.50}
\definecolor{LightCyan1}{rgb}{0.88,1.00,1.00}
\definecolor{LightCyan2}{rgb}{0.82,0.93,0.93}
\definecolor{LightCyan3}{rgb}{0.71,0.80,0.80}
\definecolor{LightCyan4}{rgb}{0.48,0.55,0.55}
\definecolor{LightCyan}{rgb}{0.88,1.00,1.00}
\definecolor{LightGoldenrod1}{rgb}{1.00,0.93,0.55}
\definecolor{LightGoldenrod2}{rgb}{0.93,0.86,0.51}
\definecolor{LightGoldenrod3}{rgb}{0.80,0.75,0.44}
\definecolor{LightGoldenrod4}{rgb}{0.55,0.51,0.30}
\definecolor{LightGoldenrodYellow}{rgb}{0.98,0.98,0.82}
\definecolor{LightGoldenrod}{rgb}{0.93,0.87,0.51}
\definecolor{LightGray}{rgb}{0.83,0.83,0.83}
\definecolor{LightGreen}{rgb}{0.56,0.93,0.56}
\definecolor{LightGrey}{rgb}{0.83,0.83,0.83}
\definecolor{LightPink1}{rgb}{1.00,0.68,0.73}
\definecolor{LightPink2}{rgb}{0.93,0.64,0.68}
\definecolor{LightPink3}{rgb}{0.80,0.55,0.58}
\definecolor{LightPink4}{rgb}{0.55,0.37,0.40}
\definecolor{LightPink}{rgb}{1.00,0.71,0.76}
\definecolor{LightSalmon1}{rgb}{1.00,0.63,0.48}
\definecolor{LightSalmon2}{rgb}{0.93,0.58,0.45}
\definecolor{LightSalmon3}{rgb}{0.80,0.51,0.38}
\definecolor{LightSalmon4}{rgb}{0.55,0.34,0.26}
\definecolor{LightSalmon}{rgb}{1.00,0.63,0.48}
\definecolor{LightSeaGreen}{rgb}{0.13,0.70,0.67}
\definecolor{LightSkyBlue1}{rgb}{0.69,0.89,1.00}
\definecolor{LightSkyBlue2}{rgb}{0.64,0.83,0.93}
\definecolor{LightSkyBlue3}{rgb}{0.55,0.71,0.80}
\definecolor{LightSkyBlue4}{rgb}{0.38,0.48,0.55}
\definecolor{LightSkyBlue}{rgb}{0.53,0.81,0.98}
\definecolor{LightSlateBlue}{rgb}{0.52,0.44,1.00}
\definecolor{LightSlateGray}{rgb}{0.47,0.53,0.60}
\definecolor{LightSlateGrey}{rgb}{0.47,0.53,0.60}
\definecolor{LightSteelBlue1}{rgb}{0.79,0.88,1.00}
\definecolor{LightSteelBlue2}{rgb}{0.74,0.82,0.93}
\definecolor{LightSteelBlue3}{rgb}{0.64,0.71,0.80}
\definecolor{LightSteelBlue4}{rgb}{0.43,0.48,0.55}
\definecolor{LightSteelBlue}{rgb}{0.69,0.77,0.87}
\definecolor{LightYellow1}{rgb}{1.00,1.00,0.88}
\definecolor{LightYellow2}{rgb}{0.93,0.93,0.82}
\definecolor{LightYellow3}{rgb}{0.80,0.80,0.71}
\definecolor{LightYellow4}{rgb}{0.55,0.55,0.48}
\definecolor{LightYellow}{rgb}{1.00,1.00,0.88}
\definecolor{LimeGreen}{rgb}{0.20,0.80,0.20}
\definecolor{MediumAquamarine}{rgb}{0.40,0.80,0.67}
\definecolor{MediumBlue}{rgb}{0.00,0.00,0.80}
\definecolor{MediumOrchid1}{rgb}{0.88,0.40,1.00}
\definecolor{MediumOrchid2}{rgb}{0.82,0.37,0.93}
\definecolor{MediumOrchid3}{rgb}{0.71,0.32,0.80}
\definecolor{MediumOrchid4}{rgb}{0.48,0.22,0.55}
\definecolor{MediumOrchid}{rgb}{0.73,0.33,0.83}
\definecolor{MediumPurple1}{rgb}{0.67,0.51,1.00}
\definecolor{MediumPurple2}{rgb}{0.62,0.47,0.93}
\definecolor{MediumPurple3}{rgb}{0.54,0.41,0.80}
\definecolor{MediumPurple4}{rgb}{0.36,0.28,0.55}
\definecolor{MediumPurple}{rgb}{0.58,0.44,0.86}
\definecolor{MediumSeaGreen}{rgb}{0.24,0.70,0.44}
\definecolor{MediumSlateBlue}{rgb}{0.48,0.41,0.93}
\definecolor{MediumSpringGreen}{rgb}{0.00,0.98,0.60}
\definecolor{MediumTurquoise}{rgb}{0.28,0.82,0.80}
\definecolor{MediumVioletRed}{rgb}{0.78,0.08,0.52}
\definecolor{MidnightBlue}{rgb}{0.10,0.10,0.44}
\definecolor{MintCream}{rgb}{0.96,1.00,0.98}
\definecolor{MistyRose1}{rgb}{1.00,0.89,0.88}
\definecolor{MistyRose2}{rgb}{0.93,0.84,0.82}
\definecolor{MistyRose3}{rgb}{0.80,0.72,0.71}
\definecolor{MistyRose4}{rgb}{0.55,0.49,0.48}
\definecolor{MistyRose}{rgb}{1.00,0.89,0.88}
\definecolor{NavajoWhite1}{rgb}{1.00,0.87,0.68}
\definecolor{NavajoWhite2}{rgb}{0.93,0.81,0.63}
\definecolor{NavajoWhite3}{rgb}{0.80,0.70,0.55}
\definecolor{NavajoWhite4}{rgb}{0.55,0.47,0.37}
\definecolor{NavajoWhite}{rgb}{1.00,0.87,0.68}
\definecolor{NavyBlue}{rgb}{0.00,0.00,0.50}
\definecolor{OldLace}{rgb}{0.99,0.96,0.90}
\definecolor{OliveDrab1}{rgb}{0.75,1.00,0.24}
\definecolor{OliveDrab2}{rgb}{0.70,0.93,0.23}
\definecolor{OliveDrab3}{rgb}{0.60,0.80,0.20}
\definecolor{OliveDrab4}{rgb}{0.41,0.55,0.13}
\definecolor{OliveDrab}{rgb}{0.42,0.56,0.14}
\definecolor{OrangeRed1}{rgb}{1.00,0.27,0.00}
\definecolor{OrangeRed2}{rgb}{0.93,0.25,0.00}
\definecolor{OrangeRed3}{rgb}{0.80,0.22,0.00}
\definecolor{OrangeRed4}{rgb}{0.55,0.15,0.00}
\definecolor{OrangeRed}{rgb}{1.00,0.27,0.00}
\definecolor{PaleGoldenrod}{rgb}{0.93,0.91,0.67}
\definecolor{PaleGreen1}{rgb}{0.60,1.00,0.60}
\definecolor{PaleGreen2}{rgb}{0.56,0.93,0.56}
\definecolor{PaleGreen3}{rgb}{0.49,0.80,0.49}
\definecolor{PaleGreen4}{rgb}{0.33,0.55,0.33}
\definecolor{PaleGreen}{rgb}{0.60,0.98,0.60}
\definecolor{PaleTurquoise1}{rgb}{0.73,1.00,1.00}
\definecolor{PaleTurquoise2}{rgb}{0.68,0.93,0.93}
\definecolor{PaleTurquoise3}{rgb}{0.59,0.80,0.80}
\definecolor{PaleTurquoise4}{rgb}{0.40,0.55,0.55}
\definecolor{PaleTurquoise}{rgb}{0.69,0.93,0.93}
\definecolor{PaleVioletRed1}{rgb}{1.00,0.51,0.67}
\definecolor{PaleVioletRed2}{rgb}{0.93,0.47,0.62}
\definecolor{PaleVioletRed3}{rgb}{0.80,0.41,0.54}
\definecolor{PaleVioletRed4}{rgb}{0.55,0.28,0.36}
\definecolor{PaleVioletRed}{rgb}{0.86,0.44,0.58}
\definecolor{PapayaWhip}{rgb}{1.00,0.94,0.84}
\definecolor{PeachPuff1}{rgb}{1.00,0.85,0.73}
\definecolor{PeachPuff2}{rgb}{0.93,0.80,0.68}
\definecolor{PeachPuff3}{rgb}{0.80,0.69,0.58}
\definecolor{PeachPuff4}{rgb}{0.55,0.47,0.40}
\definecolor{PeachPuff}{rgb}{1.00,0.85,0.73}
\definecolor{PowderBlue}{rgb}{0.69,0.88,0.90}
\definecolor{RosyBrown1}{rgb}{1.00,0.76,0.76}
\definecolor{RosyBrown2}{rgb}{0.93,0.71,0.71}
\definecolor{RosyBrown3}{rgb}{0.80,0.61,0.61}
\definecolor{RosyBrown4}{rgb}{0.55,0.41,0.41}
\definecolor{RosyBrown}{rgb}{0.74,0.56,0.56}
\definecolor{RoyalBlue1}{rgb}{0.28,0.46,1.00}
\definecolor{RoyalBlue2}{rgb}{0.26,0.43,0.93}
\definecolor{RoyalBlue3}{rgb}{0.23,0.37,0.80}
\definecolor{RoyalBlue4}{rgb}{0.15,0.25,0.55}
\definecolor{RoyalBlue}{rgb}{0.25,0.41,0.88}
\definecolor{SaddleBrown}{rgb}{0.55,0.27,0.07}
\definecolor{SandyBrown}{rgb}{0.96,0.64,0.38}
\definecolor{SeaGreen1}{rgb}{0.33,1.00,0.62}
\definecolor{SeaGreen2}{rgb}{0.31,0.93,0.58}
\definecolor{SeaGreen3}{rgb}{0.26,0.80,0.50}
\definecolor{SeaGreen4}{rgb}{0.18,0.55,0.34}
\definecolor{SeaGreen}{rgb}{0.18,0.55,0.34}
\definecolor{SkyBlue1}{rgb}{0.53,0.81,1.00}
\definecolor{SkyBlue2}{rgb}{0.49,0.75,0.93}
\definecolor{SkyBlue3}{rgb}{0.42,0.65,0.80}
\definecolor{SkyBlue4}{rgb}{0.29,0.44,0.55}
\definecolor{SkyBlue}{rgb}{0.53,0.81,0.92}
\definecolor{SlateBlue1}{rgb}{0.51,0.44,1.00}
\definecolor{SlateBlue2}{rgb}{0.48,0.40,0.93}
\definecolor{SlateBlue3}{rgb}{0.41,0.35,0.80}
\definecolor{SlateBlue4}{rgb}{0.28,0.24,0.55}
\definecolor{SlateBlue}{rgb}{0.42,0.35,0.80}
\definecolor{SlateGray1}{rgb}{0.78,0.89,1.00}
\definecolor{SlateGray2}{rgb}{0.73,0.83,0.93}
\definecolor{SlateGray3}{rgb}{0.62,0.71,0.80}
\definecolor{SlateGray4}{rgb}{0.42,0.48,0.55}
\definecolor{SlateGray}{rgb}{0.44,0.50,0.56}
\definecolor{SlateGrey}{rgb}{0.44,0.50,0.56}
\definecolor{SpringGreen1}{rgb}{0.00,1.00,0.50}
\definecolor{SpringGreen2}{rgb}{0.00,0.93,0.46}
\definecolor{SpringGreen3}{rgb}{0.00,0.80,0.40}
\definecolor{SpringGreen4}{rgb}{0.00,0.55,0.27}
\definecolor{SpringGreen}{rgb}{0.00,1.00,0.50}
\definecolor{SteelBlue1}{rgb}{0.39,0.72,1.00}
\definecolor{SteelBlue2}{rgb}{0.36,0.67,0.93}
\definecolor{SteelBlue3}{rgb}{0.31,0.58,0.80}
\definecolor{SteelBlue4}{rgb}{0.21,0.39,0.55}
\definecolor{SteelBlue}{rgb}{0.27,0.51,0.71}
\definecolor{VioletRed1}{rgb}{1.00,0.24,0.59}
\definecolor{VioletRed2}{rgb}{0.93,0.23,0.55}
\definecolor{VioletRed3}{rgb}{0.80,0.20,0.47}
\definecolor{VioletRed4}{rgb}{0.55,0.13,0.32}
\definecolor{VioletRed}{rgb}{0.82,0.13,0.56}
\definecolor{WhiteSmoke}{rgb}{0.96,0.96,0.96}
\definecolor{YellowGreen}{rgb}{0.60,0.80,0.20}
\definecolor{aliceblue}{rgb}{0.94,0.97,1.00}
\definecolor{antiquewhite}{rgb}{0.98,0.92,0.84}
\definecolor{aquamarine1}{rgb}{0.50,1.00,0.83}
\definecolor{aquamarine2}{rgb}{0.46,0.93,0.78}
\definecolor{aquamarine3}{rgb}{0.40,0.80,0.67}
\definecolor{aquamarine4}{rgb}{0.27,0.55,0.45}
\definecolor{aquamarine}{rgb}{0.50,1.00,0.83}
\definecolor{azure1}{rgb}{0.94,1.00,1.00}
\definecolor{azure2}{rgb}{0.88,0.93,0.93}
\definecolor{azure3}{rgb}{0.76,0.80,0.80}
\definecolor{azure4}{rgb}{0.51,0.55,0.55}
\definecolor{azure}{rgb}{0.94,1.00,1.00}
\definecolor{beige}{rgb}{0.96,0.96,0.86}
\definecolor{bisque1}{rgb}{1.00,0.89,0.77}
\definecolor{bisque2}{rgb}{0.93,0.84,0.72}
\definecolor{bisque3}{rgb}{0.80,0.72,0.62}
\definecolor{bisque4}{rgb}{0.55,0.49,0.42}
\definecolor{bisque}{rgb}{1.00,0.89,0.77}
\definecolor{black}{rgb}{0.00,0.00,0.00}
\definecolor{blanchedalmond}{rgb}{1.00,0.92,0.80}
\definecolor{blue1}{rgb}{0.00,0.00,1.00}
\definecolor{blue2}{rgb}{0.00,0.00,0.93}
\definecolor{blue3}{rgb}{0.00,0.00,0.80}
\definecolor{blue4}{rgb}{0.00,0.00,0.55}
\definecolor{blueviolet}{rgb}{0.54,0.17,0.89}
\definecolor{blue}{rgb}{0.00,0.00,1.00}
\definecolor{brown1}{rgb}{1.00,0.25,0.25}
\definecolor{brown2}{rgb}{0.93,0.23,0.23}
\definecolor{brown3}{rgb}{0.80,0.20,0.20}
\definecolor{brown4}{rgb}{0.55,0.14,0.14}
\definecolor{brown}{rgb}{0.65,0.16,0.16}
\definecolor{burlywood1}{rgb}{1.00,0.83,0.61}
\definecolor{burlywood2}{rgb}{0.93,0.77,0.57}
\definecolor{burlywood3}{rgb}{0.80,0.67,0.49}
\definecolor{burlywood4}{rgb}{0.55,0.45,0.33}
\definecolor{burlywood}{rgb}{0.87,0.72,0.53}
\definecolor{cadetblue}{rgb}{0.37,0.62,0.63}
\definecolor{chartreuse1}{rgb}{0.50,1.00,0.00}
\definecolor{chartreuse2}{rgb}{0.46,0.93,0.00}
\definecolor{chartreuse3}{rgb}{0.40,0.80,0.00}
\definecolor{chartreuse4}{rgb}{0.27,0.55,0.00}
\definecolor{chartreuse}{rgb}{0.50,1.00,0.00}
\definecolor{chocolate1}{rgb}{1.00,0.50,0.14}
\definecolor{chocolate2}{rgb}{0.93,0.46,0.13}
\definecolor{chocolate3}{rgb}{0.80,0.40,0.11}
\definecolor{chocolate4}{rgb}{0.55,0.27,0.07}
\definecolor{chocolate}{rgb}{0.82,0.41,0.12}
\definecolor{coral1}{rgb}{1.00,0.45,0.34}
\definecolor{coral2}{rgb}{0.93,0.42,0.31}
\definecolor{coral3}{rgb}{0.80,0.36,0.27}
\definecolor{coral4}{rgb}{0.55,0.24,0.18}
\definecolor{coral}{rgb}{1.00,0.50,0.31}
\definecolor{cornflowerblue}{rgb}{0.39,0.58,0.93}
\definecolor{cornsilk1}{rgb}{1.00,0.97,0.86}
\definecolor{cornsilk2}{rgb}{0.93,0.91,0.80}
\definecolor{cornsilk3}{rgb}{0.80,0.78,0.69}
\definecolor{cornsilk4}{rgb}{0.55,0.53,0.47}
\definecolor{cornsilk}{rgb}{1.00,0.97,0.86}
\definecolor{cyan1}{rgb}{0.00,1.00,1.00}
\definecolor{cyan2}{rgb}{0.00,0.93,0.93}
\definecolor{cyan3}{rgb}{0.00,0.80,0.80}
\definecolor{cyan4}{rgb}{0.00,0.55,0.55}
\definecolor{cyan}{rgb}{0.00,1.00,1.00}
\definecolor{darkblue}{rgb}{0.00,0.00,0.55}
\definecolor{darkcyan}{rgb}{0.00,0.55,0.55}
\definecolor{darkgoldenrod}{rgb}{0.72,0.53,0.04}
\definecolor{darkgray}{rgb}{0.66,0.66,0.66}
\definecolor{darkgreen}{rgb}{0.00,0.39,0.00}
\definecolor{darkgrey}{rgb}{0.66,0.66,0.66}
\definecolor{darkkhaki}{rgb}{0.74,0.72,0.42}
\definecolor{darkmagenta}{rgb}{0.55,0.00,0.55}
\definecolor{darkolive}{rgb}{0.33,0.42,0.18}
\definecolor{darkorange}{rgb}{1.00,0.55,0.00}
\definecolor{darkorchid}{rgb}{0.60,0.20,0.80}
\definecolor{darkred}{rgb}{0.55,0.00,0.00}
\definecolor{darksalmon}{rgb}{0.91,0.59,0.48}
\definecolor{darksea}{rgb}{0.56,0.74,0.56}
\definecolor{darkslate}{rgb}{0.18,0.31,0.31}
\definecolor{darkslate}{rgb}{0.18,0.31,0.31}
\definecolor{darkslate}{rgb}{0.28,0.24,0.55}
\definecolor{darkturquoise}{rgb}{0.00,0.81,0.82}
\definecolor{darkviolet}{rgb}{0.58,0.00,0.83}
\definecolor{deeppink}{rgb}{1.00,0.08,0.58}
\definecolor{deepsky}{rgb}{0.00,0.75,1.00}
\definecolor{dimgray}{rgb}{0.41,0.41,0.41}
\definecolor{dimgrey}{rgb}{0.41,0.41,0.41}
\definecolor{dodgerblue}{rgb}{0.12,0.56,1.00}
\definecolor{firebrick1}{rgb}{1.00,0.19,0.19}
\definecolor{firebrick2}{rgb}{0.93,0.17,0.17}
\definecolor{firebrick3}{rgb}{0.80,0.15,0.15}
\definecolor{firebrick4}{rgb}{0.55,0.10,0.10}
\definecolor{firebrick}{rgb}{0.70,0.13,0.13}
\definecolor{floralwhite}{rgb}{1.00,0.98,0.94}
\definecolor{forestgreen}{rgb}{0.13,0.55,0.13}
\definecolor{gainsboro}{rgb}{0.86,0.86,0.86}
\definecolor{ghostwhite}{rgb}{0.97,0.97,1.00}
\definecolor{gold1}{rgb}{1.00,0.84,0.00}
\definecolor{gold2}{rgb}{0.93,0.79,0.00}
\definecolor{gold3}{rgb}{0.80,0.68,0.00}
\definecolor{gold4}{rgb}{0.55,0.46,0.00}
\definecolor{goldenrod1}{rgb}{1.00,0.76,0.15}
\definecolor{goldenrod2}{rgb}{0.93,0.71,0.13}
\definecolor{goldenrod3}{rgb}{0.80,0.61,0.11}
\definecolor{goldenrod4}{rgb}{0.55,0.41,0.08}
\definecolor{goldenrod}{rgb}{0.85,0.65,0.13}
\definecolor{gold}{rgb}{1.00,0.84,0.00}
\definecolor{gray0}{rgb}{0.00,0.00,0.00}
\definecolor{gray100}{rgb}{1.00,1.00,1.00}
\definecolor{gray10}{rgb}{0.10,0.10,0.10}
\definecolor{gray11}{rgb}{0.11,0.11,0.11}
\definecolor{gray12}{rgb}{0.12,0.12,0.12}
\definecolor{gray13}{rgb}{0.13,0.13,0.13}
\definecolor{gray14}{rgb}{0.14,0.14,0.14}
\definecolor{gray15}{rgb}{0.15,0.15,0.15}
\definecolor{gray16}{rgb}{0.16,0.16,0.16}
\definecolor{gray17}{rgb}{0.17,0.17,0.17}
\definecolor{gray18}{rgb}{0.18,0.18,0.18}
\definecolor{gray19}{rgb}{0.19,0.19,0.19}
\definecolor{gray1}{rgb}{0.01,0.01,0.01}
\definecolor{gray20}{rgb}{0.20,0.20,0.20}
\definecolor{gray21}{rgb}{0.21,0.21,0.21}
\definecolor{gray22}{rgb}{0.22,0.22,0.22}
\definecolor{gray23}{rgb}{0.23,0.23,0.23}
\definecolor{gray24}{rgb}{0.24,0.24,0.24}
\definecolor{gray25}{rgb}{0.25,0.25,0.25}
\definecolor{gray26}{rgb}{0.26,0.26,0.26}
\definecolor{gray27}{rgb}{0.27,0.27,0.27}
\definecolor{gray28}{rgb}{0.28,0.28,0.28}
\definecolor{gray29}{rgb}{0.29,0.29,0.29}
\definecolor{gray2}{rgb}{0.02,0.02,0.02}
\definecolor{gray30}{rgb}{0.30,0.30,0.30}
\definecolor{gray31}{rgb}{0.31,0.31,0.31}
\definecolor{gray32}{rgb}{0.32,0.32,0.32}
\definecolor{gray33}{rgb}{0.33,0.33,0.33}
\definecolor{gray34}{rgb}{0.34,0.34,0.34}
\definecolor{gray35}{rgb}{0.35,0.35,0.35}
\definecolor{gray36}{rgb}{0.36,0.36,0.36}
\definecolor{gray37}{rgb}{0.37,0.37,0.37}
\definecolor{gray38}{rgb}{0.38,0.38,0.38}
\definecolor{gray39}{rgb}{0.39,0.39,0.39}
\definecolor{gray3}{rgb}{0.03,0.03,0.03}
\definecolor{gray40}{rgb}{0.40,0.40,0.40}
\definecolor{gray41}{rgb}{0.41,0.41,0.41}
\definecolor{gray42}{rgb}{0.42,0.42,0.42}
\definecolor{gray43}{rgb}{0.43,0.43,0.43}
\definecolor{gray44}{rgb}{0.44,0.44,0.44}
\definecolor{gray45}{rgb}{0.45,0.45,0.45}
\definecolor{gray46}{rgb}{0.46,0.46,0.46}
\definecolor{gray47}{rgb}{0.47,0.47,0.47}
\definecolor{gray48}{rgb}{0.48,0.48,0.48}
\definecolor{gray49}{rgb}{0.49,0.49,0.49}
\definecolor{gray4}{rgb}{0.04,0.04,0.04}
\definecolor{gray50}{rgb}{0.50,0.50,0.50}
\definecolor{gray51}{rgb}{0.51,0.51,0.51}
\definecolor{gray52}{rgb}{0.52,0.52,0.52}
\definecolor{gray53}{rgb}{0.53,0.53,0.53}
\definecolor{gray54}{rgb}{0.54,0.54,0.54}
\definecolor{gray55}{rgb}{0.55,0.55,0.55}
\definecolor{gray56}{rgb}{0.56,0.56,0.56}
\definecolor{gray57}{rgb}{0.57,0.57,0.57}
\definecolor{gray58}{rgb}{0.58,0.58,0.58}
\definecolor{gray59}{rgb}{0.59,0.59,0.59}
\definecolor{gray5}{rgb}{0.05,0.05,0.05}
\definecolor{gray60}{rgb}{0.60,0.60,0.60}
\definecolor{gray61}{rgb}{0.61,0.61,0.61}
\definecolor{gray62}{rgb}{0.62,0.62,0.62}
\definecolor{gray63}{rgb}{0.63,0.63,0.63}
\definecolor{gray64}{rgb}{0.64,0.64,0.64}
\definecolor{gray65}{rgb}{0.65,0.65,0.65}
\definecolor{gray66}{rgb}{0.66,0.66,0.66}
\definecolor{gray67}{rgb}{0.67,0.67,0.67}
\definecolor{gray68}{rgb}{0.68,0.68,0.68}
\definecolor{gray69}{rgb}{0.69,0.69,0.69}
\definecolor{gray6}{rgb}{0.06,0.06,0.06}
\definecolor{gray70}{rgb}{0.70,0.70,0.70}
\definecolor{gray71}{rgb}{0.71,0.71,0.71}
\definecolor{gray72}{rgb}{0.72,0.72,0.72}
\definecolor{gray73}{rgb}{0.73,0.73,0.73}
\definecolor{gray74}{rgb}{0.74,0.74,0.74}
\definecolor{gray75}{rgb}{0.75,0.75,0.75}
\definecolor{gray76}{rgb}{0.76,0.76,0.76}
\definecolor{gray77}{rgb}{0.77,0.77,0.77}
\definecolor{gray78}{rgb}{0.78,0.78,0.78}
\definecolor{gray79}{rgb}{0.79,0.79,0.79}
\definecolor{gray7}{rgb}{0.07,0.07,0.07}
\definecolor{gray80}{rgb}{0.80,0.80,0.80}
\definecolor{gray81}{rgb}{0.81,0.81,0.81}
\definecolor{gray82}{rgb}{0.82,0.82,0.82}
\definecolor{gray83}{rgb}{0.83,0.83,0.83}
\definecolor{gray84}{rgb}{0.84,0.84,0.84}
\definecolor{gray85}{rgb}{0.85,0.85,0.85}
\definecolor{gray86}{rgb}{0.86,0.86,0.86}
\definecolor{gray87}{rgb}{0.87,0.87,0.87}
\definecolor{gray88}{rgb}{0.88,0.88,0.88}
\definecolor{gray89}{rgb}{0.89,0.89,0.89}
\definecolor{gray8}{rgb}{0.08,0.08,0.08}
\definecolor{gray90}{rgb}{0.90,0.90,0.90}
\definecolor{gray91}{rgb}{0.91,0.91,0.91}
\definecolor{gray92}{rgb}{0.92,0.92,0.92}
\definecolor{gray93}{rgb}{0.93,0.93,0.93}
\definecolor{gray94}{rgb}{0.94,0.94,0.94}
\definecolor{gray95}{rgb}{0.95,0.95,0.95}
\definecolor{gray96}{rgb}{0.96,0.96,0.96}
\definecolor{gray97}{rgb}{0.97,0.97,0.97}
\definecolor{gray98}{rgb}{0.98,0.98,0.98}
\definecolor{gray99}{rgb}{0.99,0.99,0.99}
\definecolor{gray9}{rgb}{0.09,0.09,0.09}
\definecolor{gray}{rgb}{0.75,0.75,0.75}
\definecolor{green1}{rgb}{0.00,1.00,0.00}
\definecolor{green2}{rgb}{0.00,0.93,0.00}
\definecolor{green3}{rgb}{0.00,0.80,0.00}
\definecolor{green4}{rgb}{0.00,0.55,0.00}
\definecolor{greenyellow}{rgb}{0.68,1.00,0.18}
\definecolor{green}{rgb}{0.00,1.00,0.00}
\definecolor{grey0}{rgb}{0.00,0.00,0.00}
\definecolor{grey100}{rgb}{1.00,1.00,1.00}
\definecolor{grey10}{rgb}{0.10,0.10,0.10}
\definecolor{grey11}{rgb}{0.11,0.11,0.11}
\definecolor{grey12}{rgb}{0.12,0.12,0.12}
\definecolor{grey13}{rgb}{0.13,0.13,0.13}
\definecolor{grey14}{rgb}{0.14,0.14,0.14}
\definecolor{grey15}{rgb}{0.15,0.15,0.15}
\definecolor{grey16}{rgb}{0.16,0.16,0.16}
\definecolor{grey17}{rgb}{0.17,0.17,0.17}
\definecolor{grey18}{rgb}{0.18,0.18,0.18}
\definecolor{grey19}{rgb}{0.19,0.19,0.19}
\definecolor{grey1}{rgb}{0.01,0.01,0.01}
\definecolor{grey20}{rgb}{0.20,0.20,0.20}
\definecolor{grey21}{rgb}{0.21,0.21,0.21}
\definecolor{grey22}{rgb}{0.22,0.22,0.22}
\definecolor{grey23}{rgb}{0.23,0.23,0.23}
\definecolor{grey24}{rgb}{0.24,0.24,0.24}
\definecolor{grey25}{rgb}{0.25,0.25,0.25}
\definecolor{grey26}{rgb}{0.26,0.26,0.26}
\definecolor{grey27}{rgb}{0.27,0.27,0.27}
\definecolor{grey28}{rgb}{0.28,0.28,0.28}
\definecolor{grey29}{rgb}{0.29,0.29,0.29}
\definecolor{grey2}{rgb}{0.02,0.02,0.02}
\definecolor{grey30}{rgb}{0.30,0.30,0.30}
\definecolor{grey31}{rgb}{0.31,0.31,0.31}
\definecolor{grey32}{rgb}{0.32,0.32,0.32}
\definecolor{grey33}{rgb}{0.33,0.33,0.33}
\definecolor{grey34}{rgb}{0.34,0.34,0.34}
\definecolor{grey35}{rgb}{0.35,0.35,0.35}
\definecolor{grey36}{rgb}{0.36,0.36,0.36}
\definecolor{grey37}{rgb}{0.37,0.37,0.37}
\definecolor{grey38}{rgb}{0.38,0.38,0.38}
\definecolor{grey39}{rgb}{0.39,0.39,0.39}
\definecolor{grey3}{rgb}{0.03,0.03,0.03}
\definecolor{grey40}{rgb}{0.40,0.40,0.40}
\definecolor{grey41}{rgb}{0.41,0.41,0.41}
\definecolor{grey42}{rgb}{0.42,0.42,0.42}
\definecolor{grey43}{rgb}{0.43,0.43,0.43}
\definecolor{grey44}{rgb}{0.44,0.44,0.44}
\definecolor{grey45}{rgb}{0.45,0.45,0.45}
\definecolor{grey46}{rgb}{0.46,0.46,0.46}
\definecolor{grey47}{rgb}{0.47,0.47,0.47}
\definecolor{grey48}{rgb}{0.48,0.48,0.48}
\definecolor{grey49}{rgb}{0.49,0.49,0.49}
\definecolor{grey4}{rgb}{0.04,0.04,0.04}
\definecolor{grey50}{rgb}{0.50,0.50,0.50}
\definecolor{grey51}{rgb}{0.51,0.51,0.51}
\definecolor{grey52}{rgb}{0.52,0.52,0.52}
\definecolor{grey53}{rgb}{0.53,0.53,0.53}
\definecolor{grey54}{rgb}{0.54,0.54,0.54}
\definecolor{grey55}{rgb}{0.55,0.55,0.55}
\definecolor{grey56}{rgb}{0.56,0.56,0.56}
\definecolor{grey57}{rgb}{0.57,0.57,0.57}
\definecolor{grey58}{rgb}{0.58,0.58,0.58}
\definecolor{grey59}{rgb}{0.59,0.59,0.59}
\definecolor{grey5}{rgb}{0.05,0.05,0.05}
\definecolor{grey60}{rgb}{0.60,0.60,0.60}
\definecolor{grey61}{rgb}{0.61,0.61,0.61}
\definecolor{grey62}{rgb}{0.62,0.62,0.62}
\definecolor{grey63}{rgb}{0.63,0.63,0.63}
\definecolor{grey64}{rgb}{0.64,0.64,0.64}
\definecolor{grey65}{rgb}{0.65,0.65,0.65}
\definecolor{grey66}{rgb}{0.66,0.66,0.66}
\definecolor{grey67}{rgb}{0.67,0.67,0.67}
\definecolor{grey68}{rgb}{0.68,0.68,0.68}
\definecolor{grey69}{rgb}{0.69,0.69,0.69}
\definecolor{grey6}{rgb}{0.06,0.06,0.06}
\definecolor{grey70}{rgb}{0.70,0.70,0.70}
\definecolor{grey71}{rgb}{0.71,0.71,0.71}
\definecolor{grey72}{rgb}{0.72,0.72,0.72}
\definecolor{grey73}{rgb}{0.73,0.73,0.73}
\definecolor{grey74}{rgb}{0.74,0.74,0.74}
\definecolor{grey75}{rgb}{0.75,0.75,0.75}
\definecolor{grey76}{rgb}{0.76,0.76,0.76}
\definecolor{grey77}{rgb}{0.77,0.77,0.77}
\definecolor{grey78}{rgb}{0.78,0.78,0.78}
\definecolor{grey79}{rgb}{0.79,0.79,0.79}
\definecolor{grey7}{rgb}{0.07,0.07,0.07}
\definecolor{grey80}{rgb}{0.80,0.80,0.80}
\definecolor{grey81}{rgb}{0.81,0.81,0.81}
\definecolor{grey82}{rgb}{0.82,0.82,0.82}
\definecolor{grey83}{rgb}{0.83,0.83,0.83}
\definecolor{grey84}{rgb}{0.84,0.84,0.84}
\definecolor{grey85}{rgb}{0.85,0.85,0.85}
\definecolor{grey86}{rgb}{0.86,0.86,0.86}
\definecolor{grey87}{rgb}{0.87,0.87,0.87}
\definecolor{grey88}{rgb}{0.88,0.88,0.88}
\definecolor{grey89}{rgb}{0.89,0.89,0.89}
\definecolor{grey8}{rgb}{0.08,0.08,0.08}
\definecolor{grey90}{rgb}{0.90,0.90,0.90}
\definecolor{grey91}{rgb}{0.91,0.91,0.91}
\definecolor{grey92}{rgb}{0.92,0.92,0.92}
\definecolor{grey93}{rgb}{0.93,0.93,0.93}
\definecolor{grey94}{rgb}{0.94,0.94,0.94}
\definecolor{grey95}{rgb}{0.95,0.95,0.95}
\definecolor{grey96}{rgb}{0.96,0.96,0.96}
\definecolor{grey97}{rgb}{0.97,0.97,0.97}
\definecolor{grey98}{rgb}{0.98,0.98,0.98}
\definecolor{grey99}{rgb}{0.99,0.99,0.99}
\definecolor{grey9}{rgb}{0.09,0.09,0.09}
\definecolor{grey}{rgb}{0.75,0.75,0.75}
\definecolor{honeydew1}{rgb}{0.94,1.00,0.94}
\definecolor{honeydew2}{rgb}{0.88,0.93,0.88}
\definecolor{honeydew3}{rgb}{0.76,0.80,0.76}
\definecolor{honeydew4}{rgb}{0.51,0.55,0.51}
\definecolor{honeydew}{rgb}{0.94,1.00,0.94}
\definecolor{hotpink}{rgb}{1.00,0.41,0.71}
\definecolor{indianred}{rgb}{0.80,0.36,0.36}
\definecolor{ivory1}{rgb}{1.00,1.00,0.94}
\definecolor{ivory2}{rgb}{0.93,0.93,0.88}
\definecolor{ivory3}{rgb}{0.80,0.80,0.76}
\definecolor{ivory4}{rgb}{0.55,0.55,0.51}
\definecolor{ivory}{rgb}{1.00,1.00,0.94}
\definecolor{khaki1}{rgb}{1.00,0.96,0.56}
\definecolor{khaki2}{rgb}{0.93,0.90,0.52}
\definecolor{khaki3}{rgb}{0.80,0.78,0.45}
\definecolor{khaki4}{rgb}{0.55,0.53,0.31}
\definecolor{khaki}{rgb}{0.94,0.90,0.55}
\definecolor{lavenderblush}{rgb}{1.00,0.94,0.96}
\definecolor{lavender}{rgb}{0.90,0.90,0.98}
\definecolor{lawngreen}{rgb}{0.49,0.99,0.00}
\definecolor{lemonchiffon}{rgb}{1.00,0.98,0.80}
\definecolor{lightblue}{rgb}{0.68,0.85,0.90}
\definecolor{lightcoral}{rgb}{0.94,0.50,0.50}
\definecolor{lightcyan}{rgb}{0.88,1.00,1.00}
\definecolor{lightgoldenrod}{rgb}{0.93,0.87,0.51}
\definecolor{lightgoldenrod}{rgb}{0.98,0.98,0.82}
\definecolor{lightgray}{rgb}{0.83,0.83,0.83}
\definecolor{lightgreen}{rgb}{0.56,0.93,0.56}
\definecolor{lightgrey}{rgb}{0.83,0.83,0.83}
\definecolor{lightpink}{rgb}{1.00,0.71,0.76}
\definecolor{lightsalmon}{rgb}{1.00,0.63,0.48}
\definecolor{lightsea}{rgb}{0.13,0.70,0.67}
\definecolor{lightsky}{rgb}{0.53,0.81,0.98}
\definecolor{lightslate}{rgb}{0.47,0.53,0.60}
\definecolor{lightslate}{rgb}{0.47,0.53,0.60}
\definecolor{lightslate}{rgb}{0.52,0.44,1.00}
\definecolor{lightsteel}{rgb}{0.69,0.77,0.87}
\definecolor{lightyellow}{rgb}{1.00,1.00,0.88}
\definecolor{limegreen}{rgb}{0.20,0.80,0.20}
\definecolor{linen}{rgb}{0.98,0.94,0.90}
\definecolor{magenta1}{rgb}{1.00,0.00,1.00}
\definecolor{magenta2}{rgb}{0.93,0.00,0.93}
\definecolor{magenta3}{rgb}{0.80,0.00,0.80}
\definecolor{magenta4}{rgb}{0.55,0.00,0.55}
\definecolor{magenta}{rgb}{1.00,0.00,1.00}
\definecolor{maroon1}{rgb}{1.00,0.20,0.70}
\definecolor{maroon2}{rgb}{0.93,0.19,0.65}
\definecolor{maroon3}{rgb}{0.80,0.16,0.56}
\definecolor{maroon4}{rgb}{0.55,0.11,0.38}
\definecolor{maroon}{rgb}{0.69,0.19,0.38}
\definecolor{mediumaquamarine}{rgb}{0.40,0.80,0.67}
\definecolor{mediumblue}{rgb}{0.00,0.00,0.80}
\definecolor{mediumorchid}{rgb}{0.73,0.33,0.83}
\definecolor{mediumpurple}{rgb}{0.58,0.44,0.86}
\definecolor{mediumsea}{rgb}{0.24,0.70,0.44}
\definecolor{mediumslate}{rgb}{0.48,0.41,0.93}
\definecolor{mediumspring}{rgb}{0.00,0.98,0.60}
\definecolor{mediumturquoise}{rgb}{0.28,0.82,0.80}
\definecolor{mediumviolet}{rgb}{0.78,0.08,0.52}
\definecolor{midnightblue}{rgb}{0.10,0.10,0.44}
\definecolor{mintcream}{rgb}{0.96,1.00,0.98}
\definecolor{mistyrose}{rgb}{1.00,0.89,0.88}
\definecolor{moccasin}{rgb}{1.00,0.89,0.71}
\definecolor{navajowhite}{rgb}{1.00,0.87,0.68}
\definecolor{navyblue}{rgb}{0.00,0.00,0.50}
\definecolor{navy}{rgb}{0.00,0.00,0.50}
\definecolor{oldlace}{rgb}{0.99,0.96,0.90}
\definecolor{olivedrab}{rgb}{0.42,0.56,0.14}
\definecolor{orange1}{rgb}{1.00,0.65,0.00}
\definecolor{orange2}{rgb}{0.93,0.60,0.00}
\definecolor{orange3}{rgb}{0.80,0.52,0.00}
\definecolor{orange4}{rgb}{0.55,0.35,0.00}
\definecolor{orangered}{rgb}{1.00,0.27,0.00}
\definecolor{orange}{rgb}{1.00,0.65,0.00}
\definecolor{orchid1}{rgb}{1.00,0.51,0.98}
\definecolor{orchid2}{rgb}{0.93,0.48,0.91}
\definecolor{orchid3}{rgb}{0.80,0.41,0.79}
\definecolor{orchid4}{rgb}{0.55,0.28,0.54}
\definecolor{orchid}{rgb}{0.85,0.44,0.84}
\definecolor{palegoldenrod}{rgb}{0.93,0.91,0.67}
\definecolor{palegreen}{rgb}{0.60,0.98,0.60}
\definecolor{paleturquoise}{rgb}{0.69,0.93,0.93}
\definecolor{paleviolet}{rgb}{0.86,0.44,0.58}
\definecolor{papayawhip}{rgb}{1.00,0.94,0.84}
\definecolor{peachpuff}{rgb}{1.00,0.85,0.73}
\definecolor{peru}{rgb}{0.80,0.52,0.25}
\definecolor{pink1}{rgb}{1.00,0.71,0.77}
\definecolor{pink2}{rgb}{0.93,0.66,0.72}
\definecolor{pink3}{rgb}{0.80,0.57,0.62}
\definecolor{pink4}{rgb}{0.55,0.39,0.42}
\definecolor{pink}{rgb}{1.00,0.75,0.80}
\definecolor{plum1}{rgb}{1.00,0.73,1.00}
\definecolor{plum2}{rgb}{0.93,0.68,0.93}
\definecolor{plum3}{rgb}{0.80,0.59,0.80}
\definecolor{plum4}{rgb}{0.55,0.40,0.55}
\definecolor{plum}{rgb}{0.87,0.63,0.87}
\definecolor{powderblue}{rgb}{0.69,0.88,0.90}
\definecolor{purple1}{rgb}{0.61,0.19,1.00}
\definecolor{purple2}{rgb}{0.57,0.17,0.93}
\definecolor{purple3}{rgb}{0.49,0.15,0.80}
\definecolor{purple4}{rgb}{0.33,0.10,0.55}
\definecolor{purple}{rgb}{0.63,0.13,0.94}
\definecolor{red1}{rgb}{1.00,0.00,0.00}
\definecolor{red2}{rgb}{0.93,0.00,0.00}
\definecolor{red3}{rgb}{0.80,0.00,0.00}
\definecolor{red4}{rgb}{0.55,0.00,0.00}
\definecolor{red}{rgb}{1.00,0.00,0.00}
\definecolor{rosybrown}{rgb}{0.74,0.56,0.56}
\definecolor{royalblue}{rgb}{0.25,0.41,0.88}
\definecolor{saddlebrown}{rgb}{0.55,0.27,0.07}
\definecolor{salmon1}{rgb}{1.00,0.55,0.41}
\definecolor{salmon2}{rgb}{0.93,0.51,0.38}
\definecolor{salmon3}{rgb}{0.80,0.44,0.33}
\definecolor{salmon4}{rgb}{0.55,0.30,0.22}
\definecolor{salmon}{rgb}{0.98,0.50,0.45}
\definecolor{sandybrown}{rgb}{0.96,0.64,0.38}
\definecolor{seagreen}{rgb}{0.18,0.55,0.34}
\definecolor{seashell1}{rgb}{1.00,0.96,0.93}
\definecolor{seashell2}{rgb}{0.93,0.90,0.87}
\definecolor{seashell3}{rgb}{0.80,0.77,0.75}
\definecolor{seashell4}{rgb}{0.55,0.53,0.51}
\definecolor{seashell}{rgb}{1.00,0.96,0.93}
\definecolor{sienna1}{rgb}{1.00,0.51,0.28}
\definecolor{sienna2}{rgb}{0.93,0.47,0.26}
\definecolor{sienna3}{rgb}{0.80,0.41,0.22}
\definecolor{sienna4}{rgb}{0.55,0.28,0.15}
\definecolor{sienna}{rgb}{0.63,0.32,0.18}
\definecolor{skyblue}{rgb}{0.53,0.81,0.92}
\definecolor{slateblue}{rgb}{0.42,0.35,0.80}
\definecolor{slategray}{rgb}{0.44,0.50,0.56}
\definecolor{slategrey}{rgb}{0.44,0.50,0.56}
\definecolor{snow1}{rgb}{1.00,0.98,0.98}
\definecolor{snow2}{rgb}{0.93,0.91,0.91}
\definecolor{snow3}{rgb}{0.80,0.79,0.79}
\definecolor{snow4}{rgb}{0.55,0.54,0.54}
\definecolor{snow}{rgb}{1.00,0.98,0.98}
\definecolor{springgreen}{rgb}{0.00,1.00,0.50}
\definecolor{steelblue}{rgb}{0.27,0.51,0.71}
\definecolor{tan1}{rgb}{1.00,0.65,0.31}
\definecolor{tan2}{rgb}{0.93,0.60,0.29}
\definecolor{tan3}{rgb}{0.80,0.52,0.25}
\definecolor{tan4}{rgb}{0.55,0.35,0.17}
\definecolor{tan}{rgb}{0.82,0.71,0.55}
\definecolor{thistle1}{rgb}{1.00,0.88,1.00}
\definecolor{thistle2}{rgb}{0.93,0.82,0.93}
\definecolor{thistle3}{rgb}{0.80,0.71,0.80}
\definecolor{thistle4}{rgb}{0.55,0.48,0.55}
\definecolor{thistle}{rgb}{0.85,0.75,0.85}
\definecolor{tomato1}{rgb}{1.00,0.39,0.28}
\definecolor{tomato2}{rgb}{0.93,0.36,0.26}
\definecolor{tomato3}{rgb}{0.80,0.31,0.22}
\definecolor{tomato4}{rgb}{0.55,0.21,0.15}
\definecolor{tomato}{rgb}{1.00,0.39,0.28}
\definecolor{turquoise1}{rgb}{0.00,0.96,1.00}
\definecolor{turquoise2}{rgb}{0.00,0.90,0.93}
\definecolor{turquoise3}{rgb}{0.00,0.77,0.80}
\definecolor{turquoise4}{rgb}{0.00,0.53,0.55}
\definecolor{turquoise}{rgb}{0.25,0.88,0.82}
\definecolor{violetred}{rgb}{0.82,0.13,0.56}
\definecolor{violet}{rgb}{0.93,0.51,0.93}
\definecolor{wheat1}{rgb}{1.00,0.91,0.73}
\definecolor{wheat2}{rgb}{0.93,0.85,0.68}
\definecolor{wheat3}{rgb}{0.80,0.73,0.59}
\definecolor{wheat4}{rgb}{0.55,0.49,0.40}
\definecolor{wheat}{rgb}{0.96,0.87,0.70}
\definecolor{whitesmoke}{rgb}{0.96,0.96,0.96}
\definecolor{white}{rgb}{1.00,1.00,1.00}
\definecolor{yellow1}{rgb}{1.00,1.00,0.00}
\definecolor{yellow2}{rgb}{0.93,0.93,0.00}
\definecolor{yellow3}{rgb}{0.80,0.80,0.00}
\definecolor{yellow4}{rgb}{0.55,0.55,0.00}
\definecolor{yellowgreen}{rgb}{0.60,0.80,0.20}
\definecolor{yellow}{rgb}{1.00,1.00,0.00}
 
\usepackage[usenames,dvipsnames]{xcolor}
\usepackage[psamsfonts]{amsfonts}
\usepackage{graphicx}
\usepackage{amssymb}
\usepackage{shortvrb} 
 
\usepackage{hyperref} 
\hypersetup{pdfpagemode=FullScreen,   
backref,
colorlinks=true,  
citecolor=Bittersweet,
linkcolor=Bittersweet, 
urlcolor=Bittersweet
}
 
\setattribute{journal}{name}{}

\RequirePackage{hypernat}

\usepackage{hyperref} 
\hypersetup{pdfpagemode=FullScreen,   
backref,
colorlinks=true, 
citecolor=Bittersweet,
linkcolor=Bittersweet, 
urlcolor=Bittersweet
}

\usepackage[psamsfonts]{amsfonts}

\usepackage{pdflscape}

\usepackage{dsfont}
 
\defcitealias{Cut2015}{CPV15}

\newtheorem{Lem}{Lemma}[section]

\newtheorem{Theor}{Theorem}[section]

\newcommand{\cqfd}{\hfill $\square$}

\newcommand{\R}{\mathbb R}

\newcommand{\n}{^{(n)}}

\newcommand{\Xb}{\mathbf{X}}
\newcommand{\Sb}{\mathbf{S}}

\newcommand{\Vb}{\mathbf{V}}

\newcommand{\ub}{\ensuremath{\mathbf{u}}}
\newcommand{\Tb}{\ensuremath{\mathbf{T}}}
\newcommand{\vb}{\ensuremath{\mathbf{v}}}

\newcommand{\xb}{\ensuremath{\mathbf{x}}}

\newcommand{\Ab}{\ensuremath{\mathbf{A}}}

\newcommand{\Ub}{\ensuremath{\mathbf{U}}}
\newcommand{\Mb}{\ensuremath{\mathbf{M}}}
\newcommand{\Wb}{\ensuremath{\mathbf{W}}}

\newcommand{\Gb}{\ensuremath{\mathbf{G}}}

\newcommand{\thetab}{{\pmb \theta}}

\newcommand{\Umb}{{\pmb \Upsilon}}

\newcommand{\Sigb}{{\pmb \Sigma}}

\newcommand{\Deltab}{{\pmb \Delta}}
\newcommand{\taub}{{\pmb \tau}}

\newcommand{\Gamb}{{\pmb \Gamma}}

\newcommand{\pr}{^{\prime}}
\newcommand{\ny}{n\rightarrow\infty}

\begin{document}
\begin{frontmatter}

\title{Sign Tests for Weak Principal Directions}
\runtitle{Sign tests for weak principal directions}

\begin{aug}
\author{\fnms{Davy} \snm{Paindaveine}$^\dagger$, \fnms{Julien} \snm{Remy}$^\ddagger$ and \fnms{Thomas} \snm{Verdebout}\thanksref{t1}}

\thankstext{t1}{Corresponding author.}
\runauthor{D. Paindaveine, J. Remy and Th. Verdebout}

\affiliation{Universit\'{e} libre de Bruxelles}

\address{$^{\dagger \ddagger *}$Universit\'{e} libre de Bruxelles\\
ECARES and D\'{e}partement de Math\'{e}matique\\
Avenue F.D. Roosevelt, 50\\
ECARES, CP114/04\\
B-1050, Brussels\\
Belgium\\
}

\address{$^\dagger$Universit\'{e} Toulouse Capitole\\
Toulouse School of Economics\\
21, All\'{e}e de Brienne\\
31015 Toulouse Cedex 6\\
France\\
}

\end{aug}
\vspace{3mm}

\begin{abstract}

We consider inference on the first principal direction of a $p$-variate elliptical distribution. We do so in challenging double asymptotic scenarios for which this direction eventually fails to be identifiable. In order to achieve robustness not only with respect to such weak identifiability but also with respect to heavy tails, we focus on sign-based statistical procedures, that is, on procedures that involve the observations only through their direction from the center of the distribution. We actually consider the generic problem of testing the null hypothesis that the first principal direction coincides with a given direction of~$\R^p$. We first focus on weak identifiability setups involving single spikes (that is, involving spectra for which the smallest eigenvalue has multiplicity~$p-1$). We show that, irrespective of the degree of weak identifiability, such setups offer local alternatives for which the corresponding sequence of statistical experiments converges in the Le Cam sense. Interestingly, the limiting experiments depend on the degree of weak identifiability. We exploit this convergence result to build optimal sign tests for the problem considered. In classical asymptotic scenarios where the spectrum is fixed, these tests are shown to be asymptotically equivalent to the sign-based likelihood ratio tests available in the literature. Unlike the latter, however, the proposed sign tests are robust to arbitrarily weak identifiability. We show that our tests meet the asymptotic level constraint irrespective of the structure of the spectrum, hence also in possibly multi-spike setups. We fully characterize the non-null asymptotic distributions of the corresponding test statistics under weak identifiability, which allows us to quantify the corresponding local asymptotic powers. Finally, Monte Carlo exercises are conducted to assess the finite-sample relevance of our asymptotic results and a real-data illustration is provided.  
\end{abstract}

\begin{keyword}[class=MSC]
\kwd[Primary ]{62F05, 62H25}
\kwd[; secondary ]{62E20}
\end{keyword}

\begin{keyword}
\kwd{Le Cam's asymptotic theory of statistical experiments}
\kwd{Local asymptotic normality}
\kwd{Principal component analysis}
\kwd{Sign tests}
\kwd{Weak identifiability}
\end{keyword}

\end{frontmatter} 

\section{Introduction} 
\label{sec:intro}

Most classical methods in multivariate statistics are based on Gaussian maximum likelihood estimators of location and scatter, that is, on the sample mean and sample covariance matrix. Irrespective of the considered problem (location or scatter problems, multivariate regression problems, principal component analysis, canonical correlation analysis, etc.), these methods exhibit poor efficiency properties when Gaussian assumptions are violated, particularly when the underlying distributions have heavy tails. Moreover, the resulting procedures are very sensitive to possible outliers in the data.

To improve on this, many robust procedures were developed. In particular, multivariate sign methods, that is, methods that use the observations only through their direction from the center of the distribution, have become increasingly popular. For location problems, multivariate sign tests were considered in~\cite{Ran89}, \cite{mooj95},~\cite{HP02} and \cite{PaiVer2016}, whereas sign procedures for scatter or shape matrices were considered in \cite{Tyl1987}, \cite{Dum98}, \cite{HP06}, \cite{DuFri15} and \cite{DuFri16}, to cite only a few. PCA techniques based on multivariate signs (and on the companion concept of ranks) were studied in \cite{HPV10}, \cite{Taski12}, \cite{HPV13} and \cite{DuTyl16}. Multivariate sign tests were also developed, e.g., for testing i.i.d.-ness against serial dependence (see~\citealp{Pai2009}), or for testing for multivariate independence (see \citealp{Tasetal2003} and \citealp{Tasetal2005}). Most references above actually focus on \emph{spatial sign procedures}, that is, on procedures that are based on the signs $\Ub_i:=\Xb_i/\|\Xb_i\|,$ $i=1,\ldots,n$, obtained by projecting the $p$-variate observations at hand onto the unit sphere~${\cal S}^{p-1}:= \{\xb \in \R^p: \|\xb\|^2=\xb\pr\xb=1\}$ of~$\R^p$ (sometimes, the projection is performed on standardized observations in order to achieve affine invariance).  Since they discard the radii~$\|\Xb_i\|$, $i=1,\ldots,n$, spatial sign procedures can deal with arbitrarily heavy tails and are robust to observations that would be far from the center of the distribution. Furthermore, they are by nature well adapted to directional or axial data for which these radii are not observed; see, e.g., \cite{MarJup2000} or \cite{LV17}. For more details on spatial sign methods, we refer to the monograph \cite{oja2010book}.

The present paper considers principal component analysis, or more precisely, inference on principal component directions. Consider the case where the observations~$\Xb_1, \ldots, \Xb_n$ at hand form a random sample from a centered $p$-variate elliptical distribution, that is, a distribution whose characteristic function is of the form~${\bf t}\mapsto \exp({\bf t}\pr\Sigb{\bf t})$ for some symmetric positive definite $p\times p$ matrix~$\Sigb$. Assume that the ordered eigenvalues of~$\Sigb$ satisfy~$\lambda_{1}> \lambda_{2}\geq  \ldots \geq \lambda_{p}$, so that the leading eigenvector~$\thetab_1$, possibly unlike the other eigenvectors~$\thetab_j$, $j=2,\ldots,p$, is identifiable (as usual, identifiability is up to an unimportant sign). We will then throughout consider the problem of testing the null hypothesis~${\cal H}_0: \thetab_1= \thetab_0$ against the alternative hypothesis~${\cal H}_1: \thetab_1 \neq \thetab_0$, where~$\thetab_0$ is a fixed unit $p$-vector. 
This testing problem has attracted much attention in the past decades, and the textbook procedure, namely the \cite{And63} Gaussian likelihood ratio test, has been extended in various directions. To mention only a few, \cite{Jo84} considered a small-sample test, whereas \cite{Flu88} proposed an extension to a larger number of eigenvectors.  \cite{Tyl81,Tyl83} robustified the \cite{And63} test to possible elliptical departures from multinormality (the original likelihood ratio test require Gaussian assumptions). \cite{Sc08} considered extensions to the case of Gaussian random matrices, and \cite{HPV10} obtained Le Cam optimal tests for the problem considered.

All asymptotic tests above assume that the eigenvalues~$\lambda_{1}> \lambda_{2}\geq  \ldots \geq \lambda_{p}$ are fixed, so that the eigenvector~$\thetab_1$ remains asymptotically identifiable. The null asymptotic distribution of the corresponding test statistics, however, may very poorly approximate their \mbox{fixed-$n$} distribution when~$\lambda_1/\lambda_2$ is close to one. In the present work, we therefore consider general asymptotic scenarios that address this issue. More precisely, we allow for scatter values~${\Sigb}={\Sigb}_n$ for which the corresponding ratio~$\lambda_{n1}/\lambda_{n2}(>1)$ converges to one as~$n$ diverges to infinity. In such asymptotic scenarios, the leading eigenvector~$\thetab_{n1}$ remains identifiable for any~$n$ but is no longer identifiable in the limit. One then says that~$\thetab_{n1}$ is \emph{weakly identifiable}. The distributional framework considered here formalizes situations that are often encountered in practice where two sample eigenvalues are close to each other so that inference about the corresponding eigenvectors is a priori difficult. Inference on weakly identified parameters has already been much considered in the literature: see, e.g., \cite{Pot2002}, \cite{ForHil2003}, \cite{Duf2006} and \cite{BeLa14}. \cite{LiSha03} and  \cite{ZhZh06} consider asymptotic inference under a total lack of identifiability. Recently, \cite{PRV18} considered Gaussian tests on weakly identified eigenvectors. While these tests can handle weak identifiability, they are based on sample covariance matrices, hence cannot deal with heavy tails and are very sensitive to possible outliers. In the present work, we tackle the same problem but develop spatial sign tests that not only inherit the robustness of sign procedures but also can deal with both heavy tails and weak identifiability.

To ensure that spatial signs are well-defined with probability one, we will restrict to elliptical distributions that do not attribute a positive probability mass to the symmetry center (the symmetry center will be assumed to be known and, without any loss of generality, to coincide with the origin of~$\R^p$---extension to the unknown location case is straightforward, as we will explain in Section~\ref{sec:wrapup}). 
Suitable spatial sign tests are to be determined in the image of the model by the projection onto the unit sphere~$\mathcal{S}^{p-1}$. Now, if the random $p$-vector~$\Xb$ is centered elliptical with scatter matrix~$\Sigb$, then the corresponding spatial sign~$\Ub:=\Xb/\|\Xb\|$ follows the \emph{angular Gaussian distribution} with shape matrix~$\Vb:=p\Sigb/{\rm tr}(\Sigb)$; see \cite{Tyl87ang}. Since~$\Sigb$ and~$\Vb$ share the same eigenvectors, the induced testing problem is still the problem of testing~$\mathcal{H}_0:\thetab_1=\thetab_0$ against~$\mathcal{H}_1:\thetab_1\neq \thetab_0$, where~$\thetab_1$ is the leading eigenvector of~$\Vb$. Also, since eigenvalues of~$\Sigb$ and~$\Vb$ are equal up to a common positive factor, (weak) identifiability occurs for the original elliptical problem if and only if it does for the induced angular Gaussian problem. These considerations explain that identifying optimal spatial sign tests under weak identifiability for the elliptical problem should be done in the setup where one observes triangular arrays of spatial signs~$\Ub_{n1},\ldots,\Ub_{nn}$, $n=1,2,\ldots$, randomly drawn from the angular Gaussian distribution with shape matrix~$\Vb_n$.  

Accordingly, we consider the problem of testing the null hypothesis~$\mathcal{H}_0:\thetab_{n1}=\thetab_0$ in an angular Gaussian double asymptotic scenario for which~$\lambda_{n1}/\lambda_{n2}(>1)$ may converge to one at an arbitrary rate, which provides weak identifiability for the leading principal direction~$\thetab_{n1}$. We first focus on sequences of single-spike shape matrices (characterized by spectra of the form~$\lambda_{n1}>\lambda_{n2}=\ldots=\lambda_{np}$). We will show that, irrespective of the rate of convergence of~$\lambda_{n1}/\lambda_{n2}$ to one, there exist suitable local alternatives for which the corresponding sequence of statistical experiments converges in the Le Cam sense. Interestingly, the limiting experiment depends on the degree of weak identifiability. This paves the way, in this single-spike setup, to optimal testing for the problem considered. Quite nicely, we will actually build tests that are Le Cam optimal in single-spike setups while remaining valid (in the sense that they meet the asymptotic level constraint) under general, ``multi-spike", spectra. We will also fully characterize the non-null behavior of these tests under weak identifiability, which will allow us to extensively quantify their local asymptotic powers. 

The outline of the paper is as follows. In Section~\ref{sec:LAN}, we introduce the notation, describe the sequence of angular Gaussian models to be considered, and derive the aforementioned results on limiting experiments. This is used to derive a sign test that can deal with weak identifiability and enjoys nice Le Cam optimality properties. However, (i) since this test involves nuisance parameters, it is an infeasible statistical procedure. Moreover, (ii) the optimal sign test requires a single-spike spectrum. In Section~\ref{sec:thetest}, we therefore derive a version of the optimal test that~(i) is feasible and~(ii) should be able to cope with multi-spike spectra. We show that, under the null hypothesis, hence also under sequences of contiguous hypotheses, the infeasible and feasible tests are asymptotically equivalent, so that the latter inherits the optimality properties of the former. We also show that, in classical asymptotic scenarios where one stays away from weak identifiability, the proposed test is asymptotically equivalent to the likelihood ratio test from \cite{Tyl87ang}. In Section~\ref{sec:asymptprop}, we investigate the asymptotic properties of the proposed sign test. We first show that, as anticipated above, this test asymptotically achieves the target null size even when the underlying scatter matrix does not have a single-spike structure.
Then, for any degree of weak identifiability, we derive the asymptotic distribution of the proposed test statistic under suitable local alternatives. In Section~\ref{sec:simus}, Monte Carlo exercises are conducted (a) to show that the proposed sign test, unlike its competitors, can deal with both heavy tails and weak identifiability and (b) to compare the finite-sample powers of the various tests. A real-data illustration is presented in Section~\ref{sec:real data}. Finally, conclusions and final comments are provided in Section~\ref{sec:wrapup}. All proofs are collected in a technical appendix.

\section{Limits of angular Gaussian experiments}
\label{sec:LAN}

In this section, our objective is to derive the form of locally asymptotically optimal tests for~$\mathcal{H}_0:\thetab_{n1}=\thetab_0$ in asymptotic scenarios under which~$\thetab_{n1}$ is weakly identifiable. To do so, we study sequences of angular Gaussian experiments indexed by shape matrices~$\Vb_n$ with eigenvalues $\lambda_{n1}>\lambda_{n2}=\ldots=\lambda_{np}$ satisfying~$\lambda_{n1}/\lambda_{n2}\to 1$. Since the optimal tests we will obtain in this section are actually infeasible, we will then construct in Section~\ref{sec:thetest} a practical sign test that achieves the same (null and non-null) asymptotic properties as the infeasible tests; the present section can therefore be seen as a stepping stone to the tests proposed in Section~\ref{sec:thetest}.

Consider a triangular array of observations $\Xb_{ni}$, $i=1,\ldots,n$, $n=1,2,\ldots,$ such that for any $n$, the random $p$-vectors $\Xb_{n1},\ldots,\Xb_{nn}$ form a random sample from a centered elliptical distribution that does not attribute a positive probability mass to the origin of~$\R^p$. For any~$n$, the leading eigenvector~$\thetab_{n1}$ of the corresponding scatter matrix~${\Sigb}_n$ is assumed to be well identified, up to a sign. In this general setup, we aim at designing suitable sign tests for the problem of testing the null hypothesis~$\mathcal{H}_0:\thetab_{n1}=\thetab_0$ against the alternative hypothesis~$\mathcal{H}_1:\thetab_{n1}\neq\thetab_0$, where~$\thetab_0$ is a fixed unit $p$-vector. As explained in the Introduction, such tests should be determined in the image of the model by the projection of the model onto the unit sphere~$\mathcal{S}^{p-1}$. For any~$n$, the projected observations~$\Ub_{n1},\ldots,\Ub_{nn}$, where~$\Ub_{ni}=\Xb_{ni}/\|\Xb_{ni}\|$ is the spatial sign of~$\Xb_{ni}$, form a random sample from the $p$-variate angular Gaussian distribution with shape matrix~${\bf V}_n= p{\Sigb}_n/{\rm tr}({\Sigb}_n)$ (shape matrices throughout are normalized to have trace~$p$), so that the density of~$\Ub_{n1}$ with respect to the surface area measure on ${\cal S}^{p-1}$ is
\begin{equation}
\label{angdens}
\ub (\in {\cal S}^{p-1})\mapsto \frac{\Gamma(\frac{p}{2})}{2 \pi^{p/2} ({\rm det}\,{\bf V}_n)^{1/2}} (\ub\pr {\bf V}_n^{-1} \ub)^{-p/2},
\end{equation}
where $\Gamma(\cdot)$ is Euler's gamma function; see \cite{Tyl87ang}. We will denote the corresponding hypothesis as~${\rm P}\n_{{\bf V}_{n}}$ (at places, ${\rm P}\n_{{\bf V}_{n}}$ will also denote the distribution of a random sample of size~$n$ from the angular Gaussian distribution with shape matrix~${\bf V}_{n}$). In this angular Gaussian framework, the aforementioned elliptical testing problem induces the problem of testing the null hypothesis~${\cal H}_0: \thetab_{n1}= \thetab_0$ against the alternative hypothesis~${\cal H}_1: \thetab_{n1} \neq \thetab_0$, where~$\thetab_{n1}$ is the leading eigenvector of~$\Vb_n$ (recall that~$\Vb_n$ and~$\Sigb_n$ share the same eigenvectors). Throughout, $\lambda_{nj}$ and~$\thetab_{nj}$ will refer to the (ordered) eigenvalues and eigenvectors of~$\Vb_n$. In the original elliptical problem, weak identifiability of~$\thetab_{n1}$, meaning that~$\thetab_{n1}$ is identified for any fixed~$n$ but is not in the limit as $\ny$, clearly occurs if and only if~$\lambda_{n1}/\lambda_{n2}\to 1$.

The goal of this section is to determine optimal sign tests, under possibly weak identifiability of~$\thetab_{n1}$, in the particular case of single-spike spectra, that is, in situations where the smallest eigenvalue has multiplicity~$p-1$. These tests will result from the asymptotic study of angular Gaussian log-likelihood ratios in Theorem \ref{TheorLAN} below. We will consider sequences of null shape matrices of the form
\begin{equation} 
\label{V0null}
{\bf V}_{0n}
:=
\Big(
1- \frac{\delta_{n}\xi}{p}
\Big) 
{\bf I}_p + \delta_{n}\xi \thetab_0 \thetab_0\pr,
\end{equation}
where $\xi>0$ is a \emph{locality parameter} and $\delta_n$ is a bounded positive sequence that may be~$o(1)$; throughout, we tacitly assume that~$\xi$ and $\delta_n$ are chosen so that ${\bf V}_{0n}$ is positive definite. It is straigthforward to check that~${\bf V}_{0n}$ indeed has a single spike: the largest eigenvalue is~$\lambda_{n1}=1+(p-1)\delta_{n}\xi/p$, with multiplicity one and corresponding eigenvector~$\thetab_0$, whereas the remaining eigenvalues are~$\lambda_{n2}=\ldots=\lambda_{np}=1-\delta_{n}\xi/p$, $j=2,\ldots,p$, with an eigenspace that is the orthogonal complement to~$\thetab_0$. In this setup, weak identifiability of~$\thetab_{n1}$ occurs if and only if~$\delta_n$ is~$o(1)$. 

To discuss optimality issues, we consider local alternatives associated with perturbations~$\thetab_0+ \nu_n \taub_n$ of~$\thetab_0$, where~$(\nu_n)$ is a positive sequence and~$(\taub_n)$ is a bounded sequence in $\R^p$ such that~$\thetab_0+ \nu_n \taub_n \in {\cal S}^{p-1}$ for any~$n$. It is easy to show that the latter condition entails that $\nu_n$ and~$\taub_n$ must satisfy 
\begin{equation} \label{pert}
\thetab_0\pr \taub_n
= 
-\frac{\nu_n}{2}  \| \taub_n \|^2
\end{equation}
for any~$n$. 
The resulting sequence of alternatives is then associated with the single-spike shape matrices 
\begin{equation} 
\label{V1alt}
{\bf V}_{1n}
:=
\Big(1- \frac{\delta_{n}\xi}{p}\Big) {\bf I}_p + \delta_{n}\xi (\thetab_0+ \nu_n \taub_n) (\thetab_0+ \nu_n \taub_n)\pr.
\end{equation}
Our construction of optimal sign tests requires studying the asymptotic behavior, under~${\rm P}\n_{{\bf V}_{0n}}$, of the log-likelihood ratios $\Lambda_n:= \log (d{\rm P}_{{\bf V}_{1n}}\n/d{\rm P}\n_{{\bf V}_{0n}})$. To do so, let 
\begin{equation} 
\label{defingamma}
\gamma_n
:=
\frac{p \delta_{n}\xi}{p+(p-1) \delta_{n}\xi}
\quad\textrm{and}\quad
\Sb_n (\Vb)
:= 
\frac{1}{n} \sum_{i=1}^n \frac{\Vb^{-1/2}\Ub_{ni}\Ub_{ni}\pr \Vb^{-1/2}}{ \| \Vb^{-1/2} \Ub_{ni}\|^2}
\cdot
\end{equation}
Note that the sequence~$(\gamma_n)$ is~$O(\delta_n)$. We then have the following result. 

\begin{Theor} 
\label{TheorLAN} 
Fix $\thetab_0 \in {\cal S}^{p-1}$. Let~$(\nu_n)$ be a positive sequence and~$(\taub_n)$ be a bounded sequence in $\R^p$ such that~$\thetab_0+\nu_n\taub_n\in\mathcal{S}^{p-1}$ for any~$n$. Let~${\bf V}_{0n}$ be as in~\eqref{V0null}. Then, as $n\to\infty$, under ${\rm P}\n_{{\bf V}_{0n}}$, we have the following: 
\begin{itemize}
\item[(i)] if $\delta_n \equiv 1$, then, for~$\nu_n=1/(\sqrt{n} \gamma_n)$,
\begin{equation}
\label{LAQi}	
\Lambda_n= \taub_n\pr \Deltab_{\thetab_0}^{(i)}- \frac{1}{2} \taub_n\pr \Gamb_{\thetab_0}^{(i)} \taub_n+ o_{\rm P}(1)
\end{equation}
as $\ny$, where we let
$$
\Gamb_{\thetab_0}^{(i)}:=\frac{p (p+(p-1)\xi)}{(p+2)(p-\xi)} ({\bf I}_p- \thetab_0 \thetab_0\pr)
,
$$
and where
$$
\Deltab_{\thetab_0}^{(i)}:=p \frac{\sqrt{p+(p-1)\xi}}{\sqrt{p-\xi}} 
({\bf I}_{p}- \thetab_0 \thetab_0\pr) 
\sqrt{n} 
\big(
\Sb_n({\bf V}_{0n})- {\textstyle{\frac{1}{p}}} {\bf I}_p
\big) 
\thetab_0
$$
is asymptotically normal with mean zero and covariance matrix $\Gamb_{\thetab_0}^{(i)}$;
\item[(ii)] if $\delta_n$ is $o(1)$ with $\sqrt{n} \delta_n \to \infty$,  then, for~$\nu_n=1/(\sqrt{n} \gamma_n)$,
\begin{equation}
\label{LAQii}	
\Lambda_n= \taub_n\pr \Deltab_{\thetab_0}^{(ii)}- \frac{1}{2} \taub_n\pr \Gamb_{\thetab_0}^{(ii)} \taub_n+ o_{\rm P}(1)
\end{equation}
as $\ny$, where we let
\begin{equation}
	\label{infoii}
\Gamb_{\thetab_0}^{(ii)}:=\frac{p}{p+2} ({\bf I}_p- \thetab_0 \thetab_0\pr)
\end{equation}
and where 
\begin{equation}
	\label{centralsii}
\Deltab_{\thetab_0}^{(ii)}:=p ({\bf I}_{p}- \thetab_0 \thetab_0\pr) \sqrt{n} 
\big(
\Sb_n({\bf V}_{0n})- {\textstyle{\frac{1}{p}}} {\bf I}_p
\big) 
\thetab_0
\end{equation}
is asymptotically normal with mean zero and covariance matrix $\Gamb_{\thetab_0}$;
\item[(iii)] if~$\delta_n=1/\sqrt{n}$, then, for~$\nu_n= 1/(\sqrt{n} \gamma_n)$ $($or equivalently $\nu_n \equiv 1/\xi)$,
\begin{equation} 
\label{contigdec}
\Lambda_n= \taub_n \pr \Umb_{\thetab_0}\n \thetab_0+ \frac{1}{2\xi}\, \taub_n\pr \Umb_{\thetab_0}\n \taub_n- \frac{p}{2(p+2)}
\Big(\| \taub_n\|^2 - \frac{1}{4\xi^2} \| \taub_n\|^4\Big)
+
o_{\rm P}(1)
,
\end{equation}
where 
$$
\Umb_{\thetab_0}\n
:=
p \sqrt{n} 
\big(
\Sb_n({\bf V}_{0n})- {\textstyle{\frac{1}{p}}} {\bf I}_p
\big) 
$$ 
is such that ${\rm vec}(\Umb_{\thetab_0}\n)$ is asymptotically normal with mean zero and covariance matrix $(p/(p+2))({\bf I}_{p^2}+ {\bf K}_p+ {\bf J}_p)- {\bf J}_p$;
\item[(iv)] if~$\sqrt{n} \delta_n \to 0$, then, even for~$\nu_n \equiv 1/\xi$, we have $\Lambda_n= o_{\rm P}(1)$.
\end{itemize} 
\end{Theor}

Theorem~\ref{TheorLAN} shows that the asymptotic behavior of~$\Lambda_n$ crucially depends on the sequence~$(\delta_n)$ and identifies four different regimes. A similar phenomenon has been obtained in \cite{Tyl83mu} when investigating the limiting behavior of eigenvalues of scatter estimators (in particular, the largest eigenvalue of the scatter estimators considered in \cite{Tyl83mu} shows a limiting behavior that  depends on~$\delta_n$ and, parallel to what we have in the present work, $\delta_n=1/\sqrt{n}$ also turns out to be an important threshold there). In the ``classical" regime (i) where~\mbox{$\delta_n \equiv 1$}, standard perturbations with~$\nu_n\sim 1/\sqrt{n}$ provide a sequence of experiments that is locally asymptotically normal (LAN), with central sequence $\Deltab_{\thetab_0}^{(i)}$ and Fisher information $\Gamb_{\thetab_0}^{(i)}$. In such a LAN setup, the locally asymptotically maximin test~$\phi\n$ (see, e.g., Section~5.2.3 from \cite{LV17} for the concept of maximin tests) for ${\cal H}_0: \thetab_{n1}= \thetab_0$ against ${\cal H}_1: \thetab_{n1} \neq \thetab_0$ rejects the null hypothesis at asymptotic level $\alpha$ when 
\begin{equation} 
\label{thetest}
T_n(\Vb_{0n})
:=
(\Deltab_{\thetab_0}^{(i)})\pr (\Gamb_{\thetab_0}^{(i)})^-\Deltab_{\thetab_0}^{(i)}=n p(p+2) 
\| 
({\bf I}_p- \thetab_0 \thetab_0\pr)  \Sb_n(\Vb_{0n}) \thetab_0
\|^2
 > \chi_{p-1,1-\alpha}^{2},
\end{equation}
where~$\Ab^-$ stands for the Moore-Penrose inverse of $\Ab$ and where~$\chi_{\ell,1-\alpha}^{2}$ denotes the upper $\alpha$-quantile of the chi-square distribution with~$\ell$ degrees of freedom. 
In regime~(ii), perturbations with~$\nu_n\sim 1/(\sqrt{n}\delta_n)$---that is, perturbations that are more severe than in the standard regime~(i)---make the sequence of experiments LAN, here with the central sequence $\Deltab_{\thetab_0}^{(ii)}$ in~(\ref{centralsii}) and the Fisher information matrix $\Gamb_{\thetab_0}^{(ii)}$ in~(\ref{infoii}). Since
$$
(\Deltab_{\thetab_0}^{(ii)})\pr (\Gamb_{\thetab_0}^{(ii)})^-\Deltab_{\thetab_0}^{(ii)}=T_n(\Vb_{0n})
,
$$ 
the test~$\phi\n$ is still locally asymptotically maximin in regime~(ii). 

The situation in regime~(iii) is quite different. While the sequence of experiments there is not LAN nor LAMN (locally asymptotically mixed normal), it still converges in the Le Cam sense. It is easy to check that, in this regime, 
$$
\taub_n \pr \Umb_{\thetab_0}\n \thetab_0+ \frac{1}{2\xi} \, \taub_n\pr \Umb_{\thetab_0}\n \taub_n
$$
is asymptotically normal with mean zero and variance 
$$ 
\frac{p}{p+2}
\Big(\| \taub_n\|^2 - \frac{1}{4\xi^2} \| \taub_n\|^4\Big)
,
$$ 
so that the Le Cam first lemma entails that the sequences of hypotheses~${\rm P}_{{\Vb}_{0n}}\n$ and~${\rm P}_{{\Vb}_{1n}}\n$ are mutually contiguous in regime~(iii), too. As we will show in the next section, the test~$\phi\n$ shows non-trivial asymptotic powers under these contiguous alternatives. Finally, in regime~(iv), Theorem~\ref{TheorLAN} shows that no test can discriminate between the null hypothesis and the alternatives associated with~$\nu_n\sim 1$, which are the most severe ones that can be considered.


\section{The proposed sign test}
\label{sec:thetest}

The test~$\phi\n$ from the previous section enjoys nice optimality properties. However, (i) it is unfortunately infeasible (its test statistic~$T_n(\Vb_{0n})$ indeed involves the population shape matrix~$\Vb_{0n}$, which is of course unknown in practice); moreover, (ii) the test~$\phi\n$ will in principle meet the asymptotic level constraint only under the, quite restrictive, single-spike shape structure in~(\ref{V0null}). In this section, we therefore construct a version~$\tilde\phi\n$ of~$\phi\n$ that (i) is feasible and that (ii) will be able to cope with more general, multi-spike, shape structures. 

To do so, let~$(\Vb_{0n})$ be a sequence of shape matrices associated with the null hypothesis. In other words, we assume that, for any~$n$, the shape~$\Vb_{0n}$ admits the spectral decomposition
\begin{equation}
	\label{toestimnuis}
\Vb_{0n}= \lambda_{n1} \thetab_{0} \thetab_{0}\pr + \sum_{j=2}^{p} \lambda_{nj} \thetab_{nj} \thetab_{nj}\pr
,
\end{equation}
where the~$\thetab_{nj}$'s form an orthonormal basis of the orthogonal complement to~$\thetab_0$ and where~$\lambda_{n1}>\lambda_{n2}\geq\ldots\geq \lambda_{np}$. In particular, the $p-1$ smallest eigenvalues here do not need to be equal, so that the number of ``spikes" may be arbitrary. With this notation, it is clear that estimating~$\Vb_{0n}$ requires estimating the eigenvalues~$\lambda_{n1}, \ldots, \lambda_{np}$ and the eigenvectors~$\thetab_{n2}, \ldots, \thetab_{np}$. To do so, we consider the \cite{Tyl1987} M-estimator, that is defined as the shape matrix satisfying~$\Sb_n(\hat{\Vb}_n)=(1/p)\mathbf{I}_p$; see~(\ref{defingamma}). Note that~$T_n(\hat{\Vb}_{n})=0$ almost surely, so that~$\hat{\Vb}_{n}$ cannot be used to estimate $\Vb_{n0}$ directly. Decompose then the M-estimator~$\hat{\Vb}_{n}$ into  
\begin{equation}
\label{DecompTyler}	
\hat{\Vb}_n
= 
\sum_{j=1}^p \hat{\lambda}_{nj} \hat{\thetab}_{nj} \hat{\thetab}_{nj}\pr
.
\end{equation}
The eigenvalues~$\hat{\lambda}_{n1},\ldots,\hat{\lambda}_{np}$ of Tyler's M-estimator provide estimates of the eigenvalues in~(\ref{toestimnuis}). Later asymptotic results, however, will require that the estimators of the eigenvectors~${\thetab}_{nj}$ are orthogonal to the null value~$\thetab_0$ of the first eigenvector~$\thetab_{n1}$, a constraint that the eigenvectors~$\hat{\thetab}_{nj}$ do not meet in general. To correct for this, we will rather use the estimators~$\tilde{\thetab}_{nj}$, $j=2,\ldots,p$, resulting from a Gram-Schmidt orthogonalization of~$\thetab_0,\hat{\thetab}_{n2}, \ldots, \hat{\thetab}_{np}$. In other words, $\tilde{\thetab}_{nj}$ is defined recursively through 
\begin{equation}
\label{GS}
\tilde{\thetab}_{nj}
:=
\frac{({\bf I}_p- \thetab_0\thetab_0'- \sum_{k=2}^{j-1}\tilde{\thetab}_{nk} \tilde{\thetab}_{nk}\pr)\hat{\thetab}_{nj}}{\| ({\bf I}_p- \thetab_0\thetab_0'- \sum_{k=2}^{j-1}\tilde{\thetab}_{nk} \tilde{\thetab}_{nk}\pr)\hat{\thetab}_{nj} \|}
,
\qquad
j=2, \ldots, p
,
\vspace{2mm}
 \end{equation}
with summation over an empty collection of indices being equal to zero.

In the rest of the paper, $\tilde{\Vb}_{0n}$ will denote the estimator of~$\Vb_{0n}$ obtained by substituting in~(\ref{toestimnuis}) the Tyler eigenvalues~$\hat{\lambda}_{nj}$ and eigenvectors~$\tilde{\thetab}_{nj}$ for the~${\lambda}_{nj}$'s and~${\thetab}_{nj}$'s. Since Tyler's M-estimator is normalized to have trace~$p$, the estimator~$\tilde{\Vb}_{0n}$ also has trace~$p$, hence is a shape matrix. Further note that if a single-spiked model as in Section \ref{sec:LAN} is assumed, which is a common practice for large~$p$, then~$\tilde{\Vb}_{0n}$ can be replaced by $\hat{\lambda}_{n1} \thetab_0 \thetab_0\pr + ((p-1)^{-1} \sum_{j=2}^p \hat{\lambda}_{nj})({\bf I}_p- \thetab_0 \thetab_0\pr)$. Quite nicely, replacing~$\Vb_{0n}$ with~$\tilde{\Vb}_{0n}$ in the test statistic~$T_n(\Vb_{0n})$ of~$\phi\n$ has no asymptotic impact in probability under the null hypothesis. More precisely, we have the following result.

\begin{Theor} 
\label{Theorasymplin}
Let $(\Vb_{0n})$ be a sequence of null shape matrices as in~(\ref{toestimnuis}). Then, under~${\rm P}_{\Vb_{0n}}\n$, $T_n(\tilde{\Vb}_{0n})=T_n(\Vb_{0n})+o_{\rm P}(1)$ as $\ny$.
\end{Theor}

Based on this result, the sign test we propose in this paper is the test~$\tilde{\phi}\n$ that rejects the null hypothesis~${\cal H}_0: \thetab_{n1}= \thetab_0$ at asymptotic level~$\alpha$ whenever
\begin{equation}
	\label{ProposedTest}
T_n(\tilde{\Vb}_{0n})
=
np(p+2)
\| 
({\bf I}_p- \thetab_0 \thetab_0\pr)  \Sb_n(\tilde{\Vb}_{0n}) \thetab_0
\|^2
>\chi_{p-1,1-\alpha}^{2}. 
\end{equation}
Unlike~${\phi}\n$, this new sign test is a feasible statistical procedure. An alternative test for the same problem is the likelihood ratio test,~$\phi_{\rm Tyl}\n$ say, that rejects the null hypothesis at asymptotic level~$\alpha$ whenever
\begin{equation}
	\label{TylerTest}
L_n
:= 
\frac{np}{p+2} 
\big(
\hat{\lambda}_{n1} \thetab_0\pr \hat{\Vb}_n^{-1}\thetab_0+ \hat{\lambda}_{n1}^{-1} \thetab_0\pr \hat{\Vb}_n \thetab_0-2
\big) 
> 
\chi_{p-1, 1- \alpha}^2
,
\end{equation}
where~$\hat{\Vb}_n$ still stands for Tyler's M-estimator. The test~$\phi_{\rm Tyl}\n$ has been proposed in \cite{Tyl87ang}. As shown by the following result, the proposed sign test~$\tilde{\phi}\n$, in classical asymptotic scenarios where one stays away from weak identifiability, is asymptotically equivalent to~$\phi_{\rm Tyl}\n$ under the null hypothesis.

\begin{Theor} 
\label{TheorequivTyler}
Let $(\Vb_{0n})$ be a sequence of null shape matrices as in~(\ref{toestimnuis}). Assume that there exists~$\eta>0$ such that both leading eigenvalues of~$\Vb_{0n}$ satisfy~$\lambda_{n1}/\lambda_{n2}\geq 1+\eta$ for~$n$ large enough. Then, under~${\rm P}_{\Vb_{0n}}\n$, $T_n(\tilde{\Vb}_{0n})=L_n+o_{\rm P}(1)$ as $\ny$.
\end{Theor}

From contiguity, this result also entails that~$\tilde{\phi}\n$ and~$\phi_{\rm Tyl}\n$ are asymptotically equivalent, hence exhibit the same asymptotic powers under the local alternatives considered in Theorem~\ref{TheorLAN}(i). There is no guarantee, however, that this asymptotic equivalence extends to the weak identifiability situations considered in Theorem~\ref{TheorLAN}(ii)-(iv). To investigate the validity of~$\tilde{\phi}\n$ under such non-standard asymptotic scenarios,  we now thoroughly study the null and non-null asymptotic properties of~$\tilde{\phi}\n$.


\section{Asymptotic properties of the proposed test}
\label{sec:asymptprop}

We first focus on the null hypothesis. Under the single-spike null hypotheses associated with the shape matrices~$\Vb_{0n}$ in~(\ref{V0null}), the test statistic~$T_n(\tilde{\Vb}_{0n})$ of the feasible test~$\tilde{\phi}\n$ is asymptotically chi-square with~$p-1$ degrees of freedom, which easily follows from the asymptotic equivalence result in Theorem~\ref{Theorasymplin} and from the fact that Theorem~\ref{TheorLAN} implies that~$T_n({\Vb}_{0n})$ is asymptotically chi-square with~$p-1$ degrees of freedom under such sequences. Since the latter theorem focuses on single-spike shape matrices, there is no guarantee, however, that this extends to more general shape matrices. The following result shows that the null asymptotic distribution of~$T_n(\tilde{\Vb}_{0n})$ remains asymptotically chi-square with~$p-1$ degrees of freedom for arbitrary sequences of null shape matrices. 

\begin{Theor}
	\label{nullTn}
Let $(\Vb_{0n})$ be a sequence of null shape matrices as in~(\ref{toestimnuis}). Then, under~${\rm P}_{\Vb_{0n}}\n$, $T_n(\tilde{\Vb}_{0n})$ is asymptotically chi-square with~$p-1$ degrees of freedom.
\end{Theor}

We stress that this result does not only allow for general, multi-spike, null shape matrices, but also for weakly identifiable eigenvectors~$\thetab_{n1}$. In particular, it applies in each of the four asymptotic scenarios considered in Theorem~\ref{TheorLAN}. Since there is no guarantee that the asymptotic equivalence in Theorem~\ref{TheorequivTyler} holds under weak identifiability, it is unclear whether or not the Tyler test~$\phi_{\rm Tyl}\n$ is, like~$\tilde{\phi}\n$, robust to weak identifiability. As we will show through simulations in the next section, it 	actually turns out that the Tyler test severely fails to be robust in this sense. 

Now, the nice robustness properties above are not sufficient, on their own, to justify resorting to the proposed sign test, as it might be the case that such robustness is obtained at the expense of power. To see whether or not this is the case, we turn to the investigation of the non-null asymptotic properties of~$\tilde{\phi}\n$. 
We have the following result. 

\begin{Theor} 
\label{Theor3rdLeCam}
Fix $\thetab_0 \in {\cal S}^{p-1}$. Let~$(\delta_n)$ be a positive sequence that is~$O(1)$ and take~$\nu_n=1/(\sqrt{n} \gamma_n)$. Let~$(\taub_n)$ be a sequence converging to~$\taub$ and such that~$\thetab_0+\nu_n\taub_n\in\mathcal{S}^{p-1}$ for any~$n$. Then, under ${\rm P}\n_{{\bf V}_{1n}}$, with~${\bf V}_{1n}$ as in~\eqref{V1alt}, we have the following as $n\to\infty$: %
\begin{enumerate}
\item[(i)] if $\delta_n \equiv 1$,  then~$T_n(\tilde{\Vb}_{0n})$ is asymptotically chi-square with $p-1$ degrees of freedom and non-centrality parameter
$$
\frac{p(p+(p-1)\xi)}{(p+2) (p-\xi)} \, \| \taub\|^2
;
$$ 
\item[(ii)] if~$\delta_n$ is $o(1)$ with~$\sqrt{n} \delta_n \to \infty$, then $T_n(\tilde{\Vb}_{0n})$ is asymptotically chi-square with $p-1$ degrees of freedom and non-centrality parameter
$$
\frac{p}{p+2}\, \| \taub\|^2
;
$$ 
\item[(iii)] if $\delta_n=1/\sqrt{n}$, then $T_n(\tilde{\Vb}_{0n})$ is asymptotically chi-square with $p-1$ degrees of freedom and non-centrality parameter
$$
\frac{p}{p+2}\,
 \| \taub\|^2 \Big(1- \frac{1}{2\xi^2} \| \taub \|^2 \Big)^2\Big(1- \frac{1}{4\xi^2} \| \taub\|^2\Big)
;
$$ 
\item[(iv)] if $\delta_n=1/\sqrt{n}$, then, under~${\rm P}\n_{\Vb_{n1}}$ with~$\Vb_{n1}$ based on~$\delta_n\equiv 1/\xi$, $T_n(\tilde{\Vb}_{0n})$ remains asymptotically chi-square with~$p-1$ degrees of freedom.
\end{enumerate}
\end{Theor}

From contiguity, Theorem~\ref{Theorasymplin} implies that the test~$\tilde{\phi}\n$ enjoys the same asymptotic behavior as~$\phi\n$ under the contiguous alternatives of Theorem~\ref{TheorLAN}(i)-(ii), hence inherits the Le Cam optimality properties of the latter test in the corresponding regimes. Consequently, the local asymptotic powers associated with the non-null results in Theorem~\ref{Theor3rdLeCam}(i)-(ii) are the maximal ones that can be achieved. In regime~(iii), Theorem~\ref{Theor3rdLeCam} entails that the test~$\tilde{\phi}\n$ is rate-optimal, in the sense that it shows non-trivial asymptotic powers against the corresponding contiguous alternatives (the non-standard nature of the limiting experiment in this regime does not allow stating stronger optimality properties, though). Finally, this test is optimal in regime~(iv), but trivially so since Theorem~\ref{TheorLAN}(iv) implies that the trivial $\alpha$-test is also optimal in this regime.

\begin{figure}[!htbp] 
\vspace*{-0.0cm}
\centering
  \includegraphics[width=\textwidth]{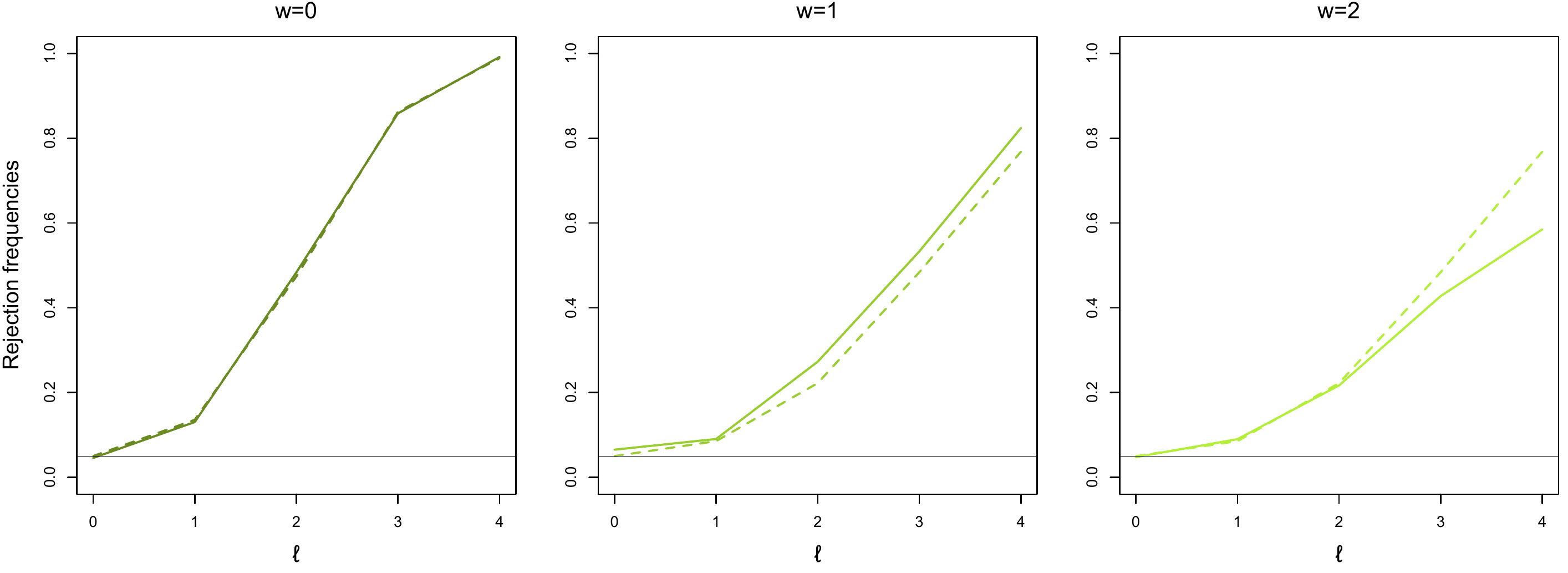}
\vspace*{-0.0cm} 
  \caption{Rejection frequencies (solid curves), under various null and non-null distributions, of the proposed sign test in~(\ref{ProposedTest}) performed at asymptotic level~$\alpha=5\%$; the value~$w=0$ corresponds to the classical asymptotic scenario where~$\thetab_{n1}$ remains asymptotically identifiable, whereas $w=1,2$ provide asymptotic scenarios involving weak identifiability (the lighter the color, the weaker the identifiability). Random samples were drawn from six-variate multinormal distributions; see Section~\ref{sec:asymptprop} for details. Theoretical asymptotic powers are also shown (dashed curves).
  }
  \label{Fig1}
\end{figure}

We performed the following simulation to check the validity of Theorem~\ref{Theor3rdLeCam}. For any combination of~$\ell\in\{0,1,2,3,4\}$ and~$w\in\{0,1,2\}$, we generated $M=2,\!500$ mutually independent random samples~$\Xb_1^{(\ell,w)}\!,\ldots,\Xb_n^{(\ell,w)}$ of size~$n=200,\!000$ from the six-variate ($p=6$) multinormal distribution with mean zero and covariance matrix 
\begin{equation}
	\label{Sigmasimupower}
\Sigb_n^{(\ell,w)} 
:= 
\bigg(
1- \frac{n^{-w/6}}{p}
\bigg) 
{\bf I}_p + 
n^{-w/6} \thetab_{n1}^{(\ell)} \thetab_{n1}^{(\ell)\prime}
,
\end{equation}
with~$\thetab_0:=(1,0,\ldots,0)'\in\R^p$ and $\thetab_{n1}^{(\ell)}:=(\cos\alpha_{n\ell}, \sin \alpha_{n\ell},0,\ldots,0)\pr\in\R^p$, where~$\alpha_{n\ell}:=2\arcsin(\ell\nu_n/2)$  is based on~$\nu_n:=1/(\sqrt{n}\gamma_n)$; here, $\gamma_n$ is as in~(\ref{defingamma}), with the values~$\delta_n=n^{-w/6}$ and~$\xi=1$ that are induced by~(\ref{Sigmasimupower}). The value~$w=0$ yields the standard asymptotic scenario in Theorem~\ref{Theor3rdLeCam}(i), while~$w=1,2$ are associated with weak identifiability situations covered by Theorem~\ref{Theor3rdLeCam}(ii). The value~$\ell=0$ corresponds to the null hypothesis~$\mathcal{H}_0:\thetab_{n1}=\thetab_0$, whereas~$\ell=1,2,3,4$ provide increasingly severe alternatives of the form~$\thetab_{n1}^{(\ell)}=\thetab_0+\nu_n\taub_{n\ell}$, where~$\taub_{n\ell}$ is some $p$-vector with norm~$\ell$; this allows us to obtain the corresponding asymptotic local powers from Theorem~\ref{Theor3rdLeCam}(i)-(ii). For each sample, we performed the proposed sign test for~${\cal H}_0: \thetab_{n1}= \thetab_0$ in~(\ref{ProposedTest}) at nominal level~$5\%$. Clearly, the resulting rejection frequencies, that are provided in Figure~\ref{Fig1}, are in agreement with the theoretical asymptotic powers computed from Theorem~\ref{Theor3rdLeCam}, except maybe for the case~$w=2$. Note, however, that at any finite sample size, empirical power curves will eventually converge to flat power curves at the nominal level~$\alpha$ for weak enough identifiability (this follows from Theorem~\ref{Theor3rdLeCam}(iv)), which explains this small deviation observed for~$w=2$. The sign nature of the proposed test makes it superfluous to also consider non-Gaussian elliptical distributions here.

\section{Finite-sample comparisons with competing tests}
\label{sec:simus}

The objective of this section is to compare the proposed sign test to some competitors. We first focus on empirical size under the null hypothesis. For any~$w\in\{0,1,2,3\}$, we generated $M=5,\!000$ mutually independent random samples~$\Xb_1^{(w)},\ldots,\Xb_n^{(w)}$ of size~$n=400$ from the six-variate ($p=6$) Gaussian distribution with mean zero and covariance matrix 
$$
\Sigb_n^{(w)}
:= 
\bigg(
1- \frac{n^{-w/4}}{p}
\bigg) 
{\bf I}_p + 
n^{-w/4} \thetab_0 \thetab_0\pr
,
$$
with $\thetab_0=(1,0,\ldots,0)'\in\R^p$. As in the simulation conducted at the end of Section~\ref{sec:asymptprop}, the value~$w=0$ is associated with classical situations where~$\thetab_{n1}$ remains asymptotically identifiable whereas the values $w=1,2,3$ provide situations where identifiability of~$\thetab_{n1}$ is weaker and weaker. In each replication, we performed four tests for ${\cal H}_0: \thetab_{n1}= \thetab_0$ at asymptotic level~$5\%$ : the classical Gaussian likelihood ratio test 
from \cite{And63}, the \cite{Tyl87ang} test 
in~(\ref{TylerTest}), the Gaussian test 
from \cite{PRV18}, and the proposed sign test in~(\ref{ProposedTest}). The same exercise was repeated with random samples drawn from multivariate $t$ distributions with~$2$, $4$ and~$6$ degrees of freedom, in each case with mean zero and scatter matrix~$\Sigb_n^{(w)}$.

\begin{figure}[!htbp]
\vspace*{-0.0cm}
\centering
  \includegraphics[width=\textwidth]{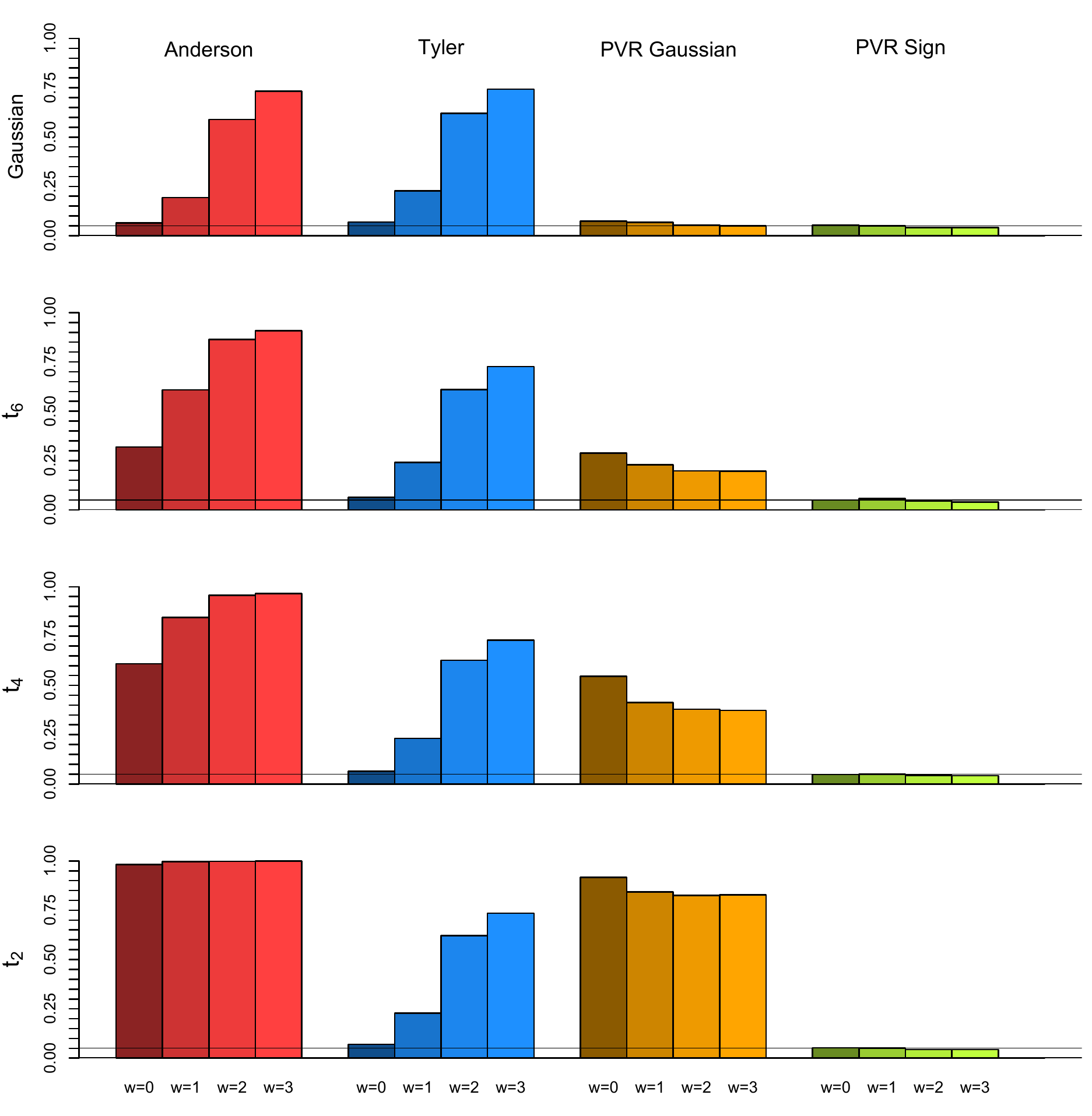}
\vspace*{-0.0cm}
  \caption{Null rejection frequencies, under six-variate Gaussian, $t_6$, $t_4$ and $t_2$ densities, of four tests for the null hypothesis~$\mathcal{H}_0:\thetab_{n1}=(1,0,\ldots,0)'(\in\R^6)$, all performed at asymptotic level~$\alpha=5\%$. The tests considered are the \cite{And63} test, the \cite{Tyl87ang} test, the Gaussian test from \cite{PRV18}, and the proposed sign test; $w=0$ corresponds to the classical asymptotic scenario where~$\thetab_{n1}$ remains asymptotically identifiable, whereas $w=1,2,3$ provide asymptotic scenarios involving weaker and weaker identifiability; see Section~\ref{sec:simus} for details.}
 \label{Fig2}
 \end{figure}

The resulting null rejection frequencies are reported in Figure \ref{Fig2}. Clearly, the \cite{And63} test meets the nominal level constraint only in the Gaussian, well identified, case. As expected, the  \cite{Tyl87ang} test, which is a sign test, can deal with heavy tails, but the results make it clear that this test strongly overrejects the null hypothesis under weak identifiability. The opposite holds for the \cite{PRV18} test, that resists weak identifiability situations in the Gaussian case but cannot deal with heavy tails. In line with the theoretical results of the previous sections, the proposed sign test resists both heavy tails and weak identifiability. 

\begin{figure}[!htbp]
\vspace*{-0.0cm}
\centering
  \includegraphics[width=\textwidth,height=130mm]{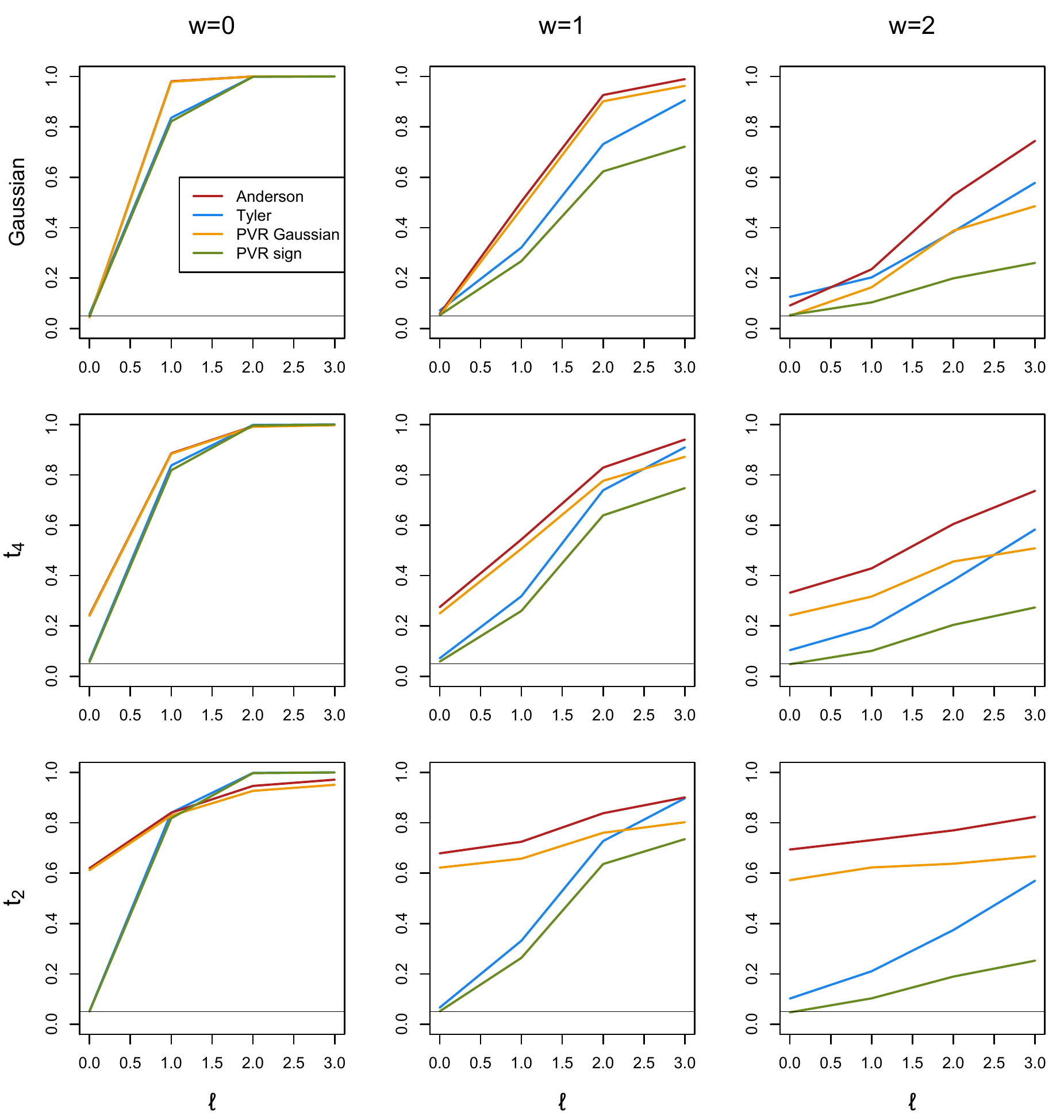}
\vspace*{-0.0cm}
  \caption{Empirical power curves, under bivariate Gaussian (top), $t_4$ densities (middle) and $t_2$ densities (bottom), of four tests for the null hypothesis~$\mathcal{H}_0:\thetab_{n1}=(1,0)'$, all performed at asymptotic level~$\alpha=5\%$. The tests considered are the \cite{And63} test, the \cite{Tyl87ang} test, the Gaussian test from \cite{PRV18}, and the proposed sign test; $w=0$ corresponds to the classical asymptotic scenario where~$\thetab_{n1}$ remains asymptotically identifiable, whereas $w=1$ and $w=2$ provide asymptotic scenarios involving weaker identifiability; see Section~\ref{sec:simus} for details.}
 \label{Fig3} 
 \end{figure}

Although the simulation exercise above shows that the proposed sign test is the only one that meets the asymptotic level constraint under heavy tails and weak identifiability, we now turn to a power comparison of the various tests. For $w\in \{0,1,2\}$, $\ell \in \{0, 1,2,3 \}$ and~$d \in \{1,2,3 \}$, we generated $M=2,\!500$ mutually independent bivariate ($p=2$) random samples~$\Xb_1^{(d,\ell,w)},\ldots,\Xb_n^{(d, \ell,w)}$ of size~$n=200$ with covariance/scatter matrix
$$
\Sigb_n^{(w, \ell)}
:= 
\bigg(
1- \frac{n^{-w/8}}{p}
\bigg) 
{\bf I}_p + 
n^{-w/8} (\thetab_0+ \taub_\ell) (\thetab_0+ \taub_\ell)\pr
,
$$
where $\thetab_0=(1,0)'$ and $\taub_\ell:=(\cos(\ell \pi/12)-1, \sin(\ell \pi/12))\pr$. The $\Xb_i^{(1,\ell,w)}$'s have a Gaussian distribution with covariance matrix $\Sigb_n^{(w, \ell)}$, the $\Xb_i^{(2,\ell,w)}$'s have a student $t_4$ distribution with scatter matrix $\Sigb_n^{(w, \ell)}$, and the $\Xb_i^{(3,\ell,w)}$'s have a student $t_2$ distribution with scatter matrix $\Sigb_n^{(w, \ell)}$. Note that the value $\ell=0$ corresponds to the null hypothesis, while the values $\ell=1,2,3$ provide increasingly severe alternatives. The parameter~$w$ provides different strengths of identifiability, in the same way as in the previous simulation exercise. We performed the same four tests for ${\cal H}_0: \thetab_{n1}= \thetab_0$ as above, still at asymptotic level~$5\%$. 
 The resulting rejection frequencies are plotted in Figure~\ref{Fig3}. Clearly, the results confirm that the proposed sign test is the only test that is robust to both weak identifiability and tail heaviness. It is also seen that this sign test shows power under weakly identified situations and that this power, as expected, does not depend on the tails of the parent distribution. Finally, note that the weaker the identifiability, the larger the sample size needs to be to provide some power, which is quite natural.

 \section{Real data illustration}
 \label{sec:real data}

 We illustrate the practical relevance of the proposed sign test on the famous Swiss banknote dataset, 
 that was also used for illustration in \cite{PRV18}, to which we refer for more details. This dataset is available in the R package \verb|uskewfactors| (\citealp{uskew}) and consists of six measurements on 100 genuine and 100 counterfeit old Swiss 1000-franc banknotes. As in \cite{Flu88} (see pp.~41--43), we restrict here to $n=85$ counterfeit bills made by the same forger and focus on four of the six available measurements: the width~$L$ of the left side of the banknote, the width~$R$ on its right side, the width~$B$ of the bottom margin and the width~$T$ of the top margin, all measured in~mm$\times 10^{-1}$. The resulting sample covariance matrix
has eigenvalues $\hat{\lambda}_1=101.48$, $\hat{\lambda}_2=12.89$, $\hat{\lambda}_3=10.11$ and $\hat{\lambda}_4=2.63$, and corresponding eigenvectors
$$
\label{eigenvectorsrealdata}
\hat{\thetab}_1 
=
\left(
\!
\begin{array}{c} 
.032 \\
 -.012 \\ 
 .820 \\ 
 -.571 
 \end{array}
\!
 \right)
\!
,
\
\hat{\thetab}_2
=
\left(
\!
\begin{array}{cccc} 
 .593  \\
 .797 \\ 
 .057 \\ 
 .097 
 \end{array}
\!
 \right)
\!
,
\
\hat{\thetab}_3
=
\left(
\!
\begin{array}{cccc} 
-.015 \\
-.129 \\ 
.566  \\ 
.814 
 \end{array}
\!
 \right)
\!
,
\textrm{ and }
\hat{\thetab}_4
=
\left(
\!
\begin{array}{cccc} 
 .804 \\
 -.590 \\ 
-.064 \\ 
 -.035 
 \end{array}
\!
 \right)
\!
.
$$
Clearly, the first principal component can be interpreted as the vertical position of the print image on the bill since it is a contrast between~$B$ and~$T$. Similarly, the second principal component could be interpreted as an aggregate of~$L$ and~$R$, that is, as the vertical size of the bill. Yet, Flury refrains from interpreting the second component in this way as the second and third roots are quite close to each other. Accordingly, he reports that the corresponding eigenvectors should be considered spherical.
 
In view of the discussion above, it is natural to test that~$L$ and~$R$ indeed contribute equally to the second component and that no other variables contribute to it. In other words, it is natural to test the null hypothesis~${\cal H}_0: \thetab_2=\thetab_2^0$, with~$\thetab_2^0:=(1,1,0,0)\pr/\sqrt{2}$. While the tests discussed in the present paper address testing problems on the first eigenvector~$\thetab_1$, obvious modifications of these tests allow performing inference on any other eigenvector~$\thetab_j$, $j=2, \ldots,p$. In \cite{PRV18}, the Anderson test~$\phi_{\rm A}$ and the Gaussian test~$\phi_{\rm PVR}$ proposed in that paper were used to test the null hypothesis~${\cal H}_0: \thetab_2=\thetab_2^0$. It is well-known, however, that several observations in the dataset may be considered as outliers (see, e.g., \citealp{Saletal2006}), 
which motivates us resorting to robust tests such as the Tyler  test~$\phi_{\rm Tyl}$ and our sign test~$\tilde{\phi}$. When testing the null  hypothesis~$\mathcal{H}_0$ above, these robust tests provided p-values~$.609$ and~$.992$, respectively, which is to be compared to the p-values~$.099$ and~$.177$, respectively provided by~$\phi_{\rm A}$ and~$\phi_{\rm PVR}$. This shows that, at level~$10\%$, only the Anderson test leads to rejection of the null hypothesis. The robustness of our sign test and the fact that the Anderson test tends to strongly overreject the null hypothesis under weak identifiability should make practitioners confident that non-rejection is the right decision in the present case.

\begin{figure}[htbp!] 
\begin{center} 
\includegraphics[width=80mm]{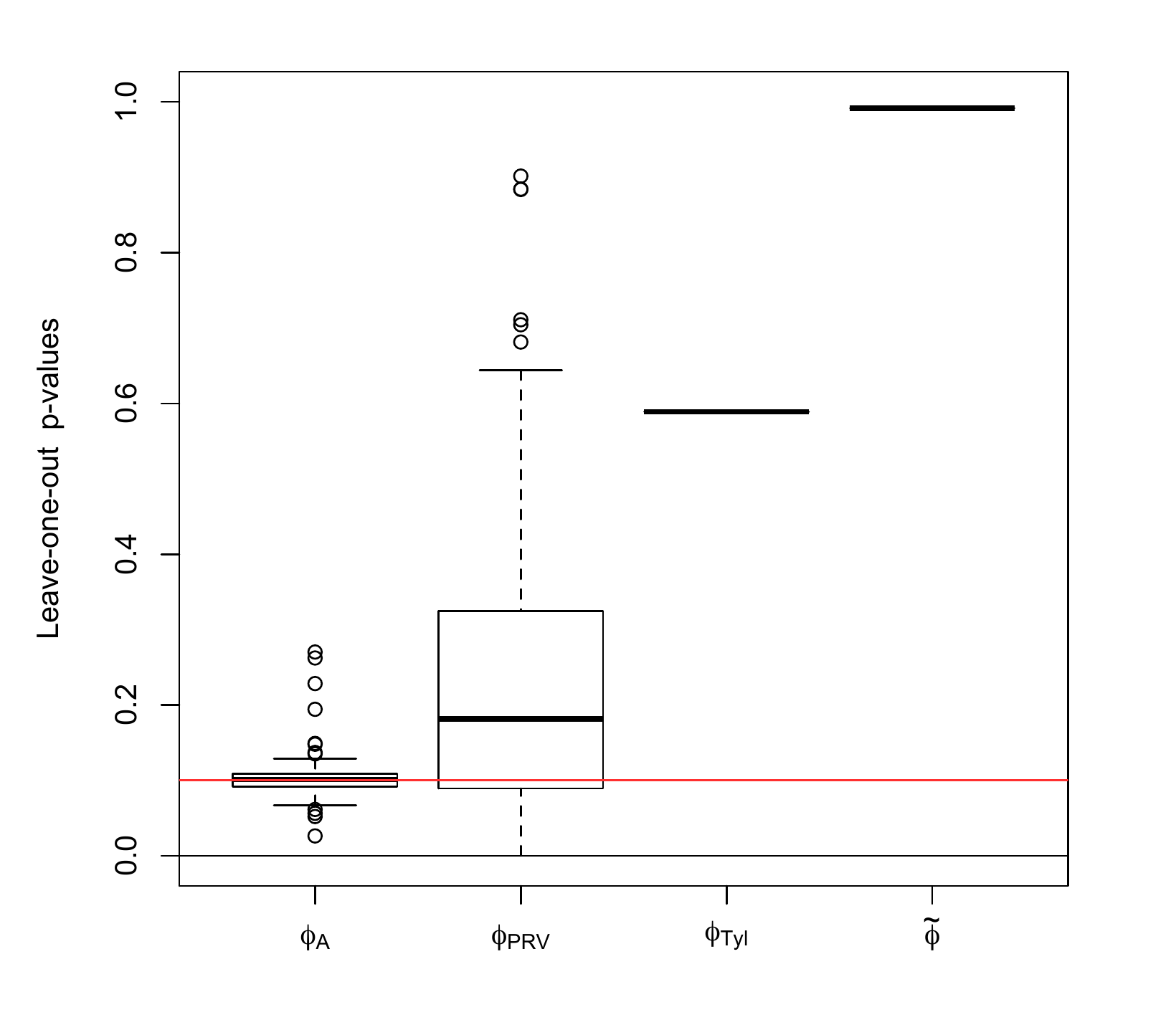}
\vspace{-5mm}
\caption{Boxplots of the 85 ``leave-one-out" $p$-values of the Anderson test~($\phi_{\rm A}$), of the Gaussian \cite{PRV18} test~($\phi_{\rm PVR}$), of the \cite{Tyl87ang} test~($\phi_{\rm Tyl}$), and of the proposed sign test~($\tilde{\phi}$), when testing the null hypothesis~$\mathcal{H}_0:\thetab_2:=(1,1,0,0)\pr/\sqrt{2}$. More precisely, these $p$-values are those obtained when applying the corresponding tests to the 85 subsample of size~$84$ obtained by removing one observation in the real data set considered in the PCA analysis of~\cite{Flu88}, pp.~41--43.}
\label{Fig4}
\end{center}
\end{figure}

To complement the analysis, we performed the same four tests on the~85 subsamples of size~$84$ obtained by removing one observation from the sample considered above. Figure~\ref{Fig4} provides, for each test, a boxplot of the resulting 85 ``leave-one-out" $p$-values. The results show that the Anderson test rejects the null hypothesis much more often than the other tests. It is remarkable that the tests based on spatial signs never led to rejection at any usual nominal level, which, arguably, is due to the natural robustness of spatial signs and of the Tyler estimator of shape.
 
\section{Conclusions and final comments} 
\label{sec:wrapup}

In this work, we considered hypothesis testing for principal directions in challenging asymptotic scenarios involving  weak identifiability. Under ellipticity assumptions, we proposed a sign test that, unlike its competitors, meets the asymptotic level constraint both under heavy tails and under weak identifiability. By resorting to Le Cam's asymptotic theory of statistical experiments, we also proved that this test enjoys strong optimality properties in the class of spatial sign tests (all optimality statements below are relative to this class of tests). In particular, it is locally asymptotically optimal in classical situations where~$\thetab_{n1}$ remains asymptotically identifiable. It follows from our results that the likelihood ratio test from \cite{Tyl87ang} satisfies the same optimality property. Our sign test, however, shows strong advantages over the Tyler test: not only does our test meet the level constraint under \emph{any} weak identifiability situation, but it also remains locally asymptotically optimal in all cases, but for the case~$\delta_n \sim 1/ \sqrt{n}$ for which our test is still rate-optimal. 
 
The following comments are in order. First, our sign test is not only robust to heavy tails and weak identifiability but also to (some) departures from ellipticity. More precisely, it should be clear that our sign test only assumes that the spatial signs~$\Ub_{ni}=\Xb_{ni}/\|\Xb_{ni}\|$ follow an angular Gaussian distribution. As a consequence, the spatial signs do not need be independent of the radii~$\|\Xb_{ni}\|$, which implies in particular that our sign test can deal with some skewed distributions. More precisely, the proposed test only requires that the observations~$\Xb_{ni}$, $i=1,\ldots,n$, form a random sample from a distribution \emph{with  elliptical directions};
see \cite{Ran00}. Second, it has been throughout assumed that the parent elliptical distribution was centered. This was mainly for the sake of readability, as our results can easily be extended to the unspecified location case. More precisely, an unspecified-location version of the proposed test can simply be obtained by replacing the spatial signs~$\Ub_{ni}=\Xb_{ni}/\|\Xb_{ni}\|$ in our sign test with centered versions~$\Ub_{ni}(\hat{\pmb\mu}_n)=(\Xb_{ni}-\hat{\pmb\mu}_n)/\|\Xb_{ni}-\hat{\pmb\mu}_n\|$, where~$\hat{\pmb \mu}_n$ is an arbitrary root-$n$ consistent estimator of the center of the underlying elliptical distribution. A natural choice, that would be root-$n$ consistent even in the large class of distributions with elliptical directions, is the affine-equivariant median from \cite{HetRan2002}. Due to the (Fisher) orthogonality between location and scatter parameters under ellipticity (see \citealp{HP06}), all asymptotic results of this paper readily extend to the resulting unspecified-location sign test.

While this work provides an overall good procedure to test for principal directions under weak identifiability, it also opens perspectives for future research. As mentioned above, the optimality of the proposed test is relative to the class of spatial sign tests. Restricting to sign tests of course is a guarantee for excellent robustness properties, yet it might be so that, if some slightly lower robustness is also acceptable, then higher asymptotic efficiency could be achieved. In particular, it should be possible to develop signed rank tests that provide a nice trade-off between efficiency and robustness. It is expected that such tests can deal with heavy tails and are robust to weak identifiability, while uniformly dominating, in terms of asymptotic relative efficiencies, parametric Gaussian tests in classical cases where the leading principal direction remains asymptotically identifiable; see \cite{Pai2006}.


\appendix

\section{Proof of Theorem~\ref{TheorLAN}} 

We start with some preliminary lemmas.

\begin{Lem} \label{LemA1}
Consider the shape matrices~$\Vb_{0n}$ and~$\Vb_{n1}$ in \eqref{V0null} and \eqref{V1alt}. Then, (i)~$\Vb_{0n}$ and~$\Vb_{n1}$ share the same determinant; (ii) for any real number~$a$, the $a$th matrix powers of~$\Vb_{0n}$ and~$\Vb_{n1}$ are given by 
\begin{eqnarray*}
	{\bf V}_{0n}^{a}
&\!\!\!=\!\!\!&
\Big(1- \frac{\delta_{n}\xi}{p} \Big)^a ({\bf I}_p-\thetab_0 \thetab_0 \pr) + \Big(1 + \frac{(p-1)\delta_{n}\xi}{p} \Big)^a \thetab_0 \thetab_0 \pr
\\[2mm]
&\!\!\!=\!\!\!&
\Big(1- \frac{\delta_{n}\xi}{p} \Big)^a {\bf I}_p + \lambda_{a,n} \thetab_0 \thetab_0 \pr
\end{eqnarray*}
and
\begin{eqnarray*}
{\bf V}_{1n}^{a}
&\!\!\!=\!\!\!&
\Big(1- \frac{\delta_{n}\xi}{p} \Big)^a ({\bf I}_p-(\thetab_0+ \nu_n \taub_n)(\thetab_0+ \nu_n \taub_n)\pr) + \lambda_{a,n} (\thetab_0+ \nu_n \taub_n)(\thetab_0+ \nu_n \taub_n)\pr
\\[2mm]
&\!\!\!=\!\!\!&
\Big(1- \frac{\delta_{n}\xi}{p} \Big)^a {\bf I}_p + \lambda_{a,n} (\thetab_0+ \nu_n \taub_n)(\thetab_0+ \nu_n \taub_n)\pr
,
\end{eqnarray*}
where we let $\lambda_{a,n}:=(1+(p-1)\delta_{n}\xi/p)^{a}-(1-\delta_{n}\xi/p)^{a}$.  
\end{Lem}

\noindent {\bf Proof of Lemma \ref{LemA1}.} 
(i) Letting~$\thetab_n:=\thetab_0+ \nu_n \taub_n$, $r:=1-\delta_{n}\xi/p$ and~$s:=1+(p-1)\delta_{n}\xi/p$, rewrite~$\Vb_{0n}$ and~$\Vb_{n1}$ as 
\begin{equation}
	\label{todescomp}
{\bf V}_{0n} = r ({\bf I}_p-\thetab_0 \thetab_0\pr) + s \thetab_0 \thetab_0\pr 
\quad
\textrm{ and } 
\quad
{\bf V}_{1n} = r ({\bf I}_p-\thetab_n\thetab_n\pr) + s \thetab_n\thetab_n\pr 
.
\end{equation}
Both these matrices have eigenvalues~$r$ with multiplicity~$p-1$ and~$s$ with multiplicity one, hence have determinant~$r^{p-1}s$. (ii) The result directly follows from the spectral decompositions in~(\ref{todescomp}). 
\cqfd 
\vspace{3mm}


\begin{Lem} \label{LemA2}  If~$\delta_n\equiv 1$, then 
$\lim_{\ny} p\delta_{n}\xi/\gamma_n=p+(p-1)\xi$, where $\gamma_n$ is defined in \eqref{defingamma}. If $\delta_n=o(1)$, then 
$\lim_{\ny} p\delta_{n}\xi/\gamma_n=p$.
\end{Lem}

\noindent {\bf Proof of Lemma \ref{LemA2}}. Since~$p\delta_{n}\xi/\gamma_n=p+(p-1) \delta_{n}\xi$,
the result follows. 
\cqfd 
\vspace{3mm}


\noindent To state the next result, we first introduce some notation. Denoting as~${\bf e}_{\ell}$ the $\ell$th vector of the canonical basis of $\R^p$ and by ${\bf A} \otimes {\bf B}$ the Kronecker product between the matrices~${\bf A}$ and ${\bf B}$, we let ${\bf K}_p:=\sum_{i,j =1}^{p}({\bf e}_{i}{\bf e}_{j}\pr)\otimes ({\bf e}_{j}{\bf e}_{i}\pr)$ stand for the $p^2 \times p^2$ \emph{commutation matrix} and define ${\bf J}_p:= \sum_{i,j =1}^{p}({\bf e}_{i}{\bf e}_{j}\pr)\otimes ({\bf e}_{i}{\bf e}_{j}\pr)=({\rm vec}\,{\bf I}_p)({\rm vec}\,{\bf I}_p)\pr$. Further let
$$
\Tb_{n}(\Vb)
:=
\frac{1}{n} \sum_{i=1}^n \, 
{\rm vec}\big( \Ub_{ni}(\Vb) \Ub_{ni}\pr(\Vb) \big) 
\big(   
{\rm vec}\big(  \Ub_{ni}(\Vb) \Ub_{ni}\pr(\Vb) \big)
\big)\pr
$$
with~$\Ub_{ni}(\Vb):= \Vb^{-1/2} \Ub_{ni}/\| \Vb^{-1/2} \Ub_{ni} \|$; note that with this notation, $\Sb_n(\Vb)$ in \eqref{defingamma} rewrites 
$$
\Sb_n(\Vb):= \frac{1}{n} \sum_{i=1}^n  \Ub_{ni}(\Vb) \Ub_{ni}\pr(\Vb)
.
$$
We then have the following result.

\begin{Lem} \label{LemA3}  Fix an arbitrary sequence of shape matrices $(\Vb_n)$. Then,  under ${\rm P}\n_{\Vb_n}$, 
$$
(i) 
\qquad 
\Tb_{n}(\Vb_n)
=
\frac{1}{p(p+2)} ({\bf I}_p+ {\bf K}_p+{\bf J}_p)
+
o_{\rm P}(1)
$$
 and 
$$
\hspace{0mm} 
(ii)
\quad
\sqrt{n}
\,
{\rm vec}\big(\Sb_n(\Vb_n)-  {\textstyle{\frac{1}{p}}} {\bf I}_p\big) 
\stackrel{\mathcal{D}}{\to} 
{\cal N}\bigg( {\bf 0} , \frac{1}{p(p+2)} ({\bf I}_p+ {\bf K}_p+{\bf J}_p)- \frac{1}{p^2}{\bf J}_p \bigg)
$$
as $\ny$.
\end{Lem}

\noindent  {\bf Proof of Lemma \ref{LemA3}}. For any $n$, $\Ub_{n1}(\Vb_n),\ldots,\Ub_{nn}(\Vb_n)$ form a random sample from the uniform distribution on~${\cal S}^{p-1}$. Therefore, the result follows from~(i) the weak law of large numbers and from~(ii) the central limit theorem, by using in both cases Lemma~\mbox{A.2} in \cite{PaiVer2016}. 
\cqfd 
\vspace{3mm}


\noindent {\bf Proof of Theorem~\ref{TheorLAN}.} First note that the quantity~$\gamma_n$ in~(\ref{defingamma}) satisfies 
\begin{equation} 
\label{eqgamma}
-\gamma_n=\Big(1- \frac{\delta_{n}\xi}{p} \Big) \lambda_{-1,n}
,
\end{equation}
where~$\lambda_{-1,n}$ was defined in Lemma~\ref{LemA1}. Part~(ii) of this lemma therefore yields
\begin{equation} 
\label{quadH0}
\Ub_{ni}\pr \Vb_{0n}^{-1}\Ub_{ni} 
= \Big(1- \frac{\delta_{n}\xi}{p} \Big)^{-1}+ \lambda_{-1,n} (\Ub_{ni}\pr \thetab_0)^2 
=
  \Big(1- \frac{\delta_{n}\xi}{p} \Big)^{-1} 
  \big(1-   \gamma_n (\Ub_{ni}\pr \thetab_0)^2\big)
\end{equation}
and, similarly, 
$$
\Ub_{ni}\pr \Vb_{n1}^{-1}\Ub_{ni} =\Big(1- \frac{\delta_{n}\xi}{p} \Big)^{-1} 
\big(1-   \gamma_n (\Ub_{ni}\pr (\thetab_0+ \nu_n \taub_n))^2\big).
$$
Recalling the angular density in~(\ref{angdens}) and using Lemma~\ref{LemA1}(i), we then obtain
\begin{eqnarray*} 
\Lambda_n 
&\!\!\!=\!\!\!& 
-\frac{p}{2} \sum_{i=1}^n \, 
\Big\{
\log (\Ub_{ni}\pr \Vb_{1n}^{-1}\Ub_{ni})- \log(\Ub_{ni}\pr \Vb_{0n}^{-1}\Ub_{ni})
\Big\}
 \nonumber\\
&\!\!\!=\!\!\!&  
-\frac{p}{2} \sum_{i=1}^n \, 
\Big\{
\log (1- \gamma_n (\Ub_{ni}\pr (\thetab_0+ \nu_n \taub_n))^2)-\log (1- \gamma_n (\Ub_{ni}\pr \thetab_0)^2) 
\Big\}
,
\end{eqnarray*}
which, by writing~$\gamma_n (\Ub_{ni}\pr (\thetab_0+ \nu_n \taub_n))^2=\gamma_n (\Ub_{ni}\pr \thetab_0)^2+ 2 \gamma_n \nu_n (\Ub_{ni}\pr \thetab_0)(\Ub_{ni}\pr \taub_n)+ \gamma_n \nu_n^2 (\Ub_{ni}\pr \taub_n)^2$, yields
\begin{eqnarray}
\Lambda_n
&\!\!\!=\!\!\!&  -\frac{p}{2} \sum_{i=1}^n \, \log\bigg(1- \frac{ 2 \gamma_n \nu_n (\Ub_{ni}\pr \thetab_0)(\Ub_{ni}\pr \taub_n)+ \gamma_n \nu_n^2 (\Ub_{ni}\pr \taub_n)^2}{1-\gamma_n (\Ub_{ni}\pr \thetab_0)^2} \bigg) \nonumber \\
&\!\!\!=\!\!\!& -\frac{p}{2} \sum_{i=1}^n \, \log\bigg(1- \gamma_n \nu_n \frac{ 2(\Ub_{ni}\pr \thetab_0)(\Ub_{ni}\pr \taub_n)+ \nu_n (\Ub_{ni}\pr \taub_n)^2}{1-\gamma_n (\Ub_{ni}\pr \thetab_0)^2} \bigg)  \nonumber \\
&\!\!\!=:\!\!\!& -\frac{p}{2} \sum_{i=1}^n \, \log\big(1- \gamma_n \nu_n  R_{ni} \big)  
.
\label{arr1}
\end{eqnarray}
A Taylor expansion yields
\begin{eqnarray*}
\Lambda_n
&\!\!\!=\!\!\!&
\bigg(
\frac{p\gamma_n\nu_n}{2}
\sum_{i=1}^n 
 R_{ni}
\bigg)
+
\bigg(
\frac{p\gamma_n^2\nu_n^2}{4}
\sum_{i=1}^n 
 R_{ni}^2
\bigg)
+
\bigg(
\frac{p\gamma_n^3\nu_n^3}{6}
\sum_{i=1}^n 
 \frac{R_{ni}^3}{(1-\gamma_n\nu_n H_{ni})^3}
\bigg)
\\[2mm]
&\!\!\!=:\!\!\!&
L_{n1}+L_{n2}+L_{n3}
,
\end{eqnarray*}
for some~$H_{ni}$ between~$0$ and~$R_{ni}$. Note that, in all cases~(i)-(iv) considered in the theorem, we have that~$\nu_n\gamma_n=O(1/\sqrt{n})$, $\nu_n=O(1)$ and that there exists~$\eta\in(0,1)$ such that~$|\gamma_n|<1-\eta$ (recall that, in case~(i), $\xi$ is chosen in such a way that~$\Vb_{0n}$ is positive definite). Consequently, the boundedness of the sequence~$(\taub_n)$, along with the fact that~$\|\Ub_{ni}\|= 1$ almost surely, ensures that there exists a positive constant~$C$ such that~$|L_{n3}| \leq C/\sqrt{n}$ almost surely, so that~$L_{n3}$ is~$o_{\rm P}(1)$ (all stochastic convergences in this proof are as~$n \to\infty$ under ${\rm P}\n_{\Vb_{0n}}$). 

Now, using~(\ref{quadH0}), rewrite~$L_{n1}$ as 
\begin{eqnarray} 
 L_{n1} 
&\!\!\!=\!\!\!&
\frac{p\gamma_n\nu_n}{2}
\sum_{i=1}^n \,
\frac{ 2 (\Ub_{ni}\pr \thetab_0)(\Ub_{ni}\pr \taub_n)+ \nu_n (\Ub_{ni}\pr \taub_n)^2}{1-\gamma_n (\Ub_{ni}\pr \thetab_0)^2}
 \nonumber 
 \\[1mm]
&\!\!\!=\!\!\!&
 p\gamma_n\nu_n
\sum_{i=1}^n \, 
\frac{\taub_n\pr \Ub_{ni}\Ub_{ni}\pr\big(\thetab_0+ {\textstyle{\frac{1}{2}}} \nu_n \taub_n\big)}{1-\gamma_n (\Ub_{ni}\pr \thetab_0)^2}
 \nonumber 
 \\[1mm]
&\!\!\!=\!\!\!&
 p\gamma_n\nu_n
\Big(1- \frac{\delta_{n}\xi}{p} \Big)^{-1}
\sum_{i=1}^n  \,
\frac{\taub_n\pr \Ub_{ni}\Ub_{ni}\pr\big(\thetab_0+ {\textstyle{\frac{1}{2}}} \nu_n \taub_n\big)}{\Ub_{ni}\pr \Vb_{0n}^{-1}\Ub_{ni}}
 \nonumber 
 \\[1mm]
&\!\!\!=\!\!\!&
 \frac{p^2n\gamma_n\nu_n}{p-\delta_{n}\xi}
\,
 \taub_n\pr 
  {\bf V}_{0n}^{1/2} 
\Sb_n(\Vb_{0n})
 {\bf V}_{0n}^{1/2} 
 \big(\thetab_0+ {\textstyle{\frac{1}{2}}} \nu_n \taub_n\big) 
 \nonumber 
.
\end{eqnarray}
Note that, using~(\ref{pert}), we have
\begin{eqnarray*}
\lefteqn{
\taub_n\pr 
  {\bf V}_{0n}
\big(\thetab_0+ {\textstyle{\frac{1}{2}}} \nu_n \taub_n\big)
=
\Big(
1- \frac{\delta_{n}\xi}{p}
\Big) 
\taub_n\pr\big(\thetab_0+ {\textstyle{\frac{1}{2}}} \nu_n \taub_n\big)
  + 
  \delta_{n}\xi (\taub_n\pr \thetab_0) \thetab_0\pr \big(\thetab_0+ {\textstyle{\frac{1}{2}}} \nu_n \taub_n\big)
}
\\[1mm]
& & 
\hspace{10mm} 
=
  \delta_{n}\xi (\taub_n\pr \thetab_0)  
  \big( 1 + {\textstyle{\frac{1}{2}}} \nu_n \thetab_0\pr\taub_n \big)
=
-\frac{1}{2}
\nu_n
  \delta_{n}\xi \|\taub_n\|^2  
  \big( 1 + {\textstyle{\frac{1}{2}}} \nu_n \thetab_0\pr\taub_n \big)
 ,
 \end{eqnarray*}
which yields
\begin{eqnarray} 
 L_{n1} 
&\!\!\!=\!\!\!&
 \frac{p^2n\gamma_n\nu_n}{p-\delta_{n}\xi}
\,
 \taub_n\pr 
  {\bf V}_{0n}^{1/2} 
\big( \Sb_n(\Vb_{0n}) - {\textstyle{\frac{1}{p}}} \mathbf{I}_p \big)
 {\bf V}_{0n}^{1/2} 
 \big(\thetab_0+ {\textstyle{\frac{1}{2}}} \nu_n \taub_n\big) 
 \nonumber 
+
 \frac{p n\gamma_n\nu_n}{p-\delta_{n}\xi}
\,
 \taub_n\pr 
  {\bf V}_{0n}
 \big(\thetab_0+ {\textstyle{\frac{1}{2}}} \nu_n \taub_n\big) 
 \nonumber 
\\[2mm]
&\!\!\!=\!\!\!&
 \frac{p^2n\gamma_n\nu_n}{p-\delta_{n}\xi}
\,
 \taub_n\pr 
  {\bf V}_{0n}^{1/2} 
\big( \Sb_n(\Vb_{0n}) - {\textstyle{\frac{1}{p}}} \mathbf{I}_p \big)
 {\bf V}_{0n}^{1/2} 
 \big(\thetab_0+ {\textstyle{\frac{1}{2}}} \nu_n \taub_n\big) 
 \nonumber 
-
 \frac{p n\gamma_n\nu_n^2\delta_{n}\xi }{2(p-\delta_{n}\xi)}
\,
  \|\taub_n\|^2  
  \big( 1 + {\textstyle{\frac{1}{2}}} \nu_n \thetab_0\pr\taub_n \big)
\label{L1nexpress}
,
\end{eqnarray}
where, from Lemma~\ref{LemA1}, we have
\begin{equation}
\label{V0undemitau}	
{\bf V}_{0n}^{1/2} \thetab_0= \Big(1 + \frac{(p-1)\delta_{n}\xi}{p} \Big)^{1/2} \thetab_0
\ \textrm{ and } \
{\bf V}_{0n}^{1/2} \taub_n
= 
\Big(1 - \frac{\delta_{n}\xi}{p} \Big)^{1/2} \taub_n + \lambda_{-1/2,n} (\thetab_0'\taub_n)\thetab_0
.
\end{equation}

Turning to~$L_{n2}$,  
\begin{eqnarray*} 
\lefteqn{
 L_{n2} 
= \frac{p\gamma_n^2\nu_n^2}{4}
\sum_{i=1}^n 
\frac{ \{2(\Ub_{ni}\pr \thetab_0)(\Ub_{ni}\pr \taub_n)+ \nu_n (\Ub_{ni}\pr \taub_n)^2\}^2}{(1-\gamma_n (\Ub_{ni}\pr \thetab_0)^2)^2}
}
\\[2mm]
& & 
\hspace{0mm} 
=  
 \frac{p\gamma_n^2\nu_n^2}{4}
  \Bigg\{
    4 (\thetab_0 \otimes \thetab_0)\pr \sum_{i=1}^n\frac{{\rm vec}(\Ub_{ni} \Ub_{ni}\pr) \big({\rm vec}(\Ub_{ni} \Ub_{ni}\pr)\big)\pr}{(1-\gamma_n (\Ub_{ni}\pr \thetab_0)^2)^2}\, (\taub_n \otimes \taub_n) 
\\[2mm]
& & 
\hspace{28mm} 
+ 4  \nu_n (\taub_n \otimes \taub_n)\pr \sum_{i=1}^n\frac{{\rm vec}(\Ub_{ni} \Ub_{ni}\pr) \big({\rm vec}(\Ub_{ni} \Ub_{ni}\pr)\big)\pr}{(1-\gamma_n (\Ub_{ni}\pr \thetab_0)^2)^2}\, (\taub_n \otimes \thetab_0) 
\\[2mm]
& & 
\hspace{28mm} 
+ \nu_n^2  (\taub_n \otimes \taub_n)\pr \sum_{i=1}^n\frac{{\rm vec}(\Ub_{ni} \Ub_{ni}\pr) \big({\rm vec}(\Ub_{ni} \Ub_{ni}\pr)\big)\pr}{(1-\gamma_n (\Ub_{ni}\pr \thetab_0)^2)^2}\, (\taub_n \otimes \taub_n) 
\Bigg\}
.
\end{eqnarray*}

Using~(\ref{quadH0}) again, Lemma~\ref{LemA3}(i), and the fact that~$(\vb\otimes\vb)\pr{\bf K}_p={\bf K}_1(\vb\otimes\vb)\pr=(\vb\otimes\vb)\pr$ for any $p$-vector~$\vb$, this yields
\begin{eqnarray} 
\lefteqn{
 L_{n2} 
=  
 \frac{np\gamma_n^2\nu_n^2}{4}
 \Big(1- \frac{\delta_{n}\xi}{p} \Big)^{-2}
  \Bigg\{
    4 
((\Vb_{0n}^{1/2}\thetab_0) \otimes (\Vb_{0n}^{1/2}\thetab_0))\pr 
    \Tb_n(\Vb_{0n})
((\Vb_{0n}^{1/2}\taub_n) \otimes (\Vb_{0n}^{1/2}\taub_n)) 
}
\nonumber
\\[2mm]
& & 
\hspace{24mm} 
+ 4  \nu_n 
((\Vb_{0n}^{1/2}\taub_n) \otimes (\Vb_{0n}^{1/2}\taub_n))\pr 
    \Tb_n(\Vb_{0n})
((\Vb_{0n}^{1/2}\taub_n) \otimes (\Vb_{0n}^{1/2}\thetab_0)) 
\nonumber
\\[2mm]
& & 
\hspace{24mm} 
+ \nu_n^2  
((\Vb_{0n}^{1/2}\taub_n) \otimes (\Vb_{0n}^{1/2}\taub_n))\pr 
    \Tb_n(\Vb_{0n})
((\Vb_{0n}^{1/2}\taub_n) \otimes (\Vb_{0n}^{1/2}\taub_n)) 
\Bigg\}
\nonumber
\\[2mm]
& & 
\hspace{-3mm} 
=  
 \frac{np^2\gamma_n^2\nu_n^2}{4(p+2)(p-\delta_{n}\xi)^2}
  \Bigg\{
    4 
((\Vb_{0n}^{1/2}\thetab_0) \otimes (\Vb_{0n}^{1/2}\thetab_0))\pr 
(2{\bf I}_p+{\bf J}_p)
((\Vb_{0n}^{1/2}\taub_n) \otimes (\Vb_{0n}^{1/2}\taub_n)) 
\nonumber
\\[1mm]
& & 
\hspace{8mm} 
+ 4  \nu_n 
((\Vb_{0n}^{1/2}\taub_n) \otimes (\Vb_{0n}^{1/2}\taub_n))\pr 
(2{\bf I}_p+{\bf J}_p)
((\Vb_{0n}^{1/2}\taub_n) \otimes (\Vb_{0n}^{1/2}\thetab_0)) 
\label{L2nexpress}
\\[2mm]
& & 
\hspace{8mm} 
+ \nu_n^2  
((\Vb_{0n}^{1/2}\taub_n) \otimes (\Vb_{0n}^{1/2}\taub_n))\pr 
(2{\bf I}_p+{\bf J}_p)
((\Vb_{0n}^{1/2}\taub_n) \otimes (\Vb_{0n}^{1/2}\taub_n))
+o_{\rm P}(1) 
\Bigg\}
.
\nonumber
\end{eqnarray}

We can now consider the cases~(i)--(iv). We start with cases~(i)--(ii) and let~$\zeta$ be equal to one if case~(i) is considered and to zero if case~(ii) is. Since~$\nu_n\gamma_n=1/\sqrt{n}$ in both cases, we have
\begin{eqnarray*} 
\lefteqn{
 L_{n1} 
=
 \frac{p^2\sqrt{n}}{p-\delta_{n}\xi}
\,
 \taub_n\pr 
  {\bf V}_{0n}^{1/2} 
\big( \Sb_n(\Vb_{0n}) - {\textstyle{\frac{1}{p}}} \mathbf{I}_p \big)
 {\bf V}_{0n}^{1/2} 
 \big(\thetab_0+ {\textstyle{\frac{1}{2}}} \nu_n \taub_n\big) 
}
 \nonumber 
\\[2mm]
& & 
\hspace{23mm} 
-
 \frac{p\delta_{n}\xi}{2\gamma_n(p-\delta_{n}\xi)}
\,
  \|\taub_n\|^2  
  \big( 1 + {\textstyle{\frac{1}{2}}} \nu_n \thetab_0\pr\taub_n \big)
\\[2mm]
& & 
\hspace{3mm} 
=
 \frac{p^2\sqrt{n}}{p-\zeta\xi}
\,
 \taub_n\pr 
  {\bf V}_{0n}^{1/2} 
\big( \Sb_n(\Vb_{0n}) - {\textstyle{\frac{1}{p}}} \mathbf{I}_p \big)
 {\bf V}_{0n}^{1/2} 
 \big(\thetab_0+ {\textstyle{\frac{1}{2}}} \nu_n \taub_n\big) 
 \nonumber 
\\[2mm]
& & 
\hspace{23mm} 
-
 \frac{p+\zeta(p-1)\xi}{2(p-\zeta\xi)}
\,
  \|\taub_n\|^2  
 +o_{\rm P}(1) 
,
\end{eqnarray*}
where the last equality follows from Lemma~\ref{LemA3}(ii), Lemma~\ref{LemA2}, and the fact that $\nu_n=o(1)$.

So, using~(\ref{V0undemitau}), Lemma~\ref{LemA3} again, and the fact that~$\thetab_0\pr\taub_n=o(1)$, we obtain
\begin{eqnarray} 
 L_{n1} 
&\!\!\!=\!\!\!&
\frac{p^2\sqrt{n}}{p-\zeta\xi}
\Big(1 - \frac{\zeta\xi}{p} \Big)^{1/2}
\Big(1 + \frac{\zeta(p-1)\xi}{p} \Big)^{1/2}
\,
 \taub_n\pr 
\big( \Sb_n(\Vb_{0n}) - {\textstyle{\frac{1}{p}}} \mathbf{I}_p \big)
\thetab_0
- \frac{p+\zeta(p-1) \xi}{2(p-\zeta\xi)}
\,
 \|\taub_n\|^2 
 +o_{\rm P}(1) 
\nonumber
\\[2mm]
&\!\!\!=\!\!\!&
\frac{p(p + \zeta(p-1)\xi)^{1/2}\sqrt{n}}{(p-\zeta\xi)^{1/2}}
\,
 \taub_n\pr 
\big( \Sb_n(\Vb_{0n}) - {\textstyle{\frac{1}{p}}} \mathbf{I}_p \big)
\thetab_0
- \frac{p+\zeta(p-1) \xi}{2(p-\zeta\xi)}
\,
 \|\taub_n\|^2 
 +o_{\rm P}(1) 
\nonumber
\\[2mm]
&\!\!\!=\!\!\!&
\frac{p(p + \zeta(p-1)\xi)^{1/2}\sqrt{n}}{(p-\zeta\xi)^{1/2}}
\,
 \taub_n\pr 
({\bf I}_{p}- \thetab_0 \thetab_0\pr) 
\big( \Sb_n(\Vb_{0n}) - {\textstyle{\frac{1}{p}}} \mathbf{I}_p \big)
\thetab_0
\nonumber
\\[2mm]
& & 
\hspace{43mm} 
- \frac{p+\zeta(p-1) \xi}{2(p-\zeta\xi)}
\,
\taub_n\pr({\bf I}_{p}- \thetab_0 \thetab_0\pr) \taub_n
 +o_{\rm P}(1)
 . 
\label{toshowia}
\end{eqnarray}

Now, in cases~(i)--(ii), since~$\nu_n=1/(\sqrt{n}\gamma_n)=o(1)$, (\ref{L2nexpress}) becomes
$$
 L_{n2} 
=  
 \frac{p^2}{(p+2)(p-\zeta\xi)^2}
((\Vb_{0n}^{1/2}\thetab_0) \otimes (\Vb_{0n}^{1/2}\thetab_0))\pr 
    (2{\bf I}_p+{\bf J}_p)
((\Vb_{0n}^{1/2}\taub_n) \otimes (\Vb_{0n}^{1/2}\taub_n)) 
+o_{\rm P}(1) 
,
$$
which, by using~(\ref{V0undemitau}), readily yields
\begin{eqnarray}
 L_{n2} 
&\!\!\!=\!\!\!&  
 \frac{p^2}{(p+2)(p-\zeta\xi)^2}
\Big(1 + \frac{(p-1)\delta_{n}\xi}{p} \Big)
\Big(1 - \frac{\delta_{n}\xi}{p} \Big)
(\thetab_0 \otimes \thetab_0)\pr 
(2{\bf I}_p+{\bf J}_p)
(\taub_n \otimes \taub_n) 
+o_{\rm P}(1) 
\nonumber
\\[2mm]
&\!\!\!=\!\!\!&  
 \frac{p + \zeta(p-1)\xi}{(p+2)(p-\zeta\xi)}
\Big\{
2(\thetab_0\pr\taub_n)^2 + \|\taub_n\|^2
\Big\}
+o_{\rm P}(1) 
\nonumber
\\[2mm]
&\!\!\!=\!\!\!&  
 \frac{p + \zeta(p-1)\xi}{(p+2)(p-\zeta\xi)}\,
 \taub_n\pr({\bf I}_{p}- \thetab_0 \thetab_0\pr) \taub_n
+o_{\rm P}(1) 
.
\label{toshowib}
\end{eqnarray}
The results in~(\ref{LAQi}) and~(\ref{LAQii}) then follow from~(\ref{toshowia})-(\ref{toshowib}), whereas the asymptotic normality results for~$\Deltab_{\thetab_0}^{(i)}$ and~$\Deltab_{\thetab_0}^{(ii)}$ are direct corollaries of Lemma~\ref{LemA3}(ii). 
\vspace{3mm}
 
We turn to case~(iii), for which~$\delta_n=1/\sqrt{n}$ and~$\nu_n=1/(\sqrt{n}\gamma_n)=1/\xi+o(1)$. By using Lemma~\ref{LemA2}, (\ref{L1nexpress}) here becomes
\begin{eqnarray} 
 L_{n1} 
&\!\!\!=\!\!\!&
 \frac{p}{p-(\xi/\sqrt{n})}
\,
 \taub_n\pr 
  {\bf V}_{0n}^{1/2} 
\Umb_{\thetab_0}\n
 {\bf V}_{0n}^{1/2} 
 \big(\thetab_0+ {\textstyle{\frac{1}{2}}} \nu_n \taub_n\big) 
-
 \frac{p\delta_{n}\xi}{2\gamma_n(p-(\xi/\sqrt{n}))}
\,
  \|\taub_n\|^2  
  \big( 1 + {\textstyle{\frac{1}{2}}} \nu_n \thetab_0\pr\taub_n \big)
\nonumber
\\[2mm]
&\!\!\!=\!\!\!&
 \taub_n\pr 
  {\bf V}_{0n}^{1/2} 
\Umb_{\thetab_0}\n
 {\bf V}_{0n}^{1/2} 
 \big(\thetab_0+ {\textstyle{\frac{1}{2\xi}}} \taub_n\big) 
 -
 \frac{1}{2}
\,
  \|\taub_n\|^2  
  \big( 1 - {\textstyle{\frac{1}{4\xi^2}}} \|\taub_n\|^2 \big)
  +
  o_{\rm P}(1)
.
\nonumber
\end{eqnarray}
Since~${\bf V}_{0n}=\mathbf{I}_p+o(1)$, this yields
\begin{equation}
\label{L1nexpressiii}	
 L_{n1} 
=
\taub_n\pr
\Umb_{\thetab_0}\n
\thetab_0
+
\frac{1}{2\xi} 
\,
\taub_n\pr
\Umb_{\thetab_0}\n
\taub_n
 -
 \frac{1}{2}
\,
  \|\taub_n\|^2  
 +
 \frac{1}{8\xi^2}
\,
  \|\taub_n\|^4
  +
  o_{\rm P}(1)
  .
\end{equation}
Turning to~$L_{n2}$, (\ref{L2nexpress}) provides
\begin{eqnarray} 
\lefteqn{
\hspace{3mm} 
 L_{n2} 
=  
 \frac{p^2}{4(p+2)(p+o(1))^2}
  \Big\{
    4 
(\thetab_0 \otimes \thetab_0)\pr 
(2{\bf I}_p+{\bf J}_p)
(\taub_n \otimes \taub_n)
}
\nonumber
\\[1mm]
& & 
\hspace{-3mm} 
+ 4  
\xi^{-1}  
(\taub_n \otimes \taub_n)\pr 
(2{\bf I}_p+{\bf J}_p)
(\taub_n \otimes \thetab_0)
+ 
\xi^{-2}  
(\taub_n \otimes \taub_n)
\pr 
(2{\bf I}_p+{\bf J}_p)
(\taub_n \otimes \taub_n)
\Big\}
+
o_{\rm P}(1)
\nonumber
\\[2mm]
& & 
\hspace{3mm} 
=  
 \frac{1}{4(p+2)}
  \Big\{
   8 (\thetab_0\pr\taub_n)^2+ 4 \|\taub_n\|^2
+ 
  12 \xi^{-1} (\thetab_0\pr\taub_n) \|\taub_n\|^2
+ 
 3  \xi^{-2} \|\taub_n\|^4
\Big\}
+
o_{\rm P}(1)
.
\nonumber
\\[2mm]
& & 
\hspace{3mm} 
=  
 \frac{1}{p+2}
 \|\taub_n\|^2
- 
 \frac{1}{4(p+2)\xi^2}
 \|\taub_n\|^4
+
o_{\rm P}(1)
.
\label{L2nexpressiii}
\end{eqnarray}
Hence, from~(\ref{L1nexpressiii})--(\ref{L2nexpressiii}), we conclude that 
$$
\Lambda_n
=
\taub_n\pr
\Umb_{\thetab_0}\n
\thetab_0
+
\frac{1}{2\xi} 
\,
\taub_n\pr
\Umb_{\thetab_0}\n
\taub_n
 -
 \frac{p}{2(p+2)}
\,
  \|\taub_n\|^2  
+ 
 \frac{p}{8(p+2) \xi^2}
 \|\taub_n\|^4
+
o_{\rm P}(1)
,
$$
as was to be shown. Again, the asymptotic normality result for~${\rm vec}(\Umb_{\thetab_0}\n)$ easily follows from Lemma~\ref{LemA3}(ii).

Finally, we consider case~(iv), under which~$\delta_n=o(1/\sqrt{n})$ and~$\nu_n=O(1)$. It directly follows from~(\ref{L1nexpress}) and~(\ref{L2nexpress}) that~$L_{n1}$ and~$L_{n2}$ are then~$o_{\rm P}(1)$, so that~$\Lambda_n$ also is.
\cqfd
\vspace{3mm}


\section{Proofs of Theorems~\ref{Theorasymplin}--\ref{TheorequivTyler}}

The proofs of this section require the following preliminary result.

\begin{Lem} 
\label{LemB1}
Let $(\Vb_{0n})$ be a sequence of null shape matrices as in~(\ref{toestimnuis}) and denote Tyler's M-estimator of scatter as~$\hat{\Vb}_n$. Then, we have the following under~${\rm P}\n_{\Vb_{0n}}$ as~$n\to\infty$:
(i) letting~$\Gb_p:={\bf I}_{p^2}-{\textstyle{\frac{1}{p+2}}}({\bf I}_{p^2}+{\bf K}_p-{\bf J}_p)$, 
$$
\Gb_p
\big(\Vb_{0n}^{-1/2} \otimes \Vb_{0n}^{-1/2}\big)
\sqrt{n}\, {\rm vec}(\hat{\Vb}_n - \Vb_{0n})
=
p\sqrt{n}\,{\rm vec}\big(\Sb_n(\Vb_{0n})-{\textstyle{\frac{1}{p}}}{\bf I}_p\big)
+
o_{\rm P}(1)
;
$$
(ii) $\sqrt{n}(\hat{\Vb}_n- \Vb_{0n})$ is $O_{\rm P}(1)$.
\end{Lem}

\noindent {\bf Proof of Lemma~\ref{LemB1}}. 
Part~(i) of the lemma follows from~(3.7)--(3.8) in \cite{Tyl1987} and Lemma~\ref{LemA3}(i), whereas Part~(ii) follows from Part~(i) and Lemma~\ref{LemA3}(ii).
\cqfd 
\vspace{3mm}


In the proofs of this section, all stochastic convergences will be as $\ny$ under ${\rm P}_{\Vb_{0n}}\n$. 
\vspace{1mm}

\noindent {\bf Proof of Theorem~\ref{Theorasymplin}}. 
Standard properties of the vec operator provide
\begin{eqnarray} 
\lefteqn{
\frac{1}{\sqrt{n}}
\sum_{i=1}^n
{\rm vec} 
\bigg( 
\frac{\tilde{\Vb}_{0n}^{-1/2}\Ub_{ni} \Ub_{ni}\pr\tilde{\Vb}_{0n}^{-1/2}}{\Ub_{ni}\pr \tilde{\Vb}_{0n}^{-1} \Ub_{ni}}
-
 \frac{\tilde{\Vb}_{0n}^{-1/2}\Ub_{ni} \Ub_{ni}\pr\tilde{\Vb}_{0n}^{-1/2}}{\Ub_{ni}\pr  \Vb_{0n}^{-1} \Ub_{ni}} 
\bigg) 
}
 \nonumber 
 \\[2mm]
& & 
\hspace{-3mm} 
=
-
\frac{1}{\sqrt{n}}
 \sum_{i=1}^n 
\frac{\Ub_{ni}\pr \tilde{\Vb}_{0n}^{-1} \Ub_{ni}-\Ub_{ni}\pr \Vb_{0n}^{-1} \Ub_{ni}}{(\Ub_{ni}\pr \tilde{\Vb}_{0n}^{-1} \Ub_{ni})(\Ub_{ni}\pr  \Vb_{0n}^{-1} \Ub_{ni})} 
\,
{\rm vec} (\tilde{\Vb}_{0n}^{-1/2}{\Ub_{ni} \Ub_{ni}\pr}\tilde{\Vb}_{0n}^{-1/2})
 \nonumber 
 \\[2mm]
& & 
\hspace{-3mm} 
=
-\bigg(
\frac{1}{n}
\sum_{i=1}^n 
\frac{{\rm vec} (\tilde{\Vb}_{0n}^{-1/2} \Ub_{ni} \Ub_{ni}\pr\tilde{\Vb}_{0n}^{-1/2}) 
({\rm vec}(\Ub_{ni} \Ub_{ni}\pr))\pr }{(\Ub_{ni}\pr \tilde{\Vb}_{0n}^{-1} \Ub_{ni})(\Ub_{ni}\pr  \Vb_{0n}^{-1} \Ub_{ni})}
\bigg)
\sqrt{n}\, {\rm vec}\big(\tilde{\Vb}_{0n}^{-1}- \Vb_{0n}^{-1}\big)
\nonumber 
\\[3mm]
& & 
\hspace{-3mm} 
=
-(\Mb_{1n}+\Mb_{2n}) 
\sqrt{n}\, {\rm vec}\big(\tilde{\Vb}_{0n}^{-1}- \Vb_{0n}^{-1}\big)
,
\nonumber
\end{eqnarray}
where
\begin{eqnarray*}
\Mb_{1n}
&\!\!\!:=\!\!\!&
\frac{1}{n}
\sum_{i=1}^n 
\frac{{\rm vec} (\tilde{\Vb}_{0n}^{-1/2} \Ub_{ni} \Ub_{ni}\pr\tilde{\Vb}_{0n}^{-1/2}) 
({\rm vec}(\Ub_{ni} \Ub_{ni}\pr))\pr }{(\Ub_{ni}\pr  \Vb_{0n}^{-1} \Ub_{ni})^2}
\nonumber 
\\[2mm]
&\!\!\!=\!\!\!&
\big( (\tilde{\Vb}_{0n}^{-1/2}\Vb_{0n}^{1/2}) \otimes  (\tilde{\Vb}_{0n}^{-1/2}\Vb_{0n}^{1/2})\big)  
\Tb_{n}(\Vb_{0n}) 
\big( \Vb_{0n}^{1/2} \otimes  \Vb_{0n}^{1/2}\big) 
\nonumber 
\end{eqnarray*}
and
\begin{eqnarray*}
\Mb_{2n}
&\!\!\!:=\!\!\!&
\frac{1}{n}
\sum_{i=1}^n 
\bigg\{
\frac{1}{\Ub_{ni}\pr \tilde{\Vb}_{0n}^{-1} \Ub_{ni}}
-
\frac{1}{\Ub_{ni}\pr {\Vb}_{0n}^{-1} \Ub_{ni}}
\bigg\}
\frac{{\rm vec} (\tilde{\Vb}_{0n}^{-1/2} \Ub_{ni} \Ub_{ni}\pr\tilde{\Vb}_{0n}^{-1/2}) 
({\rm vec}(\Ub_{ni} \Ub_{ni}\pr))\pr }{\Ub_{ni}\pr  \Vb_{0n}^{-1} \Ub_{ni}}
\nonumber 
\\[2mm]
&\!\!\!=\!\!\!&
-
\frac{1}{n}
\sum_{i=1}^n 
\frac{\Ub_{ni}\pr (\tilde{\Vb}_{0n}^{-1}-\Vb_{0n}^{-1}) \Ub_{ni}}{(\Ub_{ni}\pr \tilde{\Vb}_{0n}^{-1} \Ub_{ni})(\Ub_{ni}\pr \Vb_{0n}^{-1} \Ub_{ni})^2}
\,
{\rm vec} (\tilde{\Vb}_{0n}^{-1/2} \Ub_{ni} \Ub_{ni}\pr\tilde{\Vb}_{0n}^{-1/2}) 
({\rm vec}(\Ub_{ni} \Ub_{ni}\pr))\pr 
.
\nonumber 
\end{eqnarray*}
The squared Frobenius norm~$\|\Mb_{2n}\|^2_F={\rm tr}[\Mb_{2n}\Mb_{2n}']$ is 
\begin{eqnarray*}
\|\Mb_{2n}\|^2_F
&\!\!\!=\!\!\!&
\frac{1}{n^2}
\sum_{i,j=1}^n 
\frac{(\Ub_{ni}\pr (\tilde{\Vb}_{0n}^{-1}-\Vb_{0n}^{-1}) \Ub_{ni})(\Ub_{nj}\pr (\tilde{\Vb}_{0n}^{-1}-\Vb_{0n}^{-1}) \Ub_{nj})}{(\Ub_{ni}\pr \tilde{\Vb}_{0n}^{-1} \Ub_{ni})(\Ub_{ni}\pr \Vb_{0n}^{-1} \Ub_{ni})^2(\Ub_{nj}\pr \tilde{\Vb}_{0n}^{-1} \Ub_{nj})(\Ub_{nj}\pr \Vb_{0n}^{-1} \Ub_{nj})^2}
\nonumber 
\\[2mm]
& &
\hspace{23mm} 
\times
(\Ub_{ni}\pr \Ub_{nj})^2
(\Ub_{ni}\pr \tilde{\Vb}_{0n}^{-1} \Ub_{nj})^2
,
\nonumber 
\end{eqnarray*}
so that, denoting as~$\rho(\Ab)$ the spectral radius of~$\Ab$, the Cauchy-Schwarz inequality yields
\begin{eqnarray*}
\lefteqn{
\hspace{-6mm} 
\|\Mb_{2n}\|^2_F
\leq
\frac{1}{n^2}
\sum_{i,j=1}^n 
\frac{|\Ub_{ni}\pr (\tilde{\Vb}_{0n}^{-1}-\Vb_{0n}^{-1}) \Ub_{ni}||\Ub_{nj}\pr (\tilde{\Vb}_{0n}^{-1}-\Vb_{0n}^{-1}) \Ub_{nj}|}{(\Ub_{ni}\pr \Vb_{0n}^{-1} \Ub_{ni})^2(\Ub_{nj}\pr \Vb_{0n}^{-1} \Ub_{nj})^2}
}
\nonumber 
\\[2mm]
& &
\hspace{18mm} 
\leq
\frac{\sup_{\ub\in\mathcal{S}^{p-1}} (\ub\pr (\tilde{\Vb}_{0n}^{-1}-\Vb_{0n}^{-1}) \ub)^2}{\inf_{\ub\in\mathcal{S}^{p-1}} (\ub\pr \Vb_{0n}^{-1} \ub)^4}
=
{\lambda}_{n1}^4
\big(\rho\big(\tilde{\Vb}_{0n}^{-1}-\Vb_{0n}^{-1}\big)\big)^2
,
\nonumber 
\end{eqnarray*}
which is~$o_{\rm P}(1)$ (note indeed that Lemma~\ref{LemB1}(ii) implies that~$\sqrt{n}\,(\tilde{\Vb}_{0n} - \Vb_{0n})$, hence also~$\sqrt{n}\,(\tilde{\Vb}_{0n}^{-1}- \Vb_{0n}^{-1})$, is~$O_{\rm P}(1)$).  Consequently, we have proved that 
\begin{eqnarray} 
\lefteqn{
\frac{1}{\sqrt{n}}
\sum_{i=1}^n
{\rm vec} 
\bigg( 
\frac{\tilde{\Vb}_{0n}^{-1/2}\Ub_{ni} \Ub_{ni}\pr\tilde{\Vb}_{0n}^{-1/2}}{\Ub_{ni}\pr \tilde{\Vb}_{0n}^{-1} \Ub_{ni}}
-
 \frac{\tilde{\Vb}_{0n}^{-1/2}\Ub_{ni} \Ub_{ni}\pr\tilde{\Vb}_{0n}^{-1/2}}{\Ub_{ni}\pr  \Vb_{0n}^{-1} \Ub_{ni}} 
\bigg) 
}
 \nonumber 
 \\[2mm]
& & 
\hspace{-6mm} 
=
-
\big( (\tilde{\Vb}_{0n}^{-1/2}\Vb_{0n}^{1/2}) \otimes  (\tilde{\Vb}_{0n}^{-1/2}\Vb_{0n}^{1/2})\big)  
\Tb_{n}(\Vb_{0n}) 
\big( \Vb_{0n}^{1/2} \otimes  \Vb_{0n}^{1/2}\big)
\sqrt{n}\, {\rm vec}\big(\tilde{\Vb}_{0n}^{-1}- \Vb_{0n}^{-1}\big)
+ o_{\rm P}(1)  
.
\nonumber
\end{eqnarray}
Still using the fact that~$\sqrt{n}\,(\tilde{\Vb}_{0n}^{-1}- \Vb_{0n}^{-1})$ is~$O_{\rm P}(1)$, we then obtain from Lemma~\ref{LemA3}(i) and the continuous mapping theorem that
\begin{eqnarray} 
\lefteqn{
\hspace{3mm} 
\frac{1}{\sqrt{n}}
\sum_{i=1}^n
{\rm vec} 
\bigg( 
\frac{\tilde{\Vb}_{0n}^{-1/2}\Ub_{ni} \Ub_{ni}\pr\tilde{\Vb}_{0n}^{-1/2}}{\Ub_{ni}\pr \tilde{\Vb}_{0n}^{-1} \Ub_{ni}}
-
 \frac{\tilde{\Vb}_{0n}^{-1/2}\Ub_{ni} \Ub_{ni}\pr\tilde{\Vb}_{0n}^{-1/2}}{\Ub_{ni}\pr  \Vb_{0n}^{-1} \Ub_{ni}} 
\bigg) 
}
\nonumber 
\\[3mm]
& & 
\hspace{4mm} 
= 
-\frac{1}{p(p+2)} 
\big( (\tilde{\Vb}_{0n}^{-1/2}\Vb_{0n}^{1/2}) \otimes  (\tilde{\Vb}_{0n}^{-1/2}\Vb_{0n}^{1/2})\big)  
({\bf I}_p + {\bf K}_p+{\bf J}_p) 
\label{step1alinD}
\\[3mm]
& & 
\hspace{40mm} 
\times 
\big( \Vb_{0n}^{1/2} \otimes  \Vb_{0n}^{1/2}\big)
\sqrt{n}\,  {\rm vec}\big(\tilde{\Vb}_{0n}^{-1}- \Vb_{0n}^{-1}\big) 
+ o_{\rm P}(1)  
\nonumber 
.
\end{eqnarray}
Now, since~$\tilde{\Vb}_{0n}^{-1/2}-{\Vb}_{0n}^{-1/2}$ is~$O_{\rm P}(1/\sqrt{n})$, Lemma~\ref{LemA3}(ii) provides
\begin{eqnarray} 
\lefteqn{
\hspace{-3mm} 
\frac{1}{\sqrt{n}}
\sum_{i=1}^n
{\rm vec} 
\bigg( 
  \frac{\tilde{\Vb}_{0n}^{-1/2}\Ub_{ni} \Ub_{ni}\pr\tilde{\Vb}_{0n}^{-1/2}}{\Ub_{ni}\pr {\Vb}_{0n}^{-1} \Ub_{ni}}
  -
   \frac{{\Vb}_{0n}^{-1/2}\Ub_{ni} \Ub_{ni}\pr {\Vb}_{0n}^{-1/2}}{\Ub_{ni}\pr  \Vb_{0n}^{-1} \Ub_{ni}} 
\bigg) 
}
\nonumber 
\\[2mm]
& &  
\hspace{-0mm} 
= 
\sqrt{n}
\Big[
\big( \tilde{\Vb}_{0n}^{-1/2}\Vb_{0n}^{1/2}\otimes \tilde{\Vb}_{0n}^{-1/2}\Vb_{0n}^{1/2} \big)
-
\mathbf{I}_{p^2}
\Big]
{\rm vec} 
\big( 
\Sb(\Vb_{0n})
\big) 
\nonumber 
\\[2mm]
& &  
\hspace{-0mm} 
= 
\frac{\sqrt{n}}{p}
\Big[
\big( \tilde{\Vb}_{0n}^{-1/2}\Vb_{0n}^{1/2}\otimes \tilde{\Vb}_{0n}^{-1/2}\Vb_{0n}^{1/2} \big)
-
\mathbf{I}_{p^2}
\Big]
{\rm vec}(\mathbf{I}_p) 
+ 
o_{\rm P}(1) 
.
\label{step2alinD}
\end{eqnarray}
Therefore, \eqref{step1alinD}--\eqref{step2alinD} provide
\begin{eqnarray} 
\lefteqn{
\hspace{-0mm} 
\sqrt{n}\, {\rm vec}\big(\Sb_n(\tilde{\Vb}_{0n})-\Sb_n({\Vb}_{0n})\big) 
= 
\frac{1}{\sqrt{n}}
\sum_{i=1}^n
{\rm vec} 
\bigg( 
\frac{\tilde{\Vb}_{0n}^{-1/2}\Ub_{ni} \Ub_{ni}\pr\tilde{\Vb}_{0n}^{-1/2}}{\Ub_{ni}\pr \tilde{\Vb}_{0n}^{-1} \Ub_{ni}}
-
 \frac{{\Vb}_{0n}^{-1/2}\Ub_{ni} \Ub_{ni}\pr {\Vb}_{0n}^{-1}}{\Ub_{ni}\pr  \Vb_{0n}^{-1} \Ub_{ni}} 
 \bigg) 
}
\nonumber 
\\[2mm]
& &  
\hspace{-6mm} 
=
 -
 \frac{1}{p(p+2)} 
\big( \tilde{\Vb}_{0n}^{-1/2}\Vb_{0n}^{1/2}\otimes \tilde{\Vb}_{0n}^{-1/2}\Vb_{0n}^{1/2} \big)
({\bf I}_p + {\bf K}_p+{\bf J}_p) 
\big( \Vb_{0n}^{1/2} \otimes  \Vb_{0n}^{1/2}\big)
\sqrt{n}\,  {\rm vec}\big(\tilde{\Vb}_{0n}^{-1}- \Vb_{0n}^{-1}\big) 
\nonumber
 \\[2mm]
& &  
\hspace{5mm} 
+
\frac{\sqrt{n}}{p}
\Big[
\big( \tilde{\Vb}_{0n}^{-1/2}\Vb_{0n}^{1/2}\otimes \tilde{\Vb}_{0n}^{-1/2}\Vb_{0n}^{1/2} \big)
-
\mathbf{I}_{p^2}
\Big]
 {\rm vec}(\mathbf{I}_p) 
+ 
o_{\rm P}(1) 
.
\nonumber
\end{eqnarray}
Hence, using the fact that~$\thetab_0$ is an  eigenvector of any matrix power of~$\tilde{\Vb}_{0n}$ and~${\Vb}_{0n}$ along with the identity~$({\bf I}_p- \thetab_0\thetab_0\pr)\thetab_0={\bf 0}$, we then obtain
$$
\big(\thetab_0\pr \otimes ({\bf I}_p- \thetab_0\thetab_0\pr)\big)
\sqrt{n} \, 
{\rm vec}\big(\Sb_n(\tilde{\Vb}_{0n})-\Sb_n({\Vb}_{0n})\big)
=
o_{\rm P}(1)
,
$$
which finally proves that
$$
T_n(\tilde{\Vb}_{0n}) - T_n(\Vb_{0n})
=
n  
p(p+2)
\big\| 
\big(\thetab_0\pr \otimes ({\bf I}_p- \thetab_0\thetab_0\pr)\big)
{\rm vec}\big(\Sb_n(\tilde{\Vb}_{0n})-\Sb_n({\Vb}_{0n})\big)
\big\|^2
.
$$
is~$o_{\rm P}(1)$.
\cqfd 
\vspace{4mm}


The proof of Theorem~\ref{TheorequivTyler} still requires the following result.

\begin{Lem} 
\label{LemB2}
Let $(\Vb_{0n})$ be a sequence of null shape matrices as in~(\ref{toestimnuis}) and assume that there exists~$\eta>0$ such that both  corresponding leading eigenvalues satisfy~$\lambda_{n1}/\lambda_{n2}\geq 1+\eta$ for~$n$ large enough. Then, as~$n\to\infty$ under~${\rm P}\n_{\Vb_{0n}}$, 
	$$
	\frac{1}{\hat{\lambda}_{n1}}
	\sum_{j=2}^p \hat{\lambda}_{nj}^{-1} {\hat\thetab}_{nj}{\hat\thetab}_{nj}\pr
=
	\frac{1}{\lambda_{n1}}
	\sum_{j=2}^p  {\lambda}_{nj}^{-1} {\thetab}_{nj} {\thetab}_{nj}\pr
+
o_{\rm P}(1)
,	
$$
where the $\hat{\lambda}_{nj}$\!'s and~$ {\hat\thetab}_{nj}$\!'s refer to the spectral decomposition of Tyler's M-estimator of scatter as~$\hat{\Vb}_n$ in~(\ref{DecompTyler}).
\end{Lem}

\noindent {\bf Proof of Lemma~\ref{LemB2}}. 
It follows from Lemma~\ref{LemB1}(ii) that~$\hat{\Vb}_n^{-1}- \Vb_{0n}^{-1}$ is~$o_{\rm P}(1)$ as~$n\to\infty$ under~${\rm P}\n_{\Vb_{0n}}$. Since~$\lambda_{n1}/\lambda_{n2}$ stays away from one, we also have that~${\hat\thetab}_{n1}-\thetab_{n1}$ is~$o_{\rm P}(1)$ under the same sequence  of hypotheses. Consequently, 
$$
(\mathbf{I}_p-\hat{\thetab}_{n1}\hat{\thetab}_{n1}\pr) \hat{\Vb}_n^{-1} - (\mathbf{I}_p-{\thetab}_{0}{\thetab}_{0}\pr) \Vb_{0n}^{-1} 
=
 \sum_{j=2}^p \hat{\lambda}_{nj}^{-1} {\hat\thetab}_{nj}{\hat\thetab}_{nj}\pr
 -
 \sum_{j=2}^p  {\lambda}_{nj}^{-1} {\thetab}_{nj} {\thetab}_{nj}\pr
$$ 
is~$o_{\rm P}(1)$ as~$n\to\infty$ under~${\rm P}\n_{\Vb_{0n}}$. The result then follows from the fact that Lemma~\ref{LemB1}(ii) also implies that~$\hat{\lambda}_{n1}-\lambda_{n1}$ is~$o_{\rm P}(1)$ as~$n\to\infty$ under~${\rm P}\n_{\Vb_{0n}}$. 
\cqfd 
\vspace{3mm}


\noindent {\bf Proof of Theorem~\ref{TheorequivTyler}}. Since~$\tilde{\Vb}_{0n}\thetab_0=\hat{\lambda}_{n1}\thetab_0$, we have that
\begin{eqnarray*} 
L_n 
&\!\!\!=\!\!\!& 
\frac{np}{p+2} 
\,
\big(
\hat{\lambda}_{n1} \thetab_0\pr \hat{\Vb}_n^{-1}\thetab_0+ \hat{\lambda}_{n1}^{-1} \thetab_0\pr  \hat{\Vb}_n \thetab_0-2
\big) 
\nonumber 
\\[2mm]
&\!\!\!=\!\!\!& 
\frac{np}{p+2} 
\sum_{j=1}^p 
\Big(
\hat{\lambda}_{n1} \hat{\lambda}_{nj}^{-1} 
(\thetab_0\pr {\hat\thetab}_{nj})^2
+ 
\hat{\lambda}_{n1}^{-1} \hat{\lambda}_{nj} 
(\thetab_0\pr {\hat\thetab}_{nj})^2
-
2(\thetab_0\pr {\hat\thetab}_{nj})^2  
\Big) 
\nonumber 
\\
&\!\!\!=\!\!\!& 
\frac{np}{(p+2) \hat{\lambda}_{n1}} 
\sum_{j=2}^p 
\,
\hat{\lambda}_{nj}^{-1} 
\big\{
(\hat{\lambda}_{nj}-\hat{\lambda}_{n1}) (\thetab_0\pr {\hat\thetab}_{nj})
\big\}^2 
\nonumber 
\\[1mm]
&\!\!\!=\!\!\!& 
\frac{np}{(p+2) \hat{\lambda}_{n1}} 
\sum_{j=2}^p
\,
 \hat{\lambda}_{nj}^{-1} 
 \big\{
 {\hat\thetab}_{nj}\pr ( \hat{\Vb}_n- \tilde{\Vb}_{0n}) \thetab_0
 \big\}
 ^2
.
\end{eqnarray*}
 Since Lemma~\ref{LemB1}(ii) entails that~$\sqrt{n}(\hat{\Vb}_n- \tilde{\Vb}_{0n})$ is $O_{\rm P}(1)$ and since both~$\tilde{\Vb}_{0n}\thetab_0=\hat{\lambda}_{n1}\thetab_0$ and ${\Vb}_{0n}\thetab_0={\lambda}_{n1}\thetab_0$ are orthogonal to~$\thetab_{nj}$, $j=2,\ldots,n$, Lemma~\ref{LemB2} then provides
\begin{eqnarray*}
L_n
&\!\!\!=\!\!\!&
\frac{np}{(p+2) {\lambda}_{n1}}  \sum_{j=2}^p {\lambda}_{nj}^{-1} ({\thetab}_{nj}\pr(\hat{\Vb}_n- \tilde{\Vb}_{0n})\thetab_0)^2 
+
o_{\rm P}(1)
\\[1mm]
&\!\!\!=\!\!\!&
\frac{np}{(p+2) {\lambda}_{n1}}  \sum_{j=2}^p {\lambda}_{nj}^{-1} ({\thetab}_{nj}\pr(\hat{\Vb}_n- \Vb_{0n})\thetab_0)^2 
+
o_{\rm P}(1)
,
\end{eqnarray*}
which, letting~$\Gb_p:={\bf I}_{p^2}-(1/(p+2)) ({\bf I}_{p^2}+{\bf K}_p-{\bf J}_p)$, rewrites
\begin{eqnarray*}
L_n
&\!\!\! = \!\!\!&
\frac{p+2}{p}  \sum_{j=2}^p 
\bigg\{
\frac{p}{(p+2)\sqrt{{\lambda}_{n1}{\lambda}_{nj}}}
(\thetab_0\pr \otimes {\thetab}_{nj}\pr)
\sqrt{n}\,{\rm vec}(\hat{\Vb}_n - \Vb_{0n})
\bigg\}^2 
+
o_{\rm P}(1)
\\[2mm]
&\!\!\! = \!\!\!&
\frac{p+2}{p}  \sum_{j=2}^p 
\Big\{
(\thetab_0\pr \otimes {\thetab}_{nj}\pr)
\Gb_p
\big(\Vb_{0n}^{-1/2} \otimes \Vb_{0n}^{-1/2}\big)
\sqrt{n}\,
{\rm vec}(\hat{\Vb}_n - \Vb_{0n}) 
\Big\}^2 
+
o_{\rm P}(1)
.
\end{eqnarray*}
By using Lemma~\ref{LemB1}(i), we therefore conclude that
\begin{eqnarray*}
L_n
&\!\!\!=\!\!\!&
{p(p+2)}\sum_{j=2}^p 
\Big\{
(\thetab_0\pr \otimes {\thetab}_{nj}\pr)\sqrt{n}\,{\rm vec}\big(\Sb_n(\Vb_{0n})-{\textstyle{\frac{1}{p}}}{\bf I}_p\big)
\Big\}^2
+
o_{\rm P}(1)
\\[0mm]
&\!\!\!=\!\!\!&
n{p(p+2)}
\big( {\rm vec}\big(\Sb_n(\Vb_{0n})-{\textstyle{\frac{1}{p}}}{\bf I}_p\big) \big)'
\big( (\thetab_0\thetab_0\pr) \otimes (\mathbf{I}_p-\thetab_0\thetab_0\pr) \big)  
{\rm vec}\big(\Sb_n(\Vb_{0n})-{\textstyle{\frac{1}{p}}}{\bf I}_p\big)
+
o_{\rm P}(1)
\\[3mm]
&\!\!\!=\!\!\!&
T_n(\Vb_{0n})+o_{\rm P}(1)
,
\end{eqnarray*}
which, in view of Theorem~\ref{Theorasymplin}, establishes the result.
\cqfd \vspace{4mm}


\section{Proofs of Theorems~\ref{nullTn}--\ref{Theor3rdLeCam}}

\noindent {\bf Proof of Theorem~\ref{nullTn}}. 
It directly follows from Lemma \ref{LemA3} that,  under ${\rm P}_{\Vb_{0n}}\n$, 
\begin{eqnarray*}
	\Wb_{n}
&\!\!\!:=\!\!\!&
\big(\thetab_0\pr \otimes ({\bf I}_p- \thetab_0\thetab_0\pr)\big)
\sqrt{n} \, 
{\rm vec}
\big( \Sb_n(\Vb_{0n})\big)
\\[1mm]
&\!\!\!=\!\!\!&
\big(\thetab_0\pr \otimes ({\bf I}_p- \thetab_0\thetab_0\pr)\big)
\sqrt{n} \, 
{\rm vec}
\big( \Sb_n(\Vb_{0n}) - {\textstyle{\frac{1}{p}}} {\bf I}_p\big)
\end{eqnarray*}
is asymptotically normal with mean zero and covariance matrix ${\bf I}_{p}- \thetab_0\thetab_0\pr$. Since~${\bf I}_{p}- \thetab_0\thetab_0\pr$ is idempotent with rank~$p-1$, this implies that, under the same sequence of hypotheses, $\Tb_n(\Vb_{0n})=\Wb_{n}\pr ({\bf I}_p- \thetab_0 \thetab_0\pr)\Wb_{n}$ is asymptotically chi-square with~$p-1$ degrees of freedom. The result then follows from Theorem~\ref{Theorasymplin}. 
\cqfd 
\vspace{4mm}


\noindent {\bf Proof of Theorem~\ref{Theor3rdLeCam}}.  
The results in~(i)--(ii) follow from a routine application of the Le Cam third lemma. For~(iii), the mutual contiguity between ${\rm P}_{\Vb_{0n}}\n$ and ${\rm P}_{\Vb_{1n}}\n$---which follows by applying the Le Cam first lemma to Theorem~\ref{TheorLAN}(iii)---enables the use of the same Le Cam third lemma. Using the notation introduced in the proof of Theorem~\ref{nullTn}, the central limit theorem yields that, under ${\rm P}_{\Vb_{0n}}\n$,  
$$
\left(
\begin{array}{c}  
\Wb_{n}
\\[1mm]
\taub_n \pr \Umb_{\thetab_0}\n \thetab_0+ \frac{1}{2\xi} \taub_n\pr \Umb_{\thetab_0}\n \taub_n
\end{array} 
\right)
%
$$ is asymptotically 
 Gaussian with mean zero and covariance matrix
$$ 
\left(
\begin{array}{cc}  
{\bf I}_p- \thetab_0 \thetab_0\pr
&
{\textstyle{\frac{\sqrt{p}}{\sqrt{p+2}}}} (1- \frac{1}{2\xi^2} \| \taub \|^2 ) ({\bf I}_p- \thetab_0 \thetab_0\pr) \taub 
\\[1mm]
{\textstyle{\frac{\sqrt{p}}{\sqrt{p+2}}}}  (1- \frac{1}{2\xi^2} \| \taub \|^2 ) \taub\pr({\bf I}_p- \thetab_0 \thetab_0\pr)
 &
{\textstyle{\frac{p}{p+2}}} \|\taub_n\|^2 \big( 1 - \frac{1}{4\xi^2} \| \taub_n\|^2\big) 
\end{array} 
\right)
.
$$ 
The Le Cam third lemma therefore ensures that, under the sequence of local alternatives considered in Part~(iii) of the theorem, $\Wb_n$ is asymptotically normal with mean~${\textstyle{\frac{\sqrt{p}}{\sqrt{p+2}}}} (1- \frac{1}{2\xi^2} \| \taub \|^2 ) ({\bf I}_p- \thetab_0 \thetab_0\pr) \taub$ and covariance matrix~${\bf I}_p- \thetab_0 \thetab_0\pr$. It follows that~$\Tb_n(\Vb_{0n})=\Wb_{n}\pr ({\bf I}_p- \thetab_0 \thetab_0\pr)\Wb_{n}$ is asymptotically chi-square with~$p-1$ degrees of freedom and non-centrality parameter
$$
\frac{p}{p+2}\,
 \| \taub\|^2 \Big(1- \frac{1}{2\xi^2} \| \taub \|^2 \Big)^2\Big(1- \frac{1}{4\xi^2} \| \taub\|^2\Big)
.
$$ 
From contiguity, Theorem~\ref{Theorasymplin} implies that the same holds for~$T_n(\tilde{\Vb}_{0n})$, which establishes the result.  Finally, Part~(iv) of the result directly follows from Theorem~\ref{TheorLAN}(iv). 
\cqfd 
\vspace{4mm}


\section*{Acknowledgement}

Davy Paindaveine's research is supported by a research fellowship from the Francqui Foundation and by the Program of Concerted Research Actions (ARC) of the Universit\'{e} libre de Bruxelles. Thomas Verdebout's research is supported by the Cr\'{e}dit de Recherche~J.0134.18 of the FNRS (Fonds National pour la Recherche Scientifique), Communaut\'{e} Fran\c{c}aise de Belgique, and by the aforementioned ARC program of the Universit\'{e} libre de Bruxelles.


\bibliographystyle{imsart-nameyear.bst}
\bibliography{Paper.bib}           

\providecommand{\noopsort}[1]{}
\begin{thebibliography}{39}

\bibitem[\protect\citeauthoryear{Anderson}{1963}]{And63}
\begin{barticle}[author]
\bauthor{\bsnm{Anderson},~\bfnm{Theodore~Wilbur}\binits{T.~W.}}
(\byear{1963}).
\btitle{Asymptotic theory for principal component analysis}.
\bjournal{Ann. Math. Statist.}
\bvolume{34}
\bpages{122--148}.
\end{barticle}
\endbibitem

\bibitem[\protect\citeauthoryear{Antoine and Lavergne}{2014}]{BeLa14}
\begin{barticle}[author]
\bauthor{\bsnm{Antoine},~\bfnm{Bertille}\binits{B.}} \AND
  \bauthor{\bsnm{Lavergne},~\bfnm{Pascal}\binits{P.}}
(\byear{2014}).
\btitle{Conditional moment models under semi-strong identification}.
\bjournal{J. Econometrics}
\bvolume{182}
\bpages{59--69}.
\end{barticle}
\endbibitem

\bibitem[\protect\citeauthoryear{Dufour}{2006}]{Duf2006}
\begin{barticle}[author]
\bauthor{\bsnm{Dufour},~\bfnm{J.~M.}\binits{J.~M.}}
(\byear{2006}).
\btitle{Monte Carlo tests with nuisance parameters: a general approach to
  finite-sample inference and nonstandard asymptotics}.
\bjournal{J. Econometrics}
\bvolume{133}
\bpages{443--477}.
\end{barticle}
\endbibitem

\bibitem[\protect\citeauthoryear{D{\"u}mbgen}{1998}]{Dum98}
\begin{barticle}[author]
\bauthor{\bsnm{D{\"u}mbgen},~\bfnm{Lutz}\binits{L.}}
(\byear{1998}).
\btitle{On Tyler's M-functional of scatter in high dimension}.
\bjournal{Ann. Inst. Statist. Math.}
\bvolume{50}
\bpages{471--491}.
\end{barticle}
\endbibitem

\bibitem[\protect\citeauthoryear{D{\"u}rre, Fried and Vogel}{2017}]{DuFri16}
\begin{barticle}[author]
\bauthor{\bsnm{D{\"u}rre},~\bfnm{Alexander}\binits{A.}},
  \bauthor{\bsnm{Fried},~\bfnm{Roland}\binits{R.}} \AND
  \bauthor{\bsnm{Vogel},~\bfnm{Daniel}\binits{D.}}
(\byear{2017}).
\btitle{The spatial sign covariance matrix and its application for robust
  correlation estimation}.
\bjournal{Austrian J. Statist.}
\bvolume{46}
\bpages{13--22}.
\end{barticle}
\endbibitem

\bibitem[\protect\citeauthoryear{D{\"u}rre, Tyler and Vogel}{2016}]{DuTyl16}
\begin{barticle}[author]
\bauthor{\bsnm{D{\"u}rre},~\bfnm{Alexander}\binits{A.}},
  \bauthor{\bsnm{Tyler},~\bfnm{David~E}\binits{D.~E.}} \AND
  \bauthor{\bsnm{Vogel},~\bfnm{Daniel}\binits{D.}}
(\byear{2016}).
\btitle{On the eigenvalues of the spatial sign covariance matrix in more than
  two dimensions}.
\bjournal{Statist. Probab. Lett.}
\bvolume{111}
\bpages{80--85}.
\end{barticle}
\endbibitem

\bibitem[\protect\citeauthoryear{D{\"u}rre, Vogel and Fried}{2015}]{DuFri15}
\begin{barticle}[author]
\bauthor{\bsnm{D{\"u}rre},~\bfnm{Alexander}\binits{A.}},
  \bauthor{\bsnm{Vogel},~\bfnm{Daniel}\binits{D.}} \AND
  \bauthor{\bsnm{Fried},~\bfnm{Roland}\binits{R.}}
(\byear{2015}).
\btitle{Spatial sign correlation}.
\bjournal{J. Multivariate Anal.}
\bvolume{135}
\bpages{89--105}.
\end{barticle}
\endbibitem

\bibitem[\protect\citeauthoryear{Flury}{1988}]{Flu88}
\begin{bbook}[author]
\bauthor{\bsnm{Flury},~\bfnm{Bernhard}\binits{B.}}
(\byear{1988}).
\btitle{Common Principal Components \& Related Multivariate Models}.
\bpublisher{John Wiley \& Sons, Inc.}, \baddress{New York}.
\end{bbook}
\endbibitem

\bibitem[\protect\citeauthoryear{Forchini and Hillier}{2003}]{ForHil2003}
\begin{barticle}[author]
\bauthor{\bsnm{Forchini},~\bfnm{G.}\binits{G.}} \AND
  \bauthor{\bsnm{Hillier},~\bfnm{G.}\binits{G.}}
(\byear{2003}).
\btitle{Conditional inference for possibly unidentified structural equations}.
\bjournal{Econometric Theory}
\bvolume{19}
\bpages{707--743}.
\end{barticle}
\endbibitem

\bibitem[\protect\citeauthoryear{Hallin and Paindaveine}{2002}]{HP02}
\begin{barticle}[author]
\bauthor{\bsnm{Hallin},~\bfnm{Marc}\binits{M.}} \AND
  \bauthor{\bsnm{Paindaveine},~\bfnm{Davy}\binits{D.}}
(\byear{2002}).
\btitle{Optimal tests for multivariate location based on interdirections and
  pseudo-Mahalanobis ranks}.
\bjournal{Ann. Statist.}
\bvolume{30}
\bpages{1103--1133}.
\end{barticle}
\endbibitem

\bibitem[\protect\citeauthoryear{Hallin and Paindaveine}{2006}]{HP06}
\begin{barticle}[author]
\bauthor{\bsnm{Hallin},~\bfnm{Marc}\binits{M.}} \AND
  \bauthor{\bsnm{Paindaveine},~\bfnm{Davy}\binits{D.}}
(\byear{2006}).
\btitle{Semiparametrically efficient rank-based inference for shape. I. Optimal
  rank-based tests for sphericity}.
\bjournal{Ann. Statist.}
\bvolume{34}
\bpages{2707--2756}.
\end{barticle}
\endbibitem

\bibitem[\protect\citeauthoryear{Hallin, Paindaveine and
  Verdebout}{2010}]{HPV10}
\begin{barticle}[author]
\bauthor{\bsnm{Hallin},~\bfnm{Marc}\binits{M.}},
  \bauthor{\bsnm{Paindaveine},~\bfnm{Davy}\binits{D.}} \AND
  \bauthor{\bsnm{Verdebout},~\bfnm{Thomas}\binits{T.}}
(\byear{2010}).
\btitle{Optimal rank-based testing for principal components}.
\bjournal{Ann. Statist.}
\bvolume{38}
\bpages{3245--3299}.
\end{barticle}
\endbibitem

\bibitem[\protect\citeauthoryear{Hallin et~al.}{2013}]{HPV13}
\begin{barticle}[author]
\bauthor{\bsnm{Hallin},~\bfnm{Marc}\binits{M.}},
  \bauthor{\bsnm{Paindaveine},~\bfnm{Davy}\binits{D.}},
  \bauthor{\bsnm{Verdebout},~\bfnm{Thomas}\binits{T.}} \betal{et~al.}
(\byear{2013}).
\btitle{Optimal rank-based tests for common principal components}.
\bjournal{Bernoulli}
\bvolume{19}
\bpages{2524--2556}.
\end{barticle}
\endbibitem

\bibitem[\protect\citeauthoryear{Hettmansperger and Randles}{2002}]{HetRan2002}
\begin{barticle}[author]
\bauthor{\bsnm{Hettmansperger},~\bfnm{TP}\binits{T.}} \AND
  \bauthor{\bsnm{Randles},~\bfnm{RH}\binits{R.}}
(\byear{2002}).
\btitle{{A practical affine equivariant multivariate median}}.
\bjournal{Biometrika}
\bvolume{89}
\bpages{851}.
\end{barticle}
\endbibitem

\bibitem[\protect\citeauthoryear{Jolicoeur}{1984}]{Jo84}
\begin{barticle}[author]
\bauthor{\bsnm{Jolicoeur},~\bfnm{Pierre}\binits{P.}}
(\byear{1984}).
\btitle{Principal components, factor analysis, and multivariate allometry: a
  small-sample direction test}.
\bjournal{Biometrics}
\bvolume{40}
\bpages{685--690}.
\end{barticle}
\endbibitem

\bibitem[\protect\citeauthoryear{Ley and Verdebout}{2017}]{LV17}
\begin{bbook}[author]
\bauthor{\bsnm{Ley},~\bfnm{Christophe}\binits{C.}} \AND
  \bauthor{\bsnm{Verdebout},~\bfnm{Thomas}\binits{T.}}
(\byear{2017}).
\btitle{Modern Directional Statistics}.
\bpublisher{CRC Press}.
\end{bbook}
\endbibitem

\bibitem[\protect\citeauthoryear{Liu and Shao}{2003}]{LiSha03}
\begin{barticle}[author]
\bauthor{\bsnm{Liu},~\bfnm{Xin}\binits{X.}} \AND
  \bauthor{\bsnm{Shao},~\bfnm{Yongzhao}\binits{Y.}}
(\byear{2003}).
\btitle{Asymptotics for likelihood ratio tests under loss of identifiability}.
\bjournal{Ann. Statist.}
\bvolume{31}
\bpages{807--832}.
\end{barticle}
\endbibitem

\bibitem[\protect\citeauthoryear{Mardia and Jupp}{2000}]{MarJup2000}
\begin{bbook}[author]
\bauthor{\bsnm{Mardia},~\bfnm{Kanti~V.}\binits{K.~V.}} \AND
  \bauthor{\bsnm{Jupp},~\bfnm{Peter~E.}\binits{P.~E.}}
(\byear{2000}).
\btitle{Directional Statistics}.
\bpublisher{John Wiley \& Sons}.
\end{bbook}
\endbibitem

\bibitem[\protect\citeauthoryear{M{\"o}tt{\"o}nen and Oja}{1995}]{mooj95}
\begin{barticle}[author]
\bauthor{\bsnm{M{\"o}tt{\"o}nen},~\bfnm{Jyrki}\binits{J.}} \AND
  \bauthor{\bsnm{Oja},~\bfnm{Hannu}\binits{H.}}
(\byear{1995}).
\btitle{Multivariate spatial sign and rank methods}.
\bjournal{J. Nonparametr. Stat.}
\bvolume{5}
\bpages{201--213}.
\end{barticle}
\endbibitem

\bibitem[\protect\citeauthoryear{Murray, Browne and McNicholas}{2016}]{uskew}
\begin{bunpublished}[author]
\bauthor{\bsnm{Murray},~\bfnm{Paula~M.}\binits{P.~M.}},
  \bauthor{\bsnm{Browne},~\bfnm{Ryan~P.}\binits{R.~P.}} \AND
  \bauthor{\bsnm{McNicholas},~\bfnm{Paul~D.}\binits{P.~D.}}
(\byear{2016}).
\btitle{uskewFactors: model-based clustering via mixtures of unrestricted
  skew-t sactor analyzer models}.
\bnote{R package.
  \url{https://cran.r-project.org/web/packages/uskewFactors/index.html}}.
\end{bunpublished}
\endbibitem

\bibitem[\protect\citeauthoryear{Oja}{2010}]{oja2010book}
\begin{bbook}[author]
\bauthor{\bsnm{Oja},~\bfnm{Hannu}\binits{H.}}
(\byear{2010}).
\btitle{Multivariate Nonparametric Methods with R. An Approach Based on Spatial
  Signs and Ranks}.
\bpublisher{Springer Science \& Business Media}.
\end{bbook}
\endbibitem

\bibitem[\protect\citeauthoryear{Paindaveine}{2006}]{Pai2006}
\begin{barticle}[author]
\bauthor{\bsnm{Paindaveine},~\bfnm{Davy}\binits{D.}}
(\byear{2006}).
\btitle{A Chernoff--Savage result for shape. On the non-admissibility of
  pseudo-Gaussian methods}.
\bjournal{J. Multivariate Anal.}
\bvolume{97}
\bpages{2206--2220}.
\end{barticle}
\endbibitem

\bibitem[\protect\citeauthoryear{Paindaveine}{2009}]{Pai2009}
\begin{barticle}[author]
\bauthor{\bsnm{Paindaveine},~\bfnm{Davy}\binits{D.}}
(\byear{2009}).
\btitle{On multivariate runs tests for randomness}.
\bjournal{J. Amer. Statist. Assoc.}
\bvolume{104}
\bpages{1525-1538}.
\end{barticle}
\endbibitem

\bibitem[\protect\citeauthoryear{Paindaveine, Remy and Verdebout}{2018}]{PRV18}
\begin{barticle}[author]
\bauthor{\bsnm{Paindaveine},~\bfnm{Davy}\binits{D.}},
  \bauthor{\bsnm{Remy},~\bfnm{Julien}\binits{J.}} \AND
  \bauthor{\bsnm{Verdebout},~\bfnm{Thomas}\binits{T.}}
(\byear{2018}).
\btitle{Testing for Principal Component Directions under Weak Identifiability}.
\bjournal{arXiv preprint arXiv:1710.05291}.
\end{barticle}
\endbibitem

\bibitem[\protect\citeauthoryear{Paindaveine and Verdebout}{2016}]{PaiVer2016}
\begin{barticle}[author]
\bauthor{\bsnm{Paindaveine},~\bfnm{Davy}\binits{D.}} \AND
  \bauthor{\bsnm{Verdebout},~\bfnm{Thomas}\binits{T.}}
(\byear{2016}).
\btitle{On high-dimensional sign tests}.
\bjournal{Bernoulli}
\bvolume{22}
\bpages{1745--1769}.
\end{barticle}
\endbibitem

\bibitem[\protect\citeauthoryear{P\"{o}tscher}{2002}]{Pot2002}
\begin{barticle}[author]
\bauthor{\bsnm{P\"{o}tscher},~\bfnm{B.~M.}\binits{B.~M.}}
(\byear{2002}).
\btitle{Lower risk bounds and properties of confidence sets for ill-posed
  estimation problems with applications to spectral density and persistence
  estimation, unit roots, and estimation of long memory parameters}.
\bjournal{Econometrica}
\bvolume{70}
\bpages{1035--1065}.
\end{barticle}
\endbibitem

\bibitem[\protect\citeauthoryear{Randles}{1989}]{Ran89}
\begin{barticle}[author]
\bauthor{\bsnm{Randles},~\bfnm{Ronald~H}\binits{R.~H.}}
(\byear{1989}).
\btitle{A distribution-free multivariate sign test based on interdirections}.
\bjournal{J. Amer. Statist. Assoc.}
\bvolume{84}
\bpages{1045--1050}.
\end{barticle}
\endbibitem

\bibitem[\protect\citeauthoryear{Randles}{2000}]{Ran00}
\begin{barticle}[author]
\bauthor{\bsnm{Randles},~\bfnm{Ronald~H}\binits{R.~H.}}
(\byear{2000}).
\btitle{A simpler, affine-invariant, multivariate, distribution-free sign
  test}.
\bjournal{J. Amer. Statist. Assoc.}
\bvolume{95}
\bpages{1263--1268}.
\end{barticle}
\endbibitem

\bibitem[\protect\citeauthoryear{Salibi\'{a}n-Barrera, Van~Aelst and
  Willems}{2006}]{Saletal2006}
\begin{barticle}[author]
\bauthor{\bsnm{Salibi\'{a}n-Barrera},~\bfnm{Mat\'{\i}as}\binits{M.}},
  \bauthor{\bsnm{Van~Aelst},~\bfnm{Stefan}\binits{S.}} \AND
  \bauthor{\bsnm{Willems},~\bfnm{Gert}\binits{G.}}
(\byear{2006}).
\btitle{Principal components analysis based on multivariate MM estimators with
  fast and robust bootstrap}.
\bjournal{J. Amer. Statist. Assoc.}
\bvolume{101}
\bpages{1198--1211}.
\end{barticle}
\endbibitem

\bibitem[\protect\citeauthoryear{Schwartzman, Mascarenhas and
  Taylor}{2008}]{Sc08}
\begin{barticle}[author]
\bauthor{\bsnm{Schwartzman},~\bfnm{Armin}\binits{A.}},
  \bauthor{\bsnm{Mascarenhas},~\bfnm{Walter~F.}\binits{W.~F.}} \AND
  \bauthor{\bsnm{Taylor},~\bfnm{Jonathan~E.}\binits{J.~E.}}
(\byear{2008}).
\btitle{Inference for eigenvalues and eigenvectors of Gaussian symmetric
  matrices}.
\bjournal{Ann. Statist.}
\bvolume{36}
\bpages{2886--2919}.
\end{barticle}
\endbibitem

\bibitem[\protect\citeauthoryear{Taskinen, Kankainen and
  Oja}{2003}]{Tasetal2003}
\begin{barticle}[author]
\bauthor{\bsnm{Taskinen},~\bfnm{S.}\binits{S.}},
  \bauthor{\bsnm{Kankainen},~\bfnm{A.}\binits{A.}} \AND
  \bauthor{\bsnm{Oja},~\bfnm{Hannu}\binits{H.}}
(\byear{2003}).
\btitle{Sign test of independence between two random vectors}.
\bjournal{Statist. Probab. Lett.}
\bvolume{62}
\bpages{9--21}.
\end{barticle}
\endbibitem

\bibitem[\protect\citeauthoryear{Taskinen, Koch and Oja}{2012}]{Taski12}
\begin{barticle}[author]
\bauthor{\bsnm{Taskinen},~\bfnm{Sara}\binits{S.}},
  \bauthor{\bsnm{Koch},~\bfnm{Inge}\binits{I.}} \AND
  \bauthor{\bsnm{Oja},~\bfnm{Hannu}\binits{H.}}
(\byear{2012}).
\btitle{Robustifying principal component analysis with spatial sign vectors}.
\bjournal{Statist. Probab. Lett.}
\bvolume{82}
\bpages{765--774}.
\end{barticle}
\endbibitem

\bibitem[\protect\citeauthoryear{Taskinen, Oja and Randles}{2005}]{Tasetal2005}
\begin{barticle}[author]
\bauthor{\bsnm{Taskinen},~\bfnm{Sara}\binits{S.}},
  \bauthor{\bsnm{Oja},~\bfnm{Hannu}\binits{H.}} \AND
  \bauthor{\bsnm{Randles},~\bfnm{Ronald~H}\binits{R.~H.}}
(\byear{2005}).
\btitle{Multivariate nonparametric tests of independence}.
\bjournal{J. Amer. Statist. Assoc.}
\bvolume{100}
\bpages{916--925}.
\end{barticle}
\endbibitem

\bibitem[\protect\citeauthoryear{Tyler}{1981}]{Tyl81}
\begin{barticle}[author]
\bauthor{\bsnm{Tyler},~\bfnm{David~E}\binits{D.~E.}}
(\byear{1981}).
\btitle{Asymptotic inference for eigenvectors}.
\bjournal{Ann. Statist.}
\bvolume{9}
\bpages{725--736}.
\end{barticle}
\endbibitem

\bibitem[\protect\citeauthoryear{Tyler}{1983a}]{Tyl83}
\begin{barticle}[author]
\bauthor{\bsnm{Tyler},~\bfnm{David~E}\binits{D.~E.}}
(\byear{1983}a).
\btitle{A class of asymptotic tests for principal component vectors}.
\bjournal{Ann. Statist.}
\bvolume{11}
\bpages{1243--1250}.
\end{barticle}
\endbibitem

\bibitem[\protect\citeauthoryear{Tyler}{1983b}]{Tyl83mu}
\begin{barticle}[author]
\bauthor{\bsnm{Tyler},~\bfnm{D.~E.}\binits{D.~E.}}
(\byear{1983}b).
\btitle{The asymptotic distribution of principal component roots under local
  alternatives to multiple roots}.
\bjournal{Ann. Statist.}
\bvolume{11}
\bpages{1232--1242}.
\end{barticle}
\endbibitem

\bibitem[\protect\citeauthoryear{Tyler}{1987a}]{Tyl1987}
\begin{barticle}[author]
\bauthor{\bsnm{Tyler},~\bfnm{D.~E.}\binits{D.~E.}}
(\byear{1987}a).
\btitle{A distribution-free M-estimator of multivariate scatter}.
\bjournal{Ann. Statist.}
\bvolume{15}
\bpages{234--251}.
\end{barticle}
\endbibitem

\bibitem[\protect\citeauthoryear{Tyler}{1987b}]{Tyl87ang}
\begin{barticle}[author]
\bauthor{\bsnm{Tyler},~\bfnm{David~E}\binits{D.~E.}}
(\byear{1987}b).
\btitle{Statistical analysis for the angular central Gaussian distribution on
  the sphere}.
\bjournal{Biometrika}
\bvolume{74}
\bpages{579--589}.
\end{barticle}
\endbibitem

\bibitem[\protect\citeauthoryear{Zhu and Zhang}{2006}]{ZhZh06}
\begin{barticle}[author]
\bauthor{\bsnm{Zhu},~\bfnm{Hongtu}\binits{H.}} \AND
  \bauthor{\bsnm{Zhang},~\bfnm{Heping}\binits{H.}}
(\byear{2006}).
\btitle{Asymptotics for estimation and testing procedures under loss of
  identifiability}.
\bjournal{J. Multivariate Anal.}
\bvolume{97}
\bpages{19--45}.
\end{barticle}
\endbibitem

\end{thebibliography}
\vspace{3mm}

\end{document}